\documentclass{article}
\usepackage{amsmath,amscd,amssymb}

\newcommand{\nbigb}{\mathcal{B}}
\newcommand{\nbigc}{\mathcal{C}}
\newcommand{\nbigd}{\mathcal{D}}
\newcommand{\nbige}{\mathcal{E}}
\newcommand{\nbigf}{\mathcal{F}}

\newcommand{\nbigh}{\mathcal{H}}

\newcommand{\nbigk}{\mathcal{K}}
\newcommand{\nbigl}{\mathcal{L}}
\newcommand{\nbigm}{\mathcal{M}}
\newcommand{\nbign}{\mathcal{N}}
\newcommand{\nbigo}{\mathcal{O}}
\newcommand{\nbigp}{\mathcal{P}}
\newcommand{\nbigq}{\mathcal{Q}}
\newcommand{\nbigr}{\mathcal{R}}
\newcommand{\nbigs}{\mathcal{S}}
\newcommand{\nbigt}{\mathcal{T}}

\newcommand{\nbigw}{\mathcal{W}}
\newcommand{\nbigx}{\mathcal{X}}
\newcommand{\nbigy}{\mathcal{Y}}

\newcommand{\proj}{\mathbb{P}}
\newcommand{\seisuu}{{\mathbb Z}}
\newcommand{\rnum}{{\mathbb Q}}

\newcommand{\cnum}{{\mathbb C}}
\newcommand{\real}{{\mathbb R}}
\newcommand{\hyperh}{\mathbb{H}}
\newcommand{\hyperr}{\mathbb{R}}

\newcommand{\newTate}{\pmb{T}}

\newcommand{\DD}{\mathbb{D}}
\newcommand{\EE}{\mathbb{E}}

\newcommand{\gbigc}{\mathfrak C}
\newcommand{\gbigd}{\mathfrak D}

\newcommand{\gbigh}{\mathfrak H}

\newcommand{\gbigk}{\mathfrak K}

\newcommand{\gbigm}{\mathfrak M}
\newcommand{\gbign}{\mathfrak N}

\newcommand{\gbigt}{\mathfrak T}

\newcommand{\gbigx}{\mathfrak X}

\newcommand{\gminib}{\mathfrak b}

\newcommand{\gminie}{\mathfrak e}

\newcommand{\vecdelta}{{\boldsymbol \delta}}

\newcommand{\veck}{{\boldsymbol k}}
\newcommand{\vecm}{{\boldsymbol m}}

\newcommand{\vecx}{{\boldsymbol x}}

\newcommand{\vecn}{{\boldsymbol n}}
\newcommand{\vecp}{{\boldsymbol p}}

\newcommand{\vecS}{{\boldsymbol S}}

\newcommand{\lrarr}{\longrightarrow}

\newcommand{\pf}{{\bf Proof}\hspace{.1in}}
\newcommand{\qed}{\hfill\mbox{\rule{1.2mm}{3mm}}}

\def\Cok{\mathop{\rm Cok}\nolimits}

\def\Image{\mathop{\rm Im}\nolimits}

\def\Gr{\mathop{\rm Gr}\nolimits}

\def\Gr{\mathop{\rm Gr}\nolimits}

\def\Res{\mathop{\rm Res}\nolimits}

\def\id{\mathop{\rm id}\nolimits}

\def\codim{\mathop{\rm codim}\nolimits}

\def\reduced{\mathop{\rm red}\nolimits}

\newcommand{\del}{\partial}
\newcommand{\delbar}{\overline{\del}}

\newcommand{\barlambda}{\overline{\lambda}}
\newcommand{\lambdabar}{\barlambda}

\newcommand{\eigenmap}{\gminie}

\newcommand{\lefttop}[1]{{}^{#1}\!}

\newcommand{\Vzero}{V^{(\lambda_0)}}

\newcommand{\tildepsi}{\widetilde{\psi}}
\newcommand{\psitilde}{\tildepsi}

\newcommand{\deldel}{\eth}

\newcommand{\distribution}{\gbigd\gminib}

\newcommand{\closedopen}[2]{[#1,#2[}

\newcommand{\ptilde}{\widetilde{p}}
\newcommand{\pitilde}{\widetilde{\pi}}

\newcommand{\nbigltilde}{\widetilde{\nbigl}}

\newcommand{\nbigrtilde}{\widetilde{\nbigr}}

\newcommand{\nbigmtilde}{\widetilde{\nbigm}}

\def\DR{\mathop{\rm DR}\nolimits}

\def\rel{\mathop{\rm rel}\nolimits}

\newcommand{\Ptilde}{\widetilde{P}}

\newcommand{\betabar}{\overline{\beta}}

\newcommand{\nbigktilde}{\widetilde{\nbigk}}

\newcommand{\Omegatilde}{\widetilde{\Omega}}

\newcommand{\vecq}{\boldsymbol q}

\newcommand{\nbigptilde}{\widetilde{\nbigp}}
\newtheorem{theorem}{Theorem}[section]
\newtheorem{corollary}[theorem]{Corollary}

\newtheorem{remark}[theorem]{Remark}
\newtheorem{lemma}[theorem]{Lemma}
\newtheorem{proposition}[theorem]{Proposition}
\newtheorem{definition}[theorem]{Definition}

\begin{document}

\title{A twistor approach to 
the Kontsevich complexes
}
\author{Takuro Mochizuki}

\date{}
\maketitle

\begin{abstract}
We study the $V$-filtration
of the mixed twistor $\nbigd$-modules
associated to algebraic meromorphic functions.
We prove that their relative de Rham complexes
are quasi-isomorphic to the family of Kontsevich complexes.
It reveals a generalized Hodge theoretic meaning of
Kontsevich complexes.
On the basis of the quasi-isomorphism,
we revisit the results on the Kontsevich complexes
due to H. Esnault, M. Kontsevich, C. Sabbah, M. Saito and J.-D. Yu
from a viewpoint of mixed twistor $\nbigd$-modules.

\vspace{.1in}
\noindent
Keywords:
Mixed twistor $\nbigd$-module,
Kontsevich complex,
$V$-filtration.

\vspace{.1in}
\noindent
14F10,
32C38,
32S35.
\end{abstract}

\section{Introduction}
\label{section;14.5.1.2}

\subsection{Mixed twistor $\nbigd$-modules}

The theory of mixed twistor $\nbigd$-modules
was developed by C. Sabbah and the author
(\cite{mochi2}, \cite{Mochizuki-wild}, \cite{Mochizuki-MTM},
\cite{sabbah2}, \cite{sabbah5}).
Very roughly, mixed twistor $\nbigd$-modules
are holonomic $\nbigd$-modules
equipped with mixed twistor structure.
We have the standard $6$-operations
on the derived category of 
algebraic mixed twistor $\nbigd$-modules
on complex algebraic manifolds,
which are compatible
with the standard $6$-operations
for algebraic holonomic $\nbigd$-modules.
For any algebraic function $f$ on 
a complex algebraic manifold $Y$,
the associated algebraic flat bundle
$(\nbigo_Y,d+df)$ is naturally enhanced to
an algebraic pure twistor $\nbigd$-module on $Y$.
As a result, we can say that
many holonomic $\nbigd$-modules
are naturally enhanced to mixed twistor $\nbigd$-modules.

One of general issues is
to describe such mixed twistor $\nbigd$-modules
as explicitly as possible.
Once we have an explicit description of 
a mixed twistor $\nbigd$-module,
we might have a chance to relate it
with a more concrete object,
and to apply a general theory of mixed twistor $\nbigd$-modules
for the study of the object.

\paragraph{$\nbigr$-modules}
Let us recall the concept of $\nbigr$-modules
which is one of the ingredients to formulate mixed twistor
$\nbigd$-modules.
Let $\cnum_{\lambda}$ denote just a complex line
with the coordinate $\lambda$.
For any complex manifold $Y$,
let $\nbigr_Y$ denote the sheaf of algebras
on $\cnum_{\lambda}\times Y$
obtained as the subalgebra of
the sheaf of holomorphic differential operators
$\nbigd_{\cnum_{\lambda}\times Y}$
generated by
$\lambda p_{\lambda}^{\ast}\Theta_Y$.
Here, $p_{\lambda}:\cnum_{\lambda}\times Y\lrarr Y$
denotes the projection,
and $\Theta_Y$ denotes the tangent sheaf of $Y$.
When we are given a local coordinate system $(x_1,\ldots,x_n)$,
then $\lambda\del_{x_i}$ are denoted by
$\deldel_{x_i}$ or $\deldel_i$.

Mixed twistor $\nbigd$-modules $\nbigt$ on $Y$
are formulated as a pair of $\nbigr_Y$-modules
$\nbigm_i$ $(i=1,2)$,
a sesqui-linear pairing $C$
of $\nbigm_1$ and $\nbigm_2$,
and a weight filtration $W$,
satisfying some conditions.
(See \cite{sabbah2} and \cite{Mochizuki-MTM}
for more details on sesqui-linear pairings,
weight filtrations, and the conditions.)
In this paper,
$\nbigm_2$ is called the underlying $\nbigr_Y$-module
of $\nbigt$.
Note that 
$\Xi_{\DR}(\nbigm_2):=
 \iota_1^{-1}\bigl(\nbigm_2/(\lambda-1)\nbigm_2\bigr)$
is naturally a $\nbigd_Y$-module,
where $\iota_1:\{1\}\times Y\lrarr \cnum\times Y$
is the inclusion.
We call 
$\Xi_{\DR}(\nbigm_2)$
the $\nbigd$-module underlying $\nbigt$.

Then, we reword the general issue as follows:
We would like to describe the $\nbigr$-modules
underlying mixed twistor $\nbigd$-modules
as explicitly as possible.

\subsection{Main result}
\label{subsection;17.1.20.20}

In this paper,
we study the $\nbigr$-modules
underlying the mixed twistor $\nbigd$-module
associated to algebraic meromorphic functions.
More precisely, 
let $X$ be a smooth complex projective manifold
with a morphism $f:X\lrarr \proj^1$.
Let $D$ be a hypersurface of $X$
such  that $f^{-1}(\infty)\subset D$.
We assume that $D$ is normal crossing.
We put $X^{(1)}:=\cnum_{\tau}\times X$
and $D^{(1)}:=\cnum_{\tau}\times D$.
We obtain the meromorphic function
$\tau f$ on $(X^{(1)},D^{(1)})$.
We have the holonomic $\nbigd_{X^{(1)}}$-module
$\nbigm:=\nbigo_{X^{(1)}}(\ast D^{(1)})\,v$
with the flat connection $\nabla$
given by 
$\nabla v=v\,d(\tau f)$.
It is naturally enhanced to a mixed twistor $\nbigd$-module
$\nbigt_{\ast}(\tau f,D^{(1)})$.
We have the underlying $\nbigr$-module
$\nbigl_{\ast}(\tau f,D^{(1)})$.
We obtain
the $\nbigr_{X}(\ast\tau)$-module
$\nbigmtilde:=\nbigl_{\ast}(\tau f,D^{(1)})(\ast\tau)$.

We consider the sheaf of subalgebras 
$\lefttop{\tau}V_0\nbigr_{X^{(1)}}$
in $\nbigr_{X^{(1)}}$
generated by
$\nbigo_{\cnum_{\lambda}\times X^{(1)}}$
and 
$\lambda p_{\lambda}^{\ast}\Theta_{X^{(1)}}(\log \tau)$,
where $\Theta_{X^{(1)}}(\log\tau)$ denote
the sheaf of vector fields which are logarithmic
along $\tau=0$.
By a general theory of mixed twistor $\nbigd$-modules,
$\nbigmtilde$
is uniquely equipped with an increasing sequence of
coherent $\lefttop{\tau}V_0\nbigr_{X^{(1)}}$-submodules
$U_{\alpha}\nbigmtilde$ $(\alpha\in\real)$
with the following property:
\begin{itemize}
\item
 $\bigcup U_{\alpha}\nbigmtilde=\nbigmtilde$.
\item
 For any $\alpha\in\real$,
 we have $\epsilon>0$
 such that
 $U_{\alpha}\nbigmtilde=U_{\alpha+\epsilon}\nbigmtilde$.
\item
We have
 $\tau U_{\alpha}\nbigmtilde
 =U_{\alpha-1}\nbigmtilde$
and 
 $\deldel_{\tau} U_{\alpha}\nbigmtilde
 \subset U_{\alpha+1}\nbigmtilde$ for any $\alpha\in\real.$
\item
 The induced endomorphisms
 $\tau\deldel_{\tau}+\lambda\alpha$
 on 
 $\Gr^{U}_{\alpha}\nbigmtilde:=
 U_{\alpha}\nbigmtilde/U_{<\alpha}\nbigmtilde$
 are nilpotent 
 for any $\alpha\in\real$.
Here, we set
 $U_{<\alpha}\nbigmtilde:=
 \bigcup_{\beta<\alpha}U_{\beta}\nbigmtilde$.
\item
 $\Gr^{U}_{\alpha}\nbigmtilde$
 are strict,
 i.e., flat over $\nbigo_{\cnum_{\lambda}}$.
\end{itemize}
Note that the $V$-filtrations of $\nbigr$-modules
underlying mixed twistor $\nbigd$-modules
are characterized by a more involved condition.
Because $\nbigmtilde$ is naturally equipped with
the action of $\lambda^2\del_{\lambda}$,
the condition is simplified as above.

We shall describe $U_{\alpha}\nbigmtilde$ explicitly.
Then, we shall show that
their de Rham complex relative to 
$X^{(1)}\lrarr \cnum_{\tau}$
is quasi-isomorphic to
an explicitly given family of complexes,
called Kontsevich complexes.

\paragraph{Kontsevich complexes}
Let us recall the concept of 
Kontsevich complexes.
Put $P:=f^{\ast}(\infty)$ as an effective divisor.
The reduced divisor is denoted by $P_{\reduced}$.
The multiplication of $df$
induces a morphism
$\Omega^k_X(\log D)
\lrarr
 \Omega^{k+1}_X(\log D)\otimes\nbigo_X(P)$.
The inverse image of
$\Omega^{k+1}_X(\log D)
 \subset
 \Omega^{k+1}_X(\log D)\otimes\nbigo_X(P)$
by $df$
is denoted by $\Omega_f^k$.
The multiplication of $df$
induces a morphism
$df:\Omega_f^k\lrarr \Omega_f^{k+1}$.
The exterior derivative induces
$d:\Omega_f^k\lrarr \Omega_f^{k+1}$.

For any $0\leq\alpha<1$,
we set 
$\Omega^{k}_f(\alpha):=
 \Omega^{k}_f\otimes
 \nbigo_X([\alpha P])$.
Here, 
for real numbers $a_i$ 
and reduced hypersurfaces $H_i$ 
$(i=1,\ldots,N)$,
we set $[a_i]:=\max\{n\in\seisuu\,|\,n\leq a_i\}$
and 
$\bigl[\sum_{i=1}^{N} a_i H_i\bigr]
:=\sum_{i=1}^N [a_i]\,H_i$.
We have the induced operators
$d:\Omega_f^k(\alpha)\lrarr \Omega_f^{k+1}(\alpha)$
and 
$df:\Omega_f^k(\alpha)\lrarr \Omega_f^{k+1}(\alpha)$.

For any $(\lambda,\tau)\in\cnum^2$,
we have the derivative
$\lambda d+\tau df:
 \Omega_f^k(\alpha)
\lrarr
 \Omega_f^{k+1}(\alpha)$.
They satisfy the integrability condition
$(\lambda d+\tau df)\circ(\lambda d+\tau df)=0$.
Thus, we obtain complexes
$\bigl(
 \Omega_f^{\bullet}(\alpha),
 \lambda d+\tau df
 \bigr)$.
They are called Kontsevich complexes.

We clearly have the family version of the complexes.
Let $q_X:\cnum_{\lambda}\times X^{(1)}\lrarr X$
denote the projection.
We set
$\Omegatilde_f^{j}:=
 \lambda^{-j}q_X^{\ast}\Omega_f^j$.
Then, we have the derivative
$d+\lambda^{-1}\tau\,df:
 \Omegatilde_f^j\lrarr
 \Omegatilde_f^{j+1}$
satisfying 
$(d+\lambda^{-1}\tau\,df)\circ(d+\lambda^{-1}\tau\,df)=0$.

\begin{theorem}[Theorem 
\ref{thm;14.5.2.20}]
\label{thm;17.1.20.10}
For $0\leq\alpha<1$,
$\bigl(
 \Omegatilde_f^{\bullet}(\alpha),d+\lambda^{-1}\tau\,df
 \bigr)$
is naturally quasi-isomorphic to
the de Rham complex of
$U_{\alpha}\nbigmtilde$
relative to the projection $X^{(1)}\lrarr \cnum_{\tau}$.
\end{theorem}

This theorem reveals a generalized Hodge theoretic property
of the Kontsevich complexes,
and provides a bridge between
the theory of mixed twistor $\nbigd$-modules
and the study of Kontsevich complexes.
We shall show that it
is useful by giving an alternative proof
of the interesting results of 
H. Esnault, M. Kontsevich, C. Sabbah, M. Saito
and J.-D. Yu,
which we will explain in the next subsection.

\subsection{An application}
\label{subsection;17.7.30.101}

The following theorem was conjectured
by M. Kontsevich.
It was proved by
H. Esnault, C. Sabbah and J.-D. Yu
\cite{Esnault-Sabbah-Yu}.
Some interesting cases were proved by
Kontsevich and M. Saito.
(See \cite{Katzarkov-Kontsevich-Pantev-2017}
and the appendix of \cite{Esnault-Sabbah-Yu}.)

\begin{theorem}
\label{thm;14.4.24.20}
The dimension of the hypercohomology groups
$\hyperh^i\bigl(
 X,
 \Omega_f^{\bullet}(\alpha),
 \lambda d+\tau df
 \bigr)$
are independent of $(\lambda,\tau)\in\cnum^2$
and any $0\leq \alpha<1$.
\end{theorem}

We set $\gbigx:=\proj^1_{\tau}\times X$.
Let $p_i$ denote the projection of $\gbigx$
onto the $i$-th components.
We identify $\nbigo_{\proj^1}(1)$
with $\nbigo_{\proj^1}(\{\infty\})$.
The constant function $1$ and
the coordinate function $\tau$
are naturally regarded as sections of $\nbigo_{\proj^1}(1)$.
We consider the bundles
$\gbigk^k_f(\alpha):=
 p_1^{\ast}\nbigo_{\proj^1}(k)
\otimes
 p_2^{\ast}\Omega^k_{f}(\alpha)$.
We have the relative differential operators
\[
 p_2^{\ast}(d)+\tau p_2^{\ast}(df):
 \gbigk_f^k(\alpha)
\lrarr
 \gbigk_f^{k+1}(\alpha).
\]
Thus, we obtain a complex of sheaves
$\gbigk_f^{\bullet}(\alpha)$
on $\gbigx$.
We have the decreasing filtration $F^{\bullet}$
of the complex $\gbigk_f^{\bullet}(\alpha)$
defined as follows:
\begin{equation}
 \label{eq;17.7.30.102}
F^j\gbigk_f^i(\alpha):=
\left\{
 \begin{array}{ll}
 0 & (i<j)\\
 \gbigk_f^i(\alpha) & (i\geq j)
 \end{array}
\right.
\end{equation}

By Theorem \ref{thm;14.4.24.20},
we obtain a vector bundle
$\nbigk^i_f(\alpha):=
 \hyperr^ip_{1\ast}\bigl(
 \gbigk^{\bullet}_f(\alpha)
 \bigr)$
on $\proj^1$.
Let 
$F^j\nbigk^i_f(\alpha)$
denote the image of
$\hyperr^ip_{1\ast}\bigl(
 F^j\gbigk_f^{\bullet}(\alpha)
 \bigr)$.
We have
\[
 \Gr_F^j\nbigk^i_f(\alpha)
\simeq
 H^{i-j}\bigl(X,\Omega_f^j(\alpha)\bigr)
\otimes
 \nbigo_{\proj^1}(j).
\]
Hence, $F^{\bullet}$ 
is the Harder-Narasimhan filtration
of $\nbigk^i_f(\alpha)$.
Moreover, Sabbah and Yu proved
the following theorem in \cite{Sabbah-Yu}.

\begin{theorem}
\label{thm;14.4.24.21}
Let $0\leq\alpha<1$.
\begin{description}
\item[\rm(i)]
The bundle $\nbigk^i_f(\alpha)$
is equipped with a naturally induced
meromorphic connection
$\nabla$
such that
\begin{equation}
 \label{eq;14.5.1.1}
 \nabla\cdot \nbigk^i_f(\alpha)
\subset
 \nbigk^i_f(\alpha)
 \otimes\Omega^1_{\proj^1}\bigl(\{0\}+2\{\infty\}\bigr).
\end{equation}
Note that the condition 
{\rm(\ref{eq;14.5.1.1})} implies that
$\nabla F^j\nbigk_f^i(\alpha)\subset 
 F^{j-1}\nbigk_f^i(\alpha)
 \otimes
 \Omega^1_{\proj^1}\bigl(\{0\}+2\{\infty\}\bigr)$.
\item[\rm(ii)]
Let $\Res_0(\nabla)$ denote the endomorphism of
$\nbigk_f^i(\alpha)_{|0}$
obtained as the residue.
Then, any eigenvalue $\beta$ of
$\Res_0(\nabla)$ is a rational number
such that
$-\alpha\leq\beta<-\alpha+1$.
\end{description}
\end{theorem}

Let $\EE_{\beta}\nbigk_f^i(\alpha)_{|0}$
denote the generalized eigen space 
of $\Res_0(\nabla)$ corresponding to $\beta$.
By setting
$U_{\gamma}\nbigk_f^i(\alpha)_{|0}
=\bigoplus_{-\beta\leq \gamma}
 \EE_{\beta}\nbigk_f^i(\alpha)_{|0}$,
we obtain an increasing filtration
$U_{\bullet}\nbigk_f^i(\alpha)_{|0}$
indexed by
$\alpha-1<\gamma\leq\alpha$.
Note that 
$\Gr^U_{\gamma}\nbigk_f^i(\alpha)_{|0}$
is naturally isomorphic to
$\EE_{-\gamma}\nbigk_f^i(\alpha)_{|0}$.
The filtration $F$ on 
$\nbigk_f^i(\alpha)_{|0}$
naturally induces filtrations
on $\Gr^U_{\gamma}\nbigk_f^i(\alpha)_{|0}$,
which are also denoted by $F$.
Because $\Res_0(\nabla)$ preserves
the filtration $U$,
we have the induced endomorphisms
$\Gr^U_{\gamma}\Res_0(\nabla)$
of $\Gr^U_{\gamma}\nbigk_f^i(\alpha)_{|0}$.
Let $N_{\gamma}$ denote the nilpotent part of
$\Gr^U_{\gamma}\Res_0(\nabla)$.
Then, Sabbah and Yu also proved the following
in \cite{Sabbah-Yu}.
\begin{theorem}
\label{thm;17.7.30.20}
For any $\alpha-1<\gamma\leq\alpha$,
$N_{\gamma}$ gives a strict morphism
\[
 (\Gr^U_{\gamma}\nbigk^i_{f}(\alpha)_{|0},F)
\lrarr
 (\Gr^U_{\gamma}\nbigk^i_{f}(\alpha)_{|0},F[-1]),
\]
where $F[-1]$ denotes the filtration
defined by $F[-1]^a:=F^{a-1}$.
Namely,
we have
$N_{\gamma}\bigl(
 F^j\Gr^U_{\gamma}\nbigk_f^i(\alpha)_{|0}
 \bigr)
=\Image(N_{\gamma})
 \cap
 F^{j-1}\Gr^U_{\gamma}\nbigk_f^i(\alpha)_{|0}$.
\end{theorem}

In \S\ref{subsection;14.12.28.10},
we shall explain how to deduce 
Theorems \ref{thm;14.4.24.20}, \ref{thm;14.4.24.21}
and \ref{thm;17.7.30.20}
from Theorem \ref{thm;17.1.20.10}
and general results for mixed twistor $\nbigd$-modules.

We remark that our argument is based 
on a generalized Hodge theory
and some concrete computations of $V$-filtrations.
Hence, eventually, 
it is not completely different
from those in \cite{Esnault-Sabbah-Yu,Sabbah-Yu}.
But, the argument in this paper looks more direct.
The author thinks that 
it would be desirable to have 
many ways to explain the theorems
because the complexes
$(\Omega_f^{\bullet}(\alpha),\lambda d+\tau df)$
and the bundles $\nbigk^i_f$ seem quite basic.
He also hopes that 
this paper might also be useful to explain
how the general theory of twistor $\nbigd$-modules
could be applied to the study of specific objects.

\paragraph{Outline of the paper}

In \S\ref{section;14.5.16.30},
we study  the $\nbigd$-module
associated to the meromorphic function
as in \S\ref{subsection;17.1.20.20}.
In particular,
we explicitly describe the relative de Rham complexes
of their $V$-filtration along $\tau$
(Propositions \ref{prop;14.5.1.30}
and \ref{prop;14.4.24.10}),
which is essentially a simpler version of
Theorem \ref{thm;17.1.20.10}.
Then, we explain how to deduce
a part of Theorem \ref{thm;14.4.24.21}.
We could include the computations for $\nbigd$-modules 
in \S\ref{section;14.5.16.30}
to the computations for $\nbigr$-modules 
in \S\ref{section;14.5.13.1}
after minor modifications.
But, the author expects that it would be useful to give 
an explanation in this simpler situation.

In \S\ref{section;14.5.13.1},
we study the $\nbigr$-module
in \S\ref{subsection;17.1.20.20},
and we revisit Theorems \ref{thm;14.4.24.20},
\ref{thm;14.4.24.21}
and \ref{thm;17.7.30.20}
from the viewpoint of mixed twistor $\nbigd$-modules.
After the preliminaries in 
\S\ref{subsection;17.1.20.201}--\ref{subsection;14.5.16.31},
we explain  
in \S\ref{subsection;17.1.20.200}
the main theorem
(Theorem \ref{thm;14.5.2.20})
and how we deduce 
Theorems \ref{thm;14.4.24.20},
\ref{thm;14.4.24.21}
and \ref{thm;17.7.30.20}
from general results of mixed twistor $\nbigd$-modules.
In \S\ref{subsection;14.5.14.1},
we give a description of the $V$-filtration
of the $\nbigr_{X^{(1)}}(\ast\tau)$-module $\nbigmtilde$
(Theorem \ref{thm;14.5.2.40}).
Besides the argument in \S\ref{section;14.5.16.30}
for $\nbigd$-modules,
we need an additional task to check the strictness.
Then, we establish Theorem \ref{thm;14.5.2.20}
in \S\ref{subsection;14.5.14.20}.

\paragraph{Acknowledgements}

The author thanks the referee for 
valuable suggestions to improve this paper.
This note is written to understand
the intriguing work of H. Esnault,
C. Sabbah, M. Saito and J.-D. Yu
on the Kontsevich complexes
\cite{Esnault-Sabbah-Yu}, \cite{Sabbah-Yu}, \cite{Yu}.
I thank Sabbah for sending 
an earlier version of \cite{Sabbah-Yu}
and for discussions on many occasions.
I also thank him for clarifying 
the statement of Theorem \ref{thm;17.7.30.20}.
I am grateful to Esnault for some discussions
and for her kindness.
I thank M.-H. Saito for his kindness and support.
I thank A. Ishii and Y. Tsuchimoto 
for their constant encouragement.

This study was partially supported by 
the Grant-in-Aid for Scientific Research 
(C) (No. 22540078), 
the Grant-in-Aid for Scientific Research 
(C) (No. 15K04843), 
the Grant-in-Aid for Scientific Research
(A) (No. 22244003),
the Grant-in-Aid for Scientific Research 
(S) (No. 24224001), 
the Grant-in-Aid for Scientific Research 
(S) (No. 17H06127) and
the Grant-in-Aid for Scientific Research 
(S) (No. 16H06335),
Japan Society for the Promotion of Science.

\section{The case of $\nbigd$-modules}
\label{section;14.5.16.30}

\subsection{Meromorphic flat bundles}

We continue to use the setting in 
\S\ref{section;14.5.1.2}.
We set
$X^{(1)}:=\cnum_{\tau}\times X$.
Here, $\cnum_{\tau}$
is just an affine line with a coordinate $\tau$.
We use the notation
$D^{(1)}$, $P^{(1)}$, etc.,
with a similar meaning.

For any complex manifold $Y$ with a hypersurface $Y_1$,
let $\nbigo_Y(\ast Y_1)$ denote the sheaf of
meromorphic functions on $Y$
whose poles are contained in $Y_1$.
For any $\nbigo_Y$-module $M$,
let $M(\ast Y_1):=M\otimes_{\nbigo_Y}\nbigo_Y(\ast Y_1)$.
If $Y_1$ is given as $\{f=0\}$ for a holomorphic function $f$,
we also use the notation $\nbigo_Y(\ast f)$
and $M(\ast f)$.

We shall consider the meromorphic flat bundle
$\nbigm:=\nbigo_{X^{(1)}}(\ast D^{(1)})\,v$
with
\[
 \nabla v=v\,d(\tau f)
\]
of rank one,
where $v$ denotes a global frame.
The corresponding  $\nbigd_{X^{(1)}}$-module
is also denoted by $\nbigm$.
By the isomorphism
$\nbigo_{X^{(1)}}(\ast D^{(1)})
\simeq
 \nbigm$;
$1\longleftrightarrow v$,
the connection $\nabla$ is identified with 
the connection $d+d(\tau f)$ on
$\nbigo_{X^{(1)}}(\ast D^{(1)})$.

\subsubsection{Local coordinate systems}
\label{subsection;14.4.23.1}

When we study $\nbigm$
locally around a point $(\tau, Q)\in D_{\reduced}^{(1)}$,
we shall use a holomorphic local coordinate neighbourhood
$(U,x_1,\ldots,x_n)$ of $Q$ in $X$
with the following property:
\begin{itemize}
\item
 $P_{\reduced}\cap U=\bigcup_{i=1}^{\ell_1}\{x_i=0\}$,
 $H\cap U=\bigcup_{i=\ell_1+1}^{\ell}\{x_i=0\}$,
 $f_{|U}=\prod_{i=1}^{\ell_1}x_i^{-k_i}$
 for $k_i>0$ $(i=1,\ldots,\ell_1)$.
\end{itemize}
We set $k_i:=0$ $(i=\ell_1+1,\ldots,\ell)$,
and the tuple $(k_i\,|\,i=1,\ldots,\ell)\in\seisuu^{\ell}$
 is denoted by $\veck$.
For any $\vecm\in\seisuu^{\ell}$,
we set $x^{\vecm}:=\prod_{i=1}^{\ell}x_i^{m_i}$.
In particular, we have $f=x^{-\veck}$ on $U$.
We set $U^{(1)}:=\cnum_{\tau}\times U$
which is equipped with a coordinate system
$(\tau,x_1,\ldots,x_n)$.

We have
$\del_i v=-k_i \tau fx_i^{-1}\,v$
for $i=1,\ldots,\ell$,
and 
$\del_i v=0$ 
for $i=\ell+1,\ldots,n$.
We also have $\del_{\tau}v=f\,v$.
We set $\vecdelta:=(1,\ldots,1)\in\seisuu^{\ell}$.
We have
\[
 \tau\del_{\tau}(x^{-\vecdelta+\vecm}v)
=\tau fx^{-\vecdelta+\vecm}v,
\quad
 \del_ix_i(x^{-\vecdelta+\vecm}v)
=(m_i-k_i\tau f)x^{-\vecdelta+\vecm}v.
\]
Hence, for $i=1,\ldots,\ell_1$,
we have
\begin{equation}
 (\tau\del_{\tau}+k_i^{-1}\del_ix_i)
 (x^{-\vecdelta+\vecm}v)
=(m_i/k_i)x^{-\vecdelta+\vecm} v.
\end{equation}

When we are given a real number $\alpha$,
we set
$[\alpha k_i]:=
 \max\bigl\{n\in\seisuu\,\big|\,n\leq \alpha k_i\bigr\}$,
and 
$[\alpha\veck]:=
 \bigl(
 [\alpha k_i]\,\big|\,
 i=1,\ldots,\ell
 \bigr)\in\seisuu^{\ell}$.

\subsection{$V$-filtration along $\tau$}

Let $\pi:X^{(1)}\lrarr X$ denote the projection.
We naturally regard
$\pi^{\ast}\nbigd_{X}$
as the subalgebra of 
$\nbigd_{X^{(1)}}$.
Let $\lefttop{\tau}V_0\nbigd_{X^{(1)}}\subset \nbigd_{X^{(1)}}$
denote the sheaf of subalgebras
generated by $\tau\del_{\tau}$
over $\pi^{\ast}\nbigd_{X}$.
We shall construct 
$\lefttop{\tau}V_{0}\nbigd_{X^{(1)}}$-submodules
$U_{\alpha}\nbigm$ $(\alpha\in\real)$
of $\nbigm$,
and we shall prove that the filtration
$U_{\bullet}\nbigm$ is a $V$-filtration of $\nbigm$.

\vspace{.1in}

For $0\leq \alpha\leq 1$,
we have an $\nbigo_{X^{(1)}}$-submodule
$\nbigo_{X^{(1)}}\bigl(D^{(1)}+[\alpha P^{(1)}]\bigr)\,v
\subset\nbigm$.
We set
\[
 U_{\alpha}\nbigm:=
 \pi^{\ast}\nbigd_X\cdot
 \bigl(
 \nbigo_{X^{(1)}}(D^{(1)}+[\alpha P^{(1)}])\,v
 \bigr)
\subset\nbigm.
\]

\begin{lemma}
\label{lem;14.5.20.10}
Outside
$\{\tau=0\}\cap P^{(1)}$,
we have
$U_{\alpha}\nbigm=\nbigm$ for $0\leq\alpha\leq 1$.
\end{lemma}
\pf
The claim is clear outside $D^{(1)}$.
Let $(\tau_0,Q)$ be any point of $D^{(1)}$.
We use a coordinate system
as in \S\ref{subsection;14.4.23.1}.
Suppose $Q\not\in P_{\reduced}$.
For any $\vecm\in\seisuu^{\ell}_{\geq 0}$,
we have
$\del_i(x^{-\vecdelta-\vecm}v)
=-(m_i+1)x^{-\vecdelta-\vecm}x_i^{-1}v$.
Hence, we easily obtain 
\[
 \pi^{\ast}\nbigd_X\cdot\bigl(
 \nbigo_{X^{(1)}}(D^{(1)})
 \bigr)v=\nbigo_{X^{(1)}}(\ast D^{(1)})v
\]
on a neighbourhood of $(\tau_0,Q)$.
Suppose $Q\in P_{\reduced}$ and $\tau_0\neq 0$.
For any $\vecm\in\seisuu^{\ell}_{\geq 0}$,
we have
$\del_i(x^{-\vecdelta-\vecm}v)
=-(1+m_i+k_i\tau x^{-\veck})x^{-\vecdelta-\vecm}x_i^{-1}v$.
Then, we easily obtain
$\nbigo(\ast D^{(1)})v
 \subset
 U_{\alpha}\nbigm$.
\qed

\begin{lemma}
\label{lem;14.5.1.11}
We have
$\del_{\tau}U_0\nbigm
\subset U_1\nbigm$
and
$\tau U_1\nbigm
\subset
 U_0\nbigm$.
\end{lemma}
\pf
We have
$\del_{\tau}(gx^{-\vecdelta}v)
=(\del_{\tau}g)x^{-\vecdelta}v
+gx^{-\vecdelta-\veck}v
\in U_1\nbigm$.
Hence, we have
$\del_{\tau}U_0\nbigm
\subset U_1\nbigm$.
Let us prove
$\tau U_1\nbigm\subset U_0\nbigm$.
We have only to check it locally around 
any point of $P^{(1)}_{\reduced}$
by using  a coordinate system as in
\S\ref{subsection;14.4.23.1}.
We have
\[
\tau \bigl(
 gx^{-\vecdelta-\veck}v
 \bigr)
=gx^{-\vecdelta}(\tau f)v
=
-k_i^{-1}\del_{i}(x_i gx^{-\vecdelta}v) 
+k_i^{-1}(\del_ig)x^{-\vecdelta}x_iv
\in U_0\nbigm.
\]
Thus, we are done.
\qed

\vspace{.1in}

If $\alpha\leq 0$,
we take the integer $n$
such that $0\leq \alpha+n<1$,
and we set
$U_{\alpha}\nbigm:=
 \tau^nU_{\alpha+n}\nbigm$.
For $\alpha\geq 1$,
we define
$U_{\alpha}\nbigm:=
 \sum_{\beta+n\leq\alpha}
 \del_{\tau}^nU_{\beta}\nbigm$.
If $\alpha\leq \alpha'$,
we have
$U_{\alpha}\nbigm\subset U_{\alpha'}\nbigm$.
By the construction,
for any $\alpha\in\real$,
we have $\epsilon>0$
such that $U_{\alpha}\nbigm=U_{\alpha+\epsilon}\nbigm$.
We define
$U_{<\alpha}\nbigm:=
 \bigcup_{\beta<\alpha}
 U_{\beta}\nbigm$.

\begin{proposition}
\label{prop;14.5.1.10}
$U_{\bullet}\nbigm$
is a $V$-filtration of $\nbigm$
along $\tau$
indexed by the rational numbers
with the standard order
(up to shift of the degree by $1$).
More precisely, the following holds:
\begin{itemize}
\item
$U_{\alpha}\nbigm$
are coherent $\lefttop{\tau}V_0\nbigd_{X^{(1)}}$-modules
such that 
$\bigcup_{\alpha} U_{\alpha}\nbigm=\nbigm$.
\item
We have
$\tau U_{\alpha}\nbigm
\subset
 U_{\alpha-1}\nbigm$
and 
$\del_{\tau} U_{\alpha}\nbigm
\subset
 U_{\alpha+1}\nbigm$.
\item
$\tau\del_{\tau}+\alpha$
is nilpotent on 
$U_{\alpha}\nbigm/U_{<\alpha}\nbigm$.
Indeed,
we have
$\bigl(
 \tau\del_{\tau}+\alpha
\bigr)^{\dim X+1}=0$
on $U_{\alpha}\nbigm/U_{<\alpha}\nbigm$.
\end{itemize}
\end{proposition}
\pf
We divide the claim into several lemmas.

\begin{lemma}
\label{lem;14.5.1.3}
We have a natural action of
$\lefttop{\tau}V_0\nbigd_{X^{(1)}}$
on $U_{\alpha}\nbigm$.
\end{lemma}
\pf
It is enough to prove
$\tau\del_{\tau}\bigl(
 U_{\alpha}\nbigm
 \bigr)
\subset U_{\alpha}\nbigm$.
We have only to consider the case $0\leq\alpha\leq 1$.
We have only to check it 
locally around any point of $P^{(1)}_{\reduced}$.
We use a local coordinate system as in 
\S\ref{subsection;14.4.23.1}.
We set $\vecp:=[\alpha\veck]$.
We have
$\nbigo_{X^{(1)}}(D^{(1)}+[\alpha P^{(1)}]) v
=\nbigo_{X^{(1)}}\,x^{-\vecdelta-\vecp}v$
on $U^{(1)}$.
We have
\[
 \tau\del_{\tau}
 (x^{-\vecdelta-\vecp}v)
+(p_i/k_i)x^{-\vecdelta-\vecp}v
=-k_i^{-1}\del_i(x_ix^{-\vecdelta-\vecp}v).
\]
Then, the claim is clear.
\qed

\begin{lemma}
\label{lem;14.5.2.30}
$\lefttop{\tau}V_0\nbigd_{X^{(1)}}$-modules
$U_{\alpha}\nbigm$
are coherent.
\end{lemma}
\pf
They are pseudo-coherent over $\nbigo_{X^{(1)}}$
and locally finitely generated over $\lefttop{\tau}V_0\nbigd_{X^{(1)}}$.
Hence, they are coherent over $\lefttop{\tau}V_0\nbigd_{X^{(1)}}$.
(See \cite{kashiwara_text}.)
\qed

\vspace{.1in}
Let $[\tau\del_{\tau}]$ denote
the endomorphism 
of $U_{\alpha}\nbigm/U_{<\alpha}\nbigm$
induced by $\tau\del_{\tau}$.

\begin{lemma}
\label{lem;14.5.2.31}
For any $0<\alpha\leq 1$,
$[\tau\del_{\tau}]+\alpha$
is nilpotent.
Indeed,
$\bigl([\tau\del_{\tau}]+\alpha\bigr)^{\dim X}=0$
on $U_{\alpha}\nbigm/U_{<\alpha}\nbigm$.
\end{lemma}
\pf
Because the support of the sheaf
$U_{\alpha}\nbigm/U_{<\alpha}\nbigm$
is contained in $P^{(1)}_{\reduced}$,
we have only to check the claim locally around
any point of $P^{(1)}_{\reduced}$
by using a coordinate system as in \S\ref{subsection;14.4.23.1}.
We set $\vecp=(p_i):=[\alpha\veck]$.
If $\alpha=p_i/k_i$,
we have
$\bigl(
 \tau\del_{\tau}+\alpha
\bigr)(x^{-\vecdelta-\vecp}v)
=-k_i^{-1}\del_i(x_ix^{-\vecdelta-\vecp}v)$.
Let $S:=\{i\,|\,\alpha=p_i/k_i\}$.
Then, we have
\[
\bigl(\tau\del_{\tau}+\alpha\bigr)^{|S|}
(x^{-\vecdelta-\vecp}v)
=\prod_{i\in S}(-k_i^{-1}\del_i)
 \Bigl(
x^{-\vecdelta-\vecp}
 \prod_{i\in S}x_i\cdot
 v\Bigr)
 \in U_{<\alpha}\nbigm.
\]

We also have the following equality
for any holomorphic function $g$:
\[
 (\tau\del_{\tau}+\alpha)(gx^{-\vecdelta-\vecp}v)
=g\cdot(\tau\del_{\tau}+\alpha)(x^{-\vecdelta-\vecp}v)
+(\tau\del_{\tau}g)(x^{-\vecdelta-\vecp}v).
\]
We have the following 
for any $i$ with $k_i\neq 0$:
\begin{multline}
 k_i
 \tau(\del_{\tau}g)x^{-\vecdelta-\vecp}v
=
-\del_i\Bigl(
 x_i(\del_{\tau}g)x^{-\vecdelta-\vecp+\veck}v
 \Bigr)
+(\del_i\del_{\tau}g)x^{-\vecdelta-\vecp+\veck}x_iv
 \\
+(\del_{\tau}g)(k_i-p_i)x^{-\vecdelta-\vecp+\veck}v
\in U_{<\alpha}\nbigm.
\end{multline}
Then, we can easily deduce the claim of the lemma.
\qed

\begin{lemma}
\label{lem;14.5.2.33}
If $N\geq \dim X+1$,
we have 
$(\tau\del_{\tau})^{N}U_0\nbigm
\subset
 \tau U_{<1}\nbigm$.
\end{lemma}
\pf
We have only to check the claim locally around
any point of $P^{(1)}_{\reduced}$
by using a coordinate system as in
\S\ref{subsection;14.4.23.1}.
Set $\vecdelta_1:=
 (\overbrace{1,\ldots,1}^{\ell_1},0,\ldots,0)\in\seisuu^{\ell}$.
We have 
$(\tau\del_{\tau})^{\ell_1}(x^{-\vecdelta}v)
=\prod_{i=1}^{\ell_1}
 \bigl(-k_i^{-1}\del_i\bigr)\cdot 
 x^{-(\vecdelta-\vecdelta_1)}v$.
Hence, we have
\begin{multline}
 (\tau\del_{\tau})^{\ell_1+1}
 (x^{-\vecdelta}v)
=\prod_{i=1}^{\ell_1}(-k_i^{-1}\del_i)
 (\tau f x^{-(\vecdelta-\vecdelta_1)}v)
 \\
=\tau
 \prod_{i=1}^{\ell_1}(-k_i^{-1}\del_i)
 (x^{-\vecdelta-(\veck-\vecdelta_1)}v)
\in \tau U_{<1}\nbigm.
\end{multline}

For any section $s$ of $U_0\nbigm$
and any holomorphic function $g$,
we have
$\tau\del_{\tau}(gs)
=\tau(\del_{\tau}g)\,s
+g\,\tau\del_{\tau}s$,
and $\tau(\del_{\tau}g)s\in \tau U_{<1}(\nbigm)$.
Then, we can deduce the claim of the lemma.
\qed

\begin{lemma}
\label{lem;14.5.13.2}
We have
$\sum_{j=0}^{\infty}\del_{\tau}^jU_1\nbigm=\nbigm$.
\end{lemma}
\pf
Let $\nbigm'$ denote the left hand side.
We clearly have $\nbigm'\subset\nbigm$.
Let $\vecm\in\seisuu_{\geq 0}^{\ell}$.
For $i=\ell_1+1,\ldots,\ell$,
we have
$\del_{i}(x^{-\vecdelta-\vecm}v)
=-(m_i+1)x^{-\vecdelta-\vecm}x_i^{-1}v$.
We also have
$\del_{\tau}(x^{-\vecdelta-\vecm}v)
=x^{-\vecdelta-\vecm-\veck}v$.
Then, we obtain
$\nbigm=\nbigo_{X^{(1)}}(\ast D^{(1)})
\subset\nbigm'$.
\qed

\vspace{.1in}
We obtain Proposition \ref{prop;14.5.1.10}
from Lemmas \ref{lem;14.5.1.3}--\ref{lem;14.5.13.2}.
\qed

\begin{remark}
We set $V_{\alpha}\nbigm:=U_{\alpha+1}\nbigm$,
then $-\del_{\tau}\tau-\alpha$ is nilpotent
on $V_{\alpha}/V_{<\alpha}$.
\end{remark}

\subsection{Primitive expression}

\subsubsection{Primitive expression
for sections of $\nbigo_{X^{(1)}}(\ast D^{(1)})$}
\label{subsection;14.5.1.20}

Let $(U,x_1,\ldots,x_n)$ be a holomorphic coordinate 
neighbourhood of $X$ as in \S\ref{subsection;14.4.23.1}.
We set $U^{(1)}:=\cnum_{\tau}\times U$.
We set 
$Y:=\{\tau=x_1=\cdots =x_{\ell}=0\}\subset U^{(1)}$.
We shall consider local sections 
$\nbigo_{X^{(1)}}(\ast D^{(1)})$ 
on a small neighbourhood of $Y$.
Such a section $s$ has
the unique Laurent expansion:
\[
 s=\sum_{\vecm\in\seisuu^{\ell}}
 \sum_{j\in\seisuu_{\geq 0}}
 h_{\vecm,j}x^{\vecm}\tau^j.
\]
Here, $h_{\vecm,j}$
are holomorphic functions on $Y$.
Note that if $N>0$ is sufficiently large, depending on $s$,
we have $h_{\vecm,j}=0$ unless 
$m_i\geq -N$ $(i=1,\ldots,\ell)$.
We also have the following unique expression:
\[
  s=\sum_{\vecm\in\seisuu^{\ell}}
 \sum_{j\in\seisuu_{\geq 0}}
 h^{(1)}_{\vecm,j}x^{\vecm}\cdot(\tau f)^j.
\]
Here, $h^{(1)}_{\vecm,j}$ are holomorphic functions
on $Y$.
Indeed, we have 
$h^{(1)}_{\vecm,j}=h_{\vecm-j\veck,j}$.

\begin{definition}
A section $s$ of $\nbigo_{X^{(1)}}(\ast D^{(1)})$
on a neighbourhood of $Y$
is called $(\vecm,j)$-primitive
if we have $s=gx^{\vecm}(\tau f)^j$
for a section $g$ of $\nbigo_{X^{(1)}}$
with $g_{|Y}\neq 0$.
\end{definition}

\begin{definition}
A primitive expression of a section $s$
of $\nbigo_{X^{(1)}}(\ast D^{(1)})$ on a neighbourhood of $Y$
is a decomposition
\[
 s=\sum_{(\vecm,j)\in\nbigs}s_{\vecm,j}
\]
where $\nbigs\subset\seisuu^{\ell}\times\seisuu_{\geq 0}$
is a finite subset,
and each $s_{\vecm,j}$ is $(\vecm,j)$-primitive.
\end{definition}

\begin{lemma}
\label{lem;14.5.1.21}
Any section $s$ of $\nbigo_{X^{(1)}}(\ast D^{(1)})$
on a neighbourhood of $Y$
has a primitive expression.
\end{lemma}
\pf
We have an expression
$s=\sum_{(\vecm,j)\in\nbigt}
 g_{\vecm,j}x^{\vecm}\tau^j$
for a finite subset
$\nbigt\subset\seisuu^{\ell}\times\seisuu_{\geq 0}$
and holomorphic functions
$g_{\vecm,j}$
with $g_{\vecm,j|Y}\neq 0$.
Then, we have
$s=\sum_{(\vecm,j)\in\nbigt}
 g_{\vecm,j}x^{\vecm+j\veck}
 (\tau f)^j$.
Each $g_{\vecm,j}x^{\vecm+j\veck}(\tau f)^j$
is $(\vecm+j\veck,j)$-primitive.
\qed

\vspace{.1in}

We consider the partial order on $\seisuu^{\ell}$
defined by
$(a_i)\leq (b_i)
\stackrel{\rm def}{\Longleftrightarrow}
a_i\leq b_i\,(\forall i)$.
For any $T\subset\seisuu^{\ell}$,
let $\min(T)$ denote the set of the minimal elements
in $T$ with respect to the partial order.
We also use a similar partial order 
on $\seisuu^{\ell}\times\seisuu_{\geq 0}$.
Let $\pi:\seisuu^{\ell}\times\seisuu_{\geq 0}
 \lrarr \seisuu^{\ell}$
denote the projection.

\begin{lemma}
\label{lem;14.5.15.1}
Let $\nbigt\subset\seisuu^{\ell}\times\seisuu_{\geq 0}$
be a finite subset.
Suppose 
$\sum_{(\vecm,j)\in\nbigt}
 g_{\vecm,j}x^{\vecm}(\tau f)^j=0$
on a neighbourhood of $Y$.
\begin{itemize}
\item
For any $\vecm\in\min\pi(\nbigt)$
and any $j\in\seisuu_{\geq 0}$,
we have $g_{\vecm,j|Y}=0$.
\item
For any $(\vecm,j)\in\min(\nbigt)$,
we have $g_{\vecm,j|Y}=0$.
\end{itemize}
\end{lemma}
\pf
For any $(\vecm,j)$,
we have the Laurent expansion
\[
 g_{\vecm,j}=\sum_{\vecn\in\seisuu_{\geq 0}^{\ell}}
 \sum_{q\in\seisuu_{\geq 0}}
 g_{\vecm,j;\vecn,q}x^{\vecn}\tau^q.
\]
Here, $g_{\vecm,j;\vecn,q}$ are holomorphic functions
on $Y$.
We have
\[
 0=\sum_{(\vecm,j)\in\nbigt}
 \sum_{\vecn\in\seisuu_{\geq 0}^{\ell}}
 \sum_{q\geq 0}
 g_{\vecm,j;\vecn,q}
 x^{\vecm+\vecn+q\veck}
(\tau f)^{j+q}.
\]
If $\vecm\in\min\pi(\nbigt)$,
the coefficient of $x^{\vecm}(\tau f)^j$ is
$g_{\vecm,j;0,0}=g_{\vecm,j|Y}$.
Hence, we obtain
$g_{\vecm,j|Y}=0$.
If $(\vecm,j)\in\min(\nbigt)$,
the coefficient of
$x^{\vecm}(\tau f)^j$ is
$g_{\vecm,j;0,0}=g_{\vecm|Y}$.
Hence, we obtain
$g_{\vecm,j|Y}=0$.
\qed

\begin{corollary}
\label{cor;14.5.2.50}
Let $s$ be a non-zero section of
$\nbigo_{X^{(1)}}(\ast D^{(1)})$
on a neighbourhood of $Y$
with a primitive expression
$s=\sum_{(\vecm,j)\in\nbigs}
 g_{\vecm,j}x^{\vecm}(\tau f)^j$.
Then, the following holds:
\begin{itemize}
\item
The set $\min\pi(\nbigs)$
is well defined for $s$.
For any $\vecm\in\min\pi(\nbigs)$
and any $j\in\seisuu_{\geq 0}$,
$g_{\vecm,j|Y}$
is well defined for $s$.
\item
The set $\min(\nbigs)$ is well defined for $s$.
For any $(\vecm,j)\in\min(\nbigs)$,
$g_{\vecm,j|Y}$ is well defined for $s$.
\qed
\end{itemize}
\end{corollary}

\subsubsection{Subsheaf $\nbigm_0^{\alpha}$}

For $0\leq \alpha\leq 1$
and for $N\in\seisuu_{\geq 0}$,
we define
\[
 G_N\nbigm_0^{\alpha}:=
 \sum_{j=0}^{N}
 \nbigo_{X^{(1)}}\bigl(D^{(1)}+[\alpha P^{(1)}]\bigr)
 (\tau f)^jv
\subset\nbigm.
\]
We set 
$\nbigm_0^{\alpha}:=
 \bigcup_{N\geq 0}G_N\nbigm_0^{\alpha}$.
Let $V_0\nbigd_X$ denote the subalgebra
of $\nbigd_X$ generated by
$\nbigo_X$ and
the logarithmic tangent sheaf
$\Theta_X(\log D)$
of $X$ with respect to $D$.
We naturally regard 
$\pi^{\ast}V_0\nbigd_X$
as a subsheaf of
$\nbigd_{X^{(1)}}$.

\begin{lemma}
\label{lem;14.5.2.100}
We have
$\nbigm_0^{\alpha}
=\pi^{\ast}V_0\nbigd_X\cdot\Bigl(
 \nbigo_{X^{(1)}}\bigl(D^{(1)}+[\alpha P^{(1)}]\bigr)
 \,v
\Bigr)$.
\end{lemma}
\pf
We have only to check the claim
locally around any point of $P^{(1)}_{\reduced}$.
We use a local coordinate system as in \S\ref{subsection;14.4.23.1}.
We set $\vecp=(p_i):=[\alpha\veck]$.
Because
$\del_ix_i\bigl(
 vx^{-\vecdelta-\vecp}(\tau f)^j
 \bigr)
=-(p_i+jk_i)v\,x^{-\vecdelta-\vecp}(\tau f)^j
-k_i v x^{-\vecdelta-\vecp}
 (\tau f)^{j+1}$,
we easily obtain the claim of the lemma.
\qed

\vspace{.1in}
By the construction,
we have
$U_{\alpha}\nbigm
=\pi^{\ast}\nbigd_X\cdot\nbigm_0^{\alpha}$
for any $0\leq\alpha\leq 1$.

\subsubsection{Primitive expression
for sections of $U_{\alpha}\nbigm$
$(0\leq\alpha\leq 1)$}
\label{subsection;17.1.19.1}

We use the notation in \S\ref{subsection;14.5.1.20}.
Let $0\leq\alpha\leq 1$.
Set $\vecp:=[\alpha\veck]\in\seisuu^{\ell}$.
For any $\vecn\in\seisuu^{\ell}_{\geq 0}$,
we set $\del^{\vecn}:=\prod_{i=1}^{\ell}\del_i^{n_i}$
and $|\vecn|=\sum n_i$.
For any $\vecm\in\seisuu^{\ell}$,
we have the unique decomposition
$\vecm=\vecm_+-\vecm_-$
such that $\vecm_{\pm}\in\seisuu_{\geq 0}^{\ell}$
and 
$\{i\,|\,m_{+,i}\neq 0\}\cap
 \{i\,|\,m_{-,i}\neq 0\}=\emptyset$.

Any section $s$ of $U_{\alpha}\nbigm$
on a neighbourhood of $Y$
has an expression
\[
 s=\sum_{\vecn\in\seisuu^{\ell}_{\geq 0}}
 \del^{\vecn}s_{\vecn}
\]
as an essentially finite sum,
where $s_{\vecn}$
are sections of $\nbigm_0^{\alpha}$.
Here, ``essentially finite'' means
there exists a finite subset $T\subset\seisuu^{\ell}_{\geq 0}$
such that $s_{\vecn}=0$ unless $\vecn\in T$.

\begin{definition}
Let $(\vecm,j)\in\seisuu^{\ell}\times\seisuu_{\geq 0}$.
A section $s$ of $U_{\alpha}\nbigm$
on a neighbourhood of $Y$
is called $(\vecm,j)$-primitive if
$s=\del^{\vecm_-}(g\vecx^{-\vecp-\vecdelta+\vecm_+}(\tau f)^jv)$
for a holomorphic function $g$ with $g_{|Y}\neq 0$,
i.e.,
$g\vecx^{\vecm_+}(\tau f)^j$ is 
$(\vecm_+,j)$-primitive
as a section of
$\nbigo_{X^{(1)}}(\ast D^{(1)})$.
\qed
\end{definition}

\begin{definition}
Let $s$ be a non-zero section of $U_{\alpha}\nbigm$
on a neighbourhood of $Y$.
A primitive expression of $s$ is a decomposition
\[
 s=\sum_{(\vecm,j)\in\nbigs}s_{\vecm,j}
\]
where 
$\nbigs\subset\seisuu^{\ell}\times\seisuu_{\geq 0}$
is  a finite set,
and $s_{\vecm,j}$ are $(\vecm,j)$-primitive
sections of $\nbigm_0^{\alpha}$.
\end{definition}

This kind of expressions have been used in
the study of pure twistor $\nbigd$-modules
\cite{mochi2},
for example.

\begin{lemma}
\label{lem;14.5.2.110}
Any non-zero section $s$ of $U_{\alpha}\nbigm$
on a neighbourhood of $Y$
has a primitive expression.
\end{lemma}
\pf
We have the following expression
as an essentially finite sum:
\[
 s=\sum_{\vecn\in\seisuu^{\ell}_{\geq 0}}
 \sum_{\vecq\in\seisuu_{\geq 0}^{\ell}}
 \sum_{j\in\seisuu_{\geq 0}}
 \del^{\vecn}
 \Bigl(
 g_{\vecn,\vecq,j}
 x^{-\vecdelta-\vecp+\vecq}
 (\tau f)^jv
 \Bigr).
\]
Here, ``essentially finite''
means that there exists a finite subset
$T\subset \seisuu^{\ell}_{\geq 0}
\times\seisuu_{\geq 0}^{\ell}
\times\seisuu_{\geq 0}$
such that 
$g_{\vecn,\vecq,j}=0$
unless $(\vecn,\vecq,j)\in T$.

We consider the following claim.
\begin{description}
\item[$(P_a)$:]
If $s$ has an expression
$s=\sum_{|\vecn|\leq a}
 \sum_{\vecq,j}
 \del^{\vecn}
 \bigl(
 g_{\vecn,\vecq,j}
 x^{-\vecdelta-\vecp+\vecq}
 (\tau f)^jv
 \bigr)$
as an essentially finite sum,
then $s$ has a primitive expression
$s=\sum_{(\vecm,j)\in\nbigs}
 s_{\vecm,j}$
such that 
$|\vecm_-|\leq a$
for any $\vecm\in\pi(\nbigs)$.
\end{description}
If $a=0$, the claim is given
by Lemma \ref{lem;14.5.1.21}.
We prove $(P_a)$ by assuming $(P_{a-1})$.
Note that if $q_i>0$,
we have 
\begin{multline}
\label{eq;14.4.24.1}
 \del_ix_i\bigl(
 g
 x^{-\vecdelta-\vecp+\vecq}x_i^{-1}
 (\tau f)^jv
 \bigr)
=
 \Bigl(
 x_i\del_ig
+(-p_i+q_i-1-k_ij)g
 \Bigr)x^{-\vecdelta-\vecp+\vecq}x_i^{-1}(\tau f)^jv
\\
-k_igx^{-\vecdelta-\vecp+\vecq}x_i^{-1}
 (\tau f)^{j+1}v.
\end{multline}
By using (\ref{eq;14.4.24.1}),
$s$ has an expression
\[
 s=\sum_{\vecn,\vecq,j}
 \del^{\vecn}\bigl(
 g^{(1)}_{\vecn,\vecq,j}
 x^{-\vecdelta-\vecp+\vecq}
 (\tau f)^jv
 \bigr)
\]
with the following property:
\begin{itemize}
\item
 $g^{(1)}_{\vecn,\vecq,j}=0$
 unless $|\vecn|\leq a$.
\item
If $|\vecn|=a$ and $g^{(1)}_{\vecn,\vecq,j}\neq 0$,
we have $\{i\,|\,n_i\neq 0,q_i\neq 0\}
=\emptyset$
and $g^{(1)}_{\vecn,\vecq,j|Y}\neq 0$.
\end{itemize}
By applying $(P_{a-1})$
to the lower term
$\sum_{|\vecn|<a}\sum_{\vecq,j}
 \del^{\vecn}\bigl(
 g^{(1)}_{\vecn,\vecq,j}
 x^{-\vecdelta-\vecp+\vecq} 
 (\tau f)^jv
 \bigr)$,
we obtain a primitive expression of $s$ as desired.
\qed

\vspace{.1in}

Suppose that we are given
a finite set $\nbigs\subset\seisuu^{\ell}\times\seisuu_{\geq 0}$
and sections $g_{\vecm,j}$
of $\nbigo_{X^{(1)}}$
on a neighbourhood of $Y$
for $(\vecm,j)\in\nbigs$,
such that
\begin{equation}
 \label{eq;17.1.20.1}
 \sum_{(\vecm,j)\in\nbigs}
 \del^{\vecm_-}
 \Bigl(
 g_{\vecm,j}x^{-\vecdelta-\vecp+\vecm_+}
 (\tau f)^jv
 \Bigr)=0.
\end{equation}

\begin{lemma}
We have $g_{\vecm,j|Y}=0$
for any $\vecm\in\min\pi(\nbigs)$
and $j\in\seisuu$.
We also have
$g_{\vecm,j|Y}=0$
for any $(\vecm,j)\in\min\nbigs$.
\end{lemma}
\pf
In general, we have the following equality
for any section $g$ of $\nbigo_{X^{(1)}}$
on a neighbourhood of $Y$:
\begin{multline}
 \del_i\bigl(gx^{-\vecdelta-\vecp+\vecn}(\tau f)^jv\bigr)=
\Bigl(
x_i\del_ig
-(1+p_i-n_i+jk_i)g
\Bigr)
 x^{-\vecdelta-\vecp+\vecn}x_i^{-1}
 (\tau f)^jv
 \\
-k_igx^{-\vecdelta-\vecp+\vecn}x_i^{-1}
 (\tau f)^{j+1}v.
\end{multline}
Hence, we have the following expression:
\[
 \del^{\vecm_-}
 \bigl(g_{\vecm,j}x^{-\vecdelta-\vecp}x^{\vecm_+}
 (\tau f)^jv
 \bigr)
=\sum_{0\leq k\leq |\vecm_-|}
 h_{\vecm,j,k}x^{-\vecdelta-\vecp+\vecm}
 (\tau f)^{j+k}v,
\]
where
$h_{\vecm,j,k}$ are sections of
$\nbigo_{X^{(1)}}$ on a neighbourhood of $Y$
such that 
 $h_{\vecm,j,k|Y}=
 C_{\vecm,j,k}\cdot g_{\vecm,j|Y}$
for some $C_{\vecm,j,k}\in\rnum$.
Because
$\{i\,|\,m_{+,i}\neq 0,m_{-,i}\neq 0\}=\emptyset$,
we have $C_{\vecm,j,0}\neq 0$.

By (\ref{eq;17.1.20.1}),
we have the following 
in $\nbigo_{X^{(1)}}(\ast D^{(1)})$:
\[
\sum_{(\vecm,j)\in\nbigs}
 \sum_{0\leq k\leq |\vecm_-|}
 h_{\vecm,j,k}x^{-\vecdelta-\vecp+\vecm}
 (\tau f)^{j+k}
=0.
\]
Take $\vecm\in\min\pi(\nbigs)$.
According to Lemma \ref{lem;14.5.15.1},
for any $p\geq 0$,
we have
\[
 \sum_{j+k=p}
 C_{\vecm,j,k}g_{\vecm,j|Y}=0.
\]
We obtain
$g_{\vecm,j|Y}=0$ by an ascending induction on $j$.
Take $(\vecm,j)\in\min\nbigs$.
By Lemma \ref{lem;14.5.15.1},
we obtain 
$h_{\vecm,j,0|Y}
=C_{\vecm,j,0}g_{\vecm,j|Y}=0$.
Hence, we obtain
$g_{\vecm,j|Y}=0$.
\qed

\vspace{.1in}

\begin{corollary}
\label{cor;14.5.2.111}
Let $s$ be a section of $U_{\alpha}\nbigm$
on a neighbourhood of $Y$
with a primitive expression
\[
s=\sum_{(\vecm,j)\in\nbigs}
 \del^{\vecm_-}\bigl(
 g_{\vecm,j}x^{-\vecdelta-\vecp+\vecm_+}
 (\tau f)^jv
\bigr). 
\]
\begin{itemize}
\item
The set $\min\pi(\nbigs)$ is well defined for $s$.
For any $\vecm\in\min\pi(\nbigs)$
and $j\in\seisuu_{\geq 0}$,
$g_{\vecm,j|Y}$ is well defined for $s$.
\item
The set $\min\nbigs$ is well defined for $s$.
For any $(\vecm,j)\in\min(\nbigs)$,
$g_{\vecm,j|Y}$ is well defined for $s$.
\qed
\end{itemize}
\end{corollary}

\subsection{Quasi-isomorphism of complexes}

\subsubsection{Statements}
\label{subsection;14.5.17.3}

We set $\nbigm_0^{\alpha}(-D^{(1)}):=
 \nbigm_0^{\alpha}\otimes \nbigo_{X^{(1)}}(-D^{(1)})$.
The action of 
$\pi^{\ast}V_0\nbigd_{X}$
on $\nbigm_0^{\alpha}(-D^{(1)})$
naturally induces a complex
$\nbigm_0^{\alpha}(-D^{(1)})
\otimes
 \Omega^{\bullet}_{X^{(1)}/\cnum_{\tau}}(\log D^{(1)})$.
We have a natural inclusion of complexes:
\begin{equation}
\label{eq;14.4.23.3}
 \nbigm_0^{\alpha}(-D^{(1)})\otimes
 \Omega^{\bullet}_{X^{(1)}/\cnum_{\tau}}(\log D^{(1)})
\lrarr
 U_{\alpha}\nbigm\otimes
 \Omega_{X^{(1)}/\cnum_{\tau}}^{\bullet}.
\end{equation}
We shall prove the following proposition
in \S\ref{subsection;14.5.17.2}.
\begin{proposition}
\label{prop;14.5.1.30}
For $0\leq\alpha\leq 1$,
the morphism
{\rm(\ref{eq;14.4.23.3})}
is a quasi-isomorphism.
\end{proposition}

We set 
$\Omega_f^k(\alpha):=
 \Omega_f^k\otimes\nbigo([\alpha P])$.
We set
$\Omega_{f,\tau}^k(\alpha):=
 \pi^{\ast}\Omega_f^k(\alpha)$.
In the case $\alpha=0$,
we also use the symbols
$\Omega_f^k$
and $\Omega_{f,\tau}^k$.
We obtain a complex
$\Omega^{\bullet}_{f,\tau}(\alpha)$
with the differential given by
$d+\tau\,df$.
We have a natural morphism of complexes
induced by the correspondence $1\longmapsto v$:
\begin{equation}
\label{eq;14.4.23.4}
\Omega^{\bullet}_{f,\tau}(\alpha)
\lrarr
 \nbigm_0^{\alpha}(-D^{(1)})
 \otimes\Omega^{\bullet}_{X^{(1)}/\cnum_{\tau}}(\log D^{(1)}).
\end{equation}
We shall prove the following proposition
in \S\ref{subsection;14.5.17.1}.
\begin{proposition}
\label{prop;14.4.24.10}
For $0\leq \alpha\leq 1$,
the morphism {\rm(\ref{eq;14.4.23.4})}
is a quasi-isomorphism.
\end{proposition}

Let $p_2:X^{(1)}\lrarr\cnum_{\tau}$ denote the projection.
\begin{corollary}
\label{cor;14.5.2.2}
$\hyperr^ip_{2\ast}\bigl(
 \Omega_f^{\bullet}(\alpha)
 \bigr)$ $(0\leq \alpha<1)$
is a locally free $\nbigo_{\cnum_{\tau}}$-module
equipped with a naturally induced logarithmic connection
$\nabla$.
Any eigenvalue $\beta$ of $\Res_0(\nabla)$
is a rational number such that
$-\alpha\leq \beta< -\alpha+1$.
\end{corollary}
\pf
By the construction,
it is easy to see that
$\hyperr^ip_{2\ast}\bigl(
 \Omega_f^{\bullet}(\alpha)
 \bigr)$
is a coherent $\nbigo_{\cnum_{\tau}}$-module.
By the above propositions,
$\hyperr^ip_{2\ast}\bigl(
 \Omega_f^{\bullet}(\alpha)
 \bigr)$
is equipped with the action of 
the differential operator $\tau\del_{\tau}$.
The restriction 
$\hyperr^ip_{2\ast}\bigl(
 \Omega_f^{\bullet}(\alpha)
 \bigr)_{|\cnum_{\tau}^{\ast}}$
is a $\nbigd_{\cnum_{\tau}^{\ast}}$-module.
Because it is coherent
as an $\nbigo_{\cnum_{\tau}^{\ast}}$-module,
it is a locally free
$\nbigo_{\cnum_{\tau}^{\ast}}$-module.

Although we may finish the proof 
by applying the functoriality of the $V$-filtration,
let us recall the proof in this easy situation
with the argument in \S3.2 of \cite{saito1}.

The multiplication of $\tau^n$ $(n\geq 0)$
induces an isomorphism
of $\pi^{\ast}\nbigd_X$-modules
$U_{\alpha}\nbigm
\simeq
 U_{\alpha-n}\nbigm$
because $\alpha<1$.
We have the exact sequence of 
$\pi^{\ast}\nbigd_X$-modules:
\[
 0\lrarr
 U_{\alpha-1}\nbigm
\lrarr
 U_{\alpha}\nbigm
\lrarr
 U_{\alpha}\nbigm/U_{\alpha-1}\nbigm
\lrarr 0.
\]
It induces the following exact sequence of
coherent $\nbigo_{\cnum_{\tau}}$-modules:
\begin{multline}
 \hyperr^ip_{2\ast}
 \bigl(
 U_{\alpha-1}\nbigm
 \otimes\Omega^{\bullet}_{X^{(1)}/\cnum_{\tau}}
 \bigr)
\lrarr
 \hyperr^ip_{2\ast}
 \bigl(
 U_{\alpha}\nbigm
 \otimes\Omega^{\bullet}_{X^{(1)}/\cnum_{\tau}}
 \bigr) \\
\lrarr
 \hyperr^ip_{2\ast}
 \bigl(
 (U_{\alpha}\nbigm/U_{\alpha-1}\nbigm)
  \otimes\Omega^{\bullet}_{X^{(1)}/\cnum_{\tau}}
 \bigr)
\stackrel{\varphi}{\lrarr}
  \hyperr^{i+1}p_{2\ast}
 \bigl(
 U_{\alpha-1}\nbigm
  \otimes\Omega^{\bullet}_{X^{(1)}/\cnum_{\tau}}
 \bigr)
\end{multline}
Suppose that $\varphi\neq 0$,
and we will deduce a contradiction.
For $N\geq 1$,
the image of the morphism
\[
 \hyperr^{i+1}p_{2\ast}
 \bigl(
 U_{\alpha-N}\nbigm
 \otimes\Omega^{\bullet}_{X^{(1)}/\cnum_{\tau}}
 \bigr)
\lrarr
 \hyperr^{i+1}p_{2\ast}
 \bigl(
 U_{\alpha-1}\nbigm
 \otimes\Omega^{\bullet}_{X^{(1)}/\cnum_{\tau}}
 \bigr)
\]
is equal to the image of the multiplication of
$\tau^{N-1}$
on 
$\hyperr^{i+1}p_{2\ast}
 \bigl(
 U_{\alpha-1}\nbigm
 \otimes\Omega^{\bullet}_{X^{(1)}/\cnum_{\tau}}
 \bigr)$.
Hence, by applying Nakayama's lemma
to $\hyperr^{i+1}p_{2\ast}
 \bigl(
 U_{\alpha-1}\nbigm
 \otimes\Omega^{\bullet}_{X^{(1)}/\cnum_{\tau}}
 \bigr)$,
there exists $N_0$ such that
the image of the composite of 
the following morphisms is non-zero:
\begin{multline}
\label{eq;14.5.19.1}
  \hyperr^ip_{2\ast}
 \bigl(
 (U_{\alpha}\nbigm/U_{\alpha-1}\nbigm)
  \otimes\Omega^{\bullet}_{X^{(1)}/\cnum_{\tau}}
 \bigr)
\stackrel{\varphi}{\lrarr}
  \hyperr^{i+1}p_{2\ast}
 \bigl(
 U_{\alpha-1}\nbigm
  \otimes\Omega^{\bullet}_{X^{(1)}/\cnum_{\tau}}
 \bigr)
\lrarr 
 \\
 C:=
 \Cok\Bigl(
  \hyperr^{i+1}p_{2\ast}
 \bigl(
 U_{\alpha-N_0}\nbigm
 \otimes\Omega^{\bullet}_{X^{(1)}/\cnum_{\tau}}
 \bigr)
\lrarr
 \hyperr^{i+1}p_{2\ast}
 \bigl(
 U_{\alpha-1}\nbigm
 \otimes\Omega^{\bullet}_{X^{(1)}/\cnum_{\tau}}
 \bigr)
 \Bigr)
\end{multline}
For any $\beta\in\real$,
$\tau\del_{\tau}+\beta$
on
$\hyperr^ip_{2\ast}\bigl(
 \Gr^U_{\beta}(\nbigm)
 \otimes\Omega^{\bullet}_{X^{(1)}/\cnum_{\tau}}
 \bigr)$
is nilpotent.
Hence, the eigenvalues of 
the endomorphism $\tau\del_{\tau}$ of 
$\hyperr^ip_{2\ast}
 \bigl(
 U_{\beta}\nbigm/U_{\gamma}\nbigm
  \otimes\Omega^{\bullet}_{X^{(1)}/\cnum_{\tau}}
 \bigr)$
are contained in 
$\closedopen{-\beta}{-\gamma}
=\bigl\{
 -\beta\leq y<-\gamma
 \bigr\}$.
In particular,
the eigenvalues of
the endomorphism $\tau\del_{\tau}$ of 
$\hyperr^ip_{2\ast}
 \bigl(
 U_{\alpha}\nbigm/U_{\alpha-1}\nbigm
  \otimes\Omega^{\bullet}_{X^{(1)}/\cnum_{\tau}}
 \bigr)$
are contained in
$\closedopen{-\alpha}{-\alpha+1}$.
Because
$C$ is contained in
$\hyperr^{i+1}p_{2\ast}
 \bigl(
 U_{\alpha-1}\nbigm/U_{\alpha-N_0}\nbigm
  \otimes\Omega^{\bullet}_{X^{(1)}/\cnum_{\tau}}
 \bigr)$,
the eigenvalues of the endomorphism
$\tau\del_{\tau}$ on $C$
are contained in
$\closedopen{-\alpha+1}{-\alpha+N_0}$.
Hence, the image of (\ref{eq;14.5.19.1})
is $0$,
and we arrived at a contradiction,
i.e.,
$\varphi=0$.
It implies that the multiplication of $\tau$
on 
$\hyperr^ip_{2\ast}
 \bigl(
 U_{\alpha}\nbigm
 \otimes\Omega^{\bullet}_{X^{(1)}/\cnum_{\tau}}
 \bigr)$
is injective.
Thus, we obtain that 
$\hyperr^ip_{2\ast}
 \bigl(
 U_{\alpha}\nbigm
 \otimes\Omega^{\bullet}_{X^{(1)}/\cnum_{\tau}}
 \bigr)$
is locally free at $\tau=0$.
Moreover,
we have
$\hyperr^ip_{2\ast}
 \bigl(
 U_{\alpha}\nbigm
 \otimes\Omega^{\bullet}_{X^{(1)}/\cnum_{\tau}}
 \bigr)_{|\tau=0}
\simeq
 \hyperr^ip_{2\ast}
 \bigl(
 (U_{\alpha}\nbigm\big/U_{\alpha-1}\nbigm)
 \otimes\Omega^{\bullet}_{X^{(1)}/\cnum_{\tau}}
 \bigr)$,
and the eigenvalues of
$\tau\del_{\tau}$ are contained in
$\closedopen{-\alpha}{-\alpha+1}$.
\qed

\subsubsection{Proof of Proposition \ref{prop;14.5.1.30}}
\label{subsection;14.5.17.2}

We have only to check the claim
around any point 
$(\tau_0,Q)$ of $D^{(1)}$.
We use a coordinate system
$(U,x_1,\ldots,x_n)$ around $Q$
as in \S\ref{subsection;14.4.23.1}.
Set $\vecp:=[\alpha\veck]$.

Let us consider the case $\tau_0\neq 0$.
We have
$\nbigm_0^{\alpha}(-D^{(1)})=
 \nbigo_{X^{(1)}}(\ast P_{\reduced}^{(1)}) v$
and 
$U_{\alpha}\nbigm
=\nbigm
=\nbigo_{X^{(1)}}(\ast D^{(1)})\,v$
around $(\tau_0,Q)$.
For $\ell_1\leq p\leq \ell$,
we set
$S(p):=\{(0,\ldots,0)\}\times\seisuu^{p-\ell_1}_{\geq 0}
 \times\{(0,\ldots,0)\}
\subset\seisuu^{\ell_1}\times\seisuu^{p-\ell_1}
 \times\seisuu^{\ell-p}
=\seisuu^{\ell}$.
We put $H^{(1)}_{\leq p}:=\bigcup_{\ell_1+1\leq i\leq p}\{x_i=0\}$
and $H^{(1)}_{>p}:=\bigcup_{p<i\leq \ell}\{x_i=0\}$
on the neighbourhood of $(\tau_0,Q)$.
We consider
\[
 \nbigm^{\leq p}:=
 \sum_{\vecn\in S(p)}
 \del^{\vecn}
 \nbigm_0^{\alpha}
=\nbigo_{X^{(1)}}(H^{(1)}_{> p})
 (\ast P^{(1)}_{\reduced})
 (\ast H_{\leq p}^{(1)})
 v.
\]
We have 
$\nbigm^{\leq \ell_1}=\nbigm_0^{\alpha}$ for any $\alpha$.

\begin{lemma}
\label{lem;14.5.14.10}
If $p\geq \ell_1+1$,
the complexes 
$x_p\nbigm^{\leq p-1}
\stackrel{\del_p}{\lrarr}
 \nbigm^{\leq p-1}$
and
$\nbigm^{\leq p}
\stackrel{\del_p}{\lrarr}
 \nbigm^{\leq p}$
are quasi-isomorphic
with respect to the inclusion.
\end{lemma}
\pf
For $j\in\seisuu_{\geq 0}$,
we consider the following:
\[
 \lefttop{p}F_j\nbigm^{\leq p}
:=\sum_{\substack{\vecn\in S(p)\\ n_p\leq j}}
 \del^{\vecn}
 \nbigm_0^{\alpha}
=\nbigo_{X^{(1)}}(H^{(1)}_{> p})
 (\ast P^{(1)}_{\reduced})
 (\ast H^{(1)}_{\leq p-1})
 x_p^{-j-1}v.
\]
We have
$\lefttop{p}F_0\nbigm^{\leq p}
=\nbigm^{\leq p-1}$.
We have the morphisms of sheaves
$\del_p:
 \lefttop{p}F_j\nbigm^{\leq p}
\lrarr
 \lefttop{p}F_{j+1}\nbigm^{\leq p}$.
We can easily check the following by direct computations.
\begin{itemize}
\item
If $j\geq 1$,
the induced morphisms of
$\nbigo_{X^{(1)}}/x_p\nbigo_{X^{(1)}}$-modules
$\lefttop{p}F_j\big/\lefttop{p}F_{j-1}
\lrarr 
\lefttop{p}F_{j+1}\big/\lefttop{p}F_{j}$
are isomorphisms.
\item
The morphism 
$\lefttop{p}F_0\nbigm^{\leq p}
=\nbigm^{\leq p-1}
\lrarr
 \lefttop{p}F_1/\lefttop{p}F_0$
is a surjection,
and the kernel is $x_p\nbigm^{\leq p-1}$.
\end{itemize}
Then, the claim of the lemma follows.
\qed

\vspace{.1in}
By Lemma \ref{lem;14.5.14.10},
the following inclusion of the complexes of sheaves
are quasi-isomorphisms:
\begin{multline}
\Bigl(
 \nbigm^{\leq p-1}(-H^{(1)}_{>p-1})
 \stackrel{a_1}{\lrarr}
 \nbigm^{\leq p-1}(-H^{(1)}_{>p-1})\cdot dx_p/x_p
\Bigr)
\lrarr
 \\
 \Bigl(
 \nbigm^{\leq p}(-H^{(1)}_{>p})
 \stackrel{a_2}{\lrarr}
 \nbigm^{\leq p}(-H^{(1)}_{>p})\cdot dx_p
\Bigr).
\end{multline}
Here, $a_i$ $(i=1,2)$
are induced by the exterior derivative
in the $x_p$-direction.
Hence,
the following inclusions of the complexes of sheaves
are quasi-isomorphisms:
\begin{equation}
\label{eq;17.7.30.1}
 \nbigm^{\leq p-1}(-H^{(1)}_{>p-1})
 \otimes
 \Omega^{\bullet}_{X^{(1)}/\cnum_{\tau}}(\log H^{(1)}_{>p-1})
\lrarr
 \nbigm^{\leq p}(-H^{(1)}_{>p})
 \otimes
 \Omega^{\bullet}_{X^{(1)}/\cnum_{\tau}}(\log H^{(1)}_{>p}).
\end{equation}

We have 
$\nbigm^{\leq \ell}(-H^{(1)}_{>\ell})
 \otimes\Omega^{\bullet}_{X^{(1)}/\cnum_{\tau}}(\log H^{(1)}_{>\ell})
=\nbigm\otimes\Omega^{\bullet}_{X^{(1)}/\cnum_{\tau}}$,
and 
\begin{multline}
\nbigm^{\leq \ell_1}(-H^{(1)}_{>\ell_1})
 \otimes\Omega^{\bullet}_{X^{(1)}/\cnum_{\tau}}(\log H^{(1)}_{>\ell_1})
=\nbigm_0^{\alpha}(-H^{(1)})
 \otimes\Omega^{\bullet}_{X^{(1)}/\cnum_{\tau}}(\log H^{(1)})
 \\
=\nbigm_0^{\alpha}(-D^{(1)})
 \otimes\Omega^{\bullet}_{X^{(1)}/\cnum_{\tau}}(\log D^{(1)}).
\end{multline}
Hence, we are done in the case $\tau_0\neq 0$.

\vspace{.1in}

Let us consider the case $\tau_0=0$.
For $0\leq p\leq \ell$,
we regard $\seisuu^{p}=
 \seisuu^{p}
 \times\{\overbrace{(0,\ldots,0)}^{\ell-p}\}
\subset
 \seisuu^{\ell}$.
We set
\[
 U_{\alpha}\nbigm^{\leq p}:=
 \sum_{\vecn\in\seisuu^p_{\geq 0}}
 \del^{\vecn}\nbigm_0^{\alpha},
\quad\quad
 \lefttop{p}F_jU_{\alpha}\nbigm^{\leq p}
 =\sum_{\substack{
 \vecn\in\seisuu^p_{\geq 0}\\
 n_p\leq j }}
 \del^{\vecn}\nbigm_0^{\alpha}.
\]
We have
$U_{\alpha}\nbigm^{\leq \ell}
=U_{\alpha}\nbigm$
and 
$\lefttop{p}F_0U_{\alpha}\nbigm^{\leq p}
=U_{\alpha}\nbigm^{\leq p-1}$.
We consider the following maps
\[
 \del_p:
 \lefttop{p}F_jU_{\alpha}\nbigm^{\leq p}
\lrarr
 \lefttop{p}F_{j+1}U_{\alpha}\nbigm^{\leq p}.
\]

The following lemma is easy to see
by Corollary \ref{cor;14.5.2.111}.
\begin{lemma}
\label{lem;14.5.1.40}
Let $s$ be a section of 
$U_{\alpha}\nbigm$
on a neighbourhood of $Y$
with a primitive expression
\[
 s=\sum_{(\vecm,j)\in\nbigs}s_{\vecm,j}.
\]
Then,
$s$ is a section of
$\lefttop{p}F_jU_{\alpha}\nbigm^{\leq p}$
if and only if
we have 
$m_i\geq 0$ $(i>p)$
and $m_p\geq -j$
for any $\vecm\in\min\pi(\nbigs)$.
\qed
\end{lemma}

\begin{lemma}
\label{lem;14.5.1.100}
If $j\geq 1$,
the following induced morphism  of sheaves
is an isomorphism:
\[
 \frac{\lefttop{p}F_jU_{\alpha}\nbigm^{\leq p}}
 {\lefttop{p}F_{j-1}U_{\alpha}\nbigm^{\leq p}}
\stackrel{\del_p}{\lrarr}
 \frac{\lefttop{p}F_{j+1}U_{\alpha}\nbigm^{\leq p}}
 {\lefttop{p}F_{j}U_{\alpha}\nbigm^{\leq p}}
\]
\end{lemma}
\pf
It is surjective by construction.
Let $s$ be a non-zero section of
$\lefttop{p}F_jU_{\alpha}\nbigm^{\leq p}$
on a neighbourhood of $Y$
with a primitive decomposition 
$s=\sum_{(\vecm,j)\in\nbigs}s_{\vecm,j}$
such that
$\del_ps$ is also a section of 
$\lefttop{p}F_jU_{\alpha}\nbigm^{\leq p}$.
We set
$s':=\sum_{m_p=-j}s_{\vecm,j}$
and 
$s'':=\sum_{m_p>-j}s_{\vecm,j}$.
Because
$\del_ps''\in \lefttop{p}F_{j}U_{\alpha}\nbigm^{\leq p}$,
we obtain 
$\del_ps'\in\lefttop{p}F_jU_{\alpha}\nbigm^{\leq p}$.
If $s'$ is non-zero,
$\del_ps'=\sum_{m_p=-j} \del_ps_{\vecm,j}$
is a primitive expression of $\del_ps'$.
We obtain
$\del_ps'\not\in \lefttop{p}F_{j}U_{\alpha}\nbigm^{\leq p}$,
and thus we arrive at a contradiction.
Hence, $s'=0$,
i.e.,
$s\in \lefttop{p}F_{j-1}U_{\alpha}\nbigm^{\leq p}$.
\qed

\begin{lemma}
\label{lem;14.5.1.101}
The kernel of the following induced surjection
is $x_pU_{\alpha}\nbigm^{\leq p-1}$:
\[
 U_{\alpha}\nbigm^{\leq p-1}
\stackrel{\del_p}{\lrarr}
\frac{\lefttop{p}F_1U_{\alpha}\nbigm^{\leq p}}
 {U_{\alpha}\nbigm^{\leq p-1}}
\]
\end{lemma}
\pf
Let $s$ be a section of 
$U_{\alpha}\nbigm^{\leq p-1}$.
We take a primitive expression
$s=\sum_{(\vecm,j)\in\nbigs}s_{\vecm,j}$.
We set
$s':=\sum_{m_p=0}s_{\vecm,j}$
and 
$s'':=\sum_{m_p>0}s_{\vecm,j}$.
We have
$s''\in x_pU_{\alpha}\nbigm^{\leq p-1}$.
Because
\[
 \del_{p}\cdot x_p\cdot
\bigl(
 g\cdot x^{-\vecdelta-\vecp}(\tau f)^jv
\bigr)
=\bigl(
 x_p\del_p(g)
-(p_p+k_pj)g
 \bigr)x^{-\vecdelta-\vecp}(\tau f)^jv
-k_pg\cdot x^{-\vecdelta-\vecp}(\tau f)^{j+1}v,
\]
we have
$\del_ps''\in U_{\alpha}\nbigm^{\leq p-1}$.
Hence, we have
$\del_ps'\in U_{\alpha}\nbigm^{\leq p-1}$.
If $s'\neq 0$,
$\del_ps'=\sum_{m_p=0}\del_ps_{\vecm,j}$
is a primitive expression of $s'$.
We obtain
$\del_ps'\not\in U_{\alpha}\nbigm^{\leq p-1}$,
and we arrive at a contradiction.
Hence, we have $s'=0$,
i.e.,
$s\in U_{\alpha}\nbigm^{\leq p-1}x_p$.
\qed

\vspace{.1in}

By Lemma \ref{lem;14.5.1.100}
and  Lemma \ref{lem;14.5.1.101},
for $1\leq p\leq \ell$,
the inclusions of the complexes
\[
 \bigl(
 x_pU_{\alpha}\nbigm^{\leq p-1}
\stackrel{\del_p}{\lrarr}
 U_{\alpha}\nbigm^{\leq p-1}
 \bigr)
\lrarr
 \bigl(
 U_{\alpha}\nbigm^{\leq p}
\stackrel{\del_p}{\lrarr}
 U_{\alpha}\nbigm^{\leq p}
 \bigr)
\]
are quasi-isomorphisms.
For $0\leq p \leq \ell$,
we set $D^{(1)}_{>p}:=\bigcup_{i=p+1}^{\ell}\{x_i=0\}$
on the neighbourhood of $(0,Q)$.
We obtain that 
the following inclusions of the complexes of sheaves
are quasi-isomorphisms:
\begin{equation}
 U_{\alpha}\nbigm^{\leq p-1}(-D^{(1)}_{>p-1})
 \otimes
 \Omega_{X^{(1)}/\cnum_{\tau}}^{\bullet}(\log D^{(1)}_{>p-1})
\lrarr 
 \\
  U_{\alpha}\nbigm^{\leq p}(-D^{(1)}_{>p})
 \otimes
 \Omega_{X^{(1)}/\cnum_{\tau}}^{\bullet}(\log D^{(1)}_{>p})
\end{equation}
We have 
$U_{\alpha}\nbigm^{\leq \ell}(-D^{(1)}_{>\ell})
 \otimes
 \Omega_{X^{(1)}/\cnum_{\tau}}^{\bullet}(\log D^{(1)}_{>\ell})
=U_{\alpha}\nbigm\otimes
 \Omega_{X^{(1)}/\cnum_{\tau}}^{\bullet}$.
We also have
\[
 U_{\alpha}\nbigm^{\leq 0}(-D^{(1)}_{>0})
 \otimes
 \Omega_{X^{(1)}/\cnum_{\tau}}^{\bullet}(\log D^{(1)}_{>0})
=\nbigm^{\alpha}_0(-D^{(1)})\otimes
 \Omega_{X^{(1)}/\cnum_{\tau}}^{\bullet}(\log D^{(1)}).
\]
Hence, Proposition \ref{prop;14.5.1.30} is proved.
\qed

\subsubsection{Proof of Proposition \ref{prop;14.4.24.10}}
\label{subsection;14.5.17.1}

We have only to check the claim
around any point of $P^{(1)}$.
We use the coordinate system as in \S\ref{subsection;14.4.23.1}.
Set $\vecp:=[\alpha\veck]$.
For any non-negative integer $N$,
we set
$G_N\bigl(
 \nbigm_0^{\alpha}(-D^{(1)})
 \bigr)
:=\sum_{j=0}^N
 \nbigo_{X^{(1)}}x^{-\vecp}
 (\tau f)^jv$.
We set
\[
 G_N\bigl(
 \nbigm_0^{\alpha}(-D^{(1)})
 \otimes\Omega^k_{X^{(1)}/\cnum_{\tau}}(\log D^{(1)})
 \bigr)
:=
 G_N\bigl(
 \nbigm_0^{\alpha}(-D^{(1)})
 \bigr)
 \otimes
 \Omega^k_{X^{(1)}/\cnum_{\tau}}(\log D^{(1)}).
\]
The derivative $d$ of the complex
$\nbigm_0^{\alpha}(-D^{(1)})
 \otimes\Omega^{\bullet}_{X^{(1)}/\cnum_{\tau}}(\log D^{(1)})$
induces
\begin{multline}
 d:
 G_N\bigl(
 \nbigm_0^{\alpha}(-D^{(1)})
 \otimes\Omega^k_{X^{(1)}/\cnum_{\tau}}(\log D^{(1)})
 \bigr)
\lrarr \\
G_{N+1}\bigl(
 \nbigm_0^{\alpha}(-D^{(1)})
 \otimes\Omega^{k+1}_{X^{(1)}/\cnum_{\tau}}(\log D^{(1)})
 \bigr).
\end{multline}
We also set
$G_{-1}\bigl(
 \nbigm_0^{\alpha}(-D^{(1)})
\otimes
 \Omega^k_{X^{(1)}/\cnum_{\tau}}(\log D^{(1)})
 \bigr)
:=\Omega^k_{f,\tau}(\alpha) v=
\Omega_{f,\tau}^{k}x^{-\vecp}v$.

Let $N\geq 0$.
Take a section
\[
 \omega
 =\sum_{j=0}^N
 \omega_j
 x^{-\vecp}(\tau f)^j\cdot v
 \in 
 G_N\bigl(
 \nbigm^{\alpha}_0(-D^{(1)})\otimes
 \Omega^k_{X^{(1)}/\cnum_{\tau}}(\log D^{(1)})
 \bigr),
\]
where
$\omega_j\in
 \Omega^k_{X^{(1)}/\cnum_{\tau}}(\log D^{(1)})$.
Suppose that
\[
 d\omega
 \in
 G_{N}\bigl(
  \nbigm_0^{\alpha}(-D^{(1)})
 \otimes\Omega^{k+1}_{X^{(1)}/\cnum_{\tau}}(\log D^{(1)})
 \bigr).
\]
Then,
$\tau df\wedge \omega_N\,(\tau f)^Nx^{-\vecp} v$
is a section of
$G_N\bigl(
 \nbigm^{\alpha}_0(-D^{(1)})
 \bigr)
\otimes
 \Omega^{k+1}_{X^{(1)}/\cnum_{\tau}}(\log D^{(1)})$.
\begin{lemma}
\label{lem;14.5.21.1}
We have
$df\wedge\omega_N
 \in
 \Omega^{k+1}_{X^{(1)}/\cnum_{\tau}}(\log D^{(1)})$,
i.e.,
$\omega_N$ is a section of
$\Omega^k_{f,\tau}$.
\end{lemma}
\pf
Let $s$ be a local section of $\nbigo_{X^{(1)}}$
such that
$(\tau f)^{N+1}s
 \in \sum_{j=0}^N \nbigo_{X^{(1)}} (\tau f)^{j}$.
We obtain that
$\tau^{N+1}s
 \in \nbigo_{X^{(1)}}f^{-1}$,
and $s\in\nbigo_{X^{(1)}}f^{-1}$.
There exists a holomorphic function $t$
such that $s=tf^{-1}$.
We obtain
$(\tau f)^{N+1}s=(\tau f)^{N} \tau t
 \in \nbigo_{X^{(1)}}(\tau f)^{N}$.

Note that
$(df/f)\wedge \omega_N\cdot (\tau f)^{N+1}$
is a section of
$\sum_{j=0}^N
 \Omega^{k+1}_{X^{(1)}/\cnum_{\tau}}(\log D^{(1)})
 \cdot (\tau f)^{j}$.
By the above argument,
we obtain that 
$(df/f)\wedge \omega_N\cdot (\tau f)^{N+1}$
is a section of
$\Omega^{k+1}_{X^{(1)}/\cnum_{\tau}}(\log D^{(1)})
 (\tau f)^{N}$.
Hence,
$\tau df\wedge\omega_N$ is a section of
$\Omega^{k+1}_{X^{(1)}/\cnum_{\tau}}(\log D^{(1)})$.
Then, we obtain that 
$df\wedge\omega_N$ is a section of
$\Omega^{k+1}_{X^{(1)}/\cnum_{\tau}}(\log D^{(1)})$.
\qed

\begin{lemma}
We have an expression
$\omega_N=
 (df/f)\wedge \kappa_1
+f^{-1}\kappa_2$,
where
$\kappa_1$ and $\kappa_2$
are local sections of
$\Omega^{k-1}_{X^{(1)}/\cnum_{\tau}}(\log D^{(1)})$
and 
$\Omega^{k}_{X^{(1)}/\cnum_{\tau}}(\log D^{(1)})$,
respectively.
\end{lemma}
\pf
Let $Q$ be any point of $P_{\reduced}$.
Let $U$ be a neighbourhood of $Q$ in $X$.
The complex
$\bigl(\Omega^{\bullet}_X(\log D),df/f\bigr)$
is acyclic on $U$
because $df/f$ is a nowhere vanishing section of
$\Omega^1_{X}(\log D)$ on $U$.
If $U$ is sufficiently small,
we can take decompositions
$\Omega_X^k(\log D) 
=\nbigb^k\oplus\nbigc^k$
such that 
the multiplication of $df/f$
induces an isomorphism
$\nbigc^k\simeq\nbigb^{k+1}$.
Then, the claim of the lemma follows.
\qed

\vspace{.1in}

If $N\geq 1$,
we have 
$f^{-1}\kappa_2(\tau f)^N
=\tau \kappa_2(\tau f)^{N-1}$.
Hence, we have
\[
 \omega-d\bigl(\kappa_1(\tau f)^{N-1}x^{-\vecp}v\bigr)
\in
 G_{N-1}\bigl(
 \nbigm_0^{\alpha}(-D^{(1)})
 \otimes\Omega^k_{X^{(1)}/\cnum_{\tau}}(\log D^{(1)})
 \bigr).
\]
We also have
$\kappa_1(\tau f)^{N-1}x^{-\vecp}v
 \in G_{N-1}$.

Let $\omega$ be a local section of
$G_N\bigl(
 \nbigm^{\alpha}_0(-D^{(1)})
 \otimes\Omega^k_{X^{(1)}/\cnum_{\tau}}(\log D^{(1)})
 \bigr)$
such that
$d\omega$ is a local section of 
\[
 G_{-1}\bigl(
 \nbigm^{\alpha}_0(-D^{(1)})
 \otimes\Omega^k_{X^{(1)}/\cnum_{\tau}}(\log D^{(1)})
 \bigr).
\]
By applying the previous argument successively,
we can find a local section
$\tau$ of
\[
 G_{N-1}\bigl(
 \nbigm^{\alpha}_0(-D^{(1)})
 \otimes\Omega^{k-1}_{X^{(1)}/\cnum_{\tau}}(\log D^{(1)})
 \bigr)
\]
such that
$\omega-d\tau$
is a local section of 
$G_{-1}\bigl(
 \nbigm^{\alpha}_0(-D^{(1)})
 \otimes\Omega^{k}_{X^{(1)}/\cnum_{\tau}}(\log D^{(1)})
 \bigr)$.

As a consequence,
we have the following.
\begin{itemize}
\item
If a local section $\omega$
of $G_N\bigl(
 \nbigm^{\alpha}_0(-D^{(1)})
 \otimes\Omega^k_{X^{(1)}/\cnum_{\tau}}(\log D^{(1)})
 \bigr)$
satisfies $d\omega=0$,
we can find 
a local section
$\tau$ of 
$G_{N-1}\bigl(
 \nbigm^{\alpha}_0(-D^{(1)})
 \otimes\Omega^{k-1}_{X^{(1)}/\cnum_{\tau}}(\log D^{(1)})
 \bigr)$
such that
$\omega-d\tau$
is a local section of 
$G_{-1}\bigl(
 \nbigm^{\alpha}_0(-D^{(1)})
 \otimes\Omega^{k}_{X^{(1)}/\cnum_{\tau}}(\log D^{(1)})
 \bigr)$.
\item
Let $\omega$ be a local section 
of 
$G_{-1}\bigl(
 \nbigm^{\alpha}_0(-D^{(1)})
 \otimes\Omega^{k}_{X^{(1)}/\cnum_{\tau}}(\log D^{(1)})
 \bigr)$
such that $d\omega=0$.
If we have a local section $\tau$
of 
$G_{N}\bigl(
 \nbigm^{\alpha}_0(-D^{(1)})
 \otimes\Omega^{k-1}_{X^{(1)}/\cnum_{\tau}}(\log D^{(1)})
 \bigr)$
such that
$\omega=d\tau$,
then we can find 
a local section $\sigma$
of 
$G_{N-1}\bigl(
 \nbigm^{\alpha}_0(-D^{(1)})
 \otimes\Omega^{k-2}_{X^{(1)}/\cnum_{\tau}}(\log D^{(1)})
 \bigr)$
such that 
$\tau-d\sigma$ is a local section of
$G_{-1}\bigl(
 \nbigm^{\alpha}_0(-D^{(1)})
 \otimes\Omega^{k-1}_{X^{(1)}/\cnum_{\tau}}(\log D^{(1)})
 \bigr)$.
We have
$\omega=d(\tau-d\sigma)$.
\end{itemize}
Then, we obtain the claim of Proposition \ref{prop;14.4.24.10}.
\qed

\section{The case of mixed twistor $\nbigd$-modules}
\label{section;14.5.13.1}

\subsection{Preliminary}
\label{subsection;17.1.20.201}

\paragraph{$\nbigr$-modules and $\nbigrtilde$-modules}

Here, let us recall the concept of $\nbigr$-modules \cite{sabbah2}.
Let $Y$ be any complex manifold.
We set $\nbigy:=\cnum_{\lambda}\times Y$.
Let $p:\nbigy\lrarr Y$ denote  the projection.
Let $\nbigd_{\nbigy}$
denote the sheaf of holomorphic differential operators
on $\nbigy$.
Let $\nbigr_Y$ denote the sheaf of subalgebras
of $\nbigd_{\nbigy}$
generated by $\lambda p^{\ast}\Theta _Y$
over $\nbigo_{\nbigy}$,
where $\Theta_Y$ denote the tangent sheaf of $Y$.
Let $\nbigrtilde_Y$ denote the sheaf of subalgebras
of $\nbigd_{\nbigy}$
generated by $\lambda^2\del_{\lambda}$ over $\nbigr_Y$.

A left $\nbigr_Y$-module is equivalent to
an $\nbigo_{\nbigy}$-module $\nbigm$
with a relative flat meromorphic connection
$\nabla^{\rel}:
 \nbigm\lrarr
 \nbigm\otimes
 \lambda^{-1}\Omega_{\nbigy/\cnum_{\lambda}}^1$.
A left $\nbigrtilde_Y$-module is equivalent to
an $\nbigo_{\nbigy}$-module $\nbigm$
with a flat meromorphic connection
$\nabla:
 \nbigm\lrarr
 \nbigm\otimes
 \lambda^{-1}
 \Omega^1_{\nbigy}(\log \lambda)$.

\paragraph{Push-forward by projection}
We recall the functoriality of $\nbigr$-modules
with respect to the push-forward by a projection
\cite{sabbah2}.
Suppose that $Y=Z\times W$ for complex manifolds
$Z$ and $W$, and that $Z$ is projective.
Let $\pi:Y\lrarr W$ denote the projection.
For any $\nbigr_Y$-module $\nbign$,
we have the $\nbigr_W$-modules
$\pi^j_{\dagger}\nbign$
$(-\dim Z\leq j\leq \dim Z)$
as follows.
Let $q_Z:\nbigy\lrarr Z$ denote the projection.
We set
$\Omegatilde^j_Z:=\lambda^{-j}q_Z^{\ast}\Omega^j_Z$.
Then, 
we have a naturally defined complex
$\Omegatilde^{\bullet}_Z\otimes\nbign$ on $\nbigy$
induced by the exterior derivative of
$\Omega_Z^{\bullet}$
and the relative flat meromorphic connection
$\nabla^{\rel}$  of $\nbign$.
We have
\[
 \pi_{\dagger}^j\nbign
\simeq
 R^{j+\dim Z}(\id_{\cnum_{\lambda}}\times\pi)_{\ast}
 \Bigl(
 \Omegatilde^{\bullet}_Z\otimes\nbign
 \Bigr).
\]

If $\nbign$ is an $\nbigrtilde_Y$-module,
then $\pi_{\dagger}^j\nbign$
are naturally $\nbigrtilde_W$-modules.

\paragraph{$V$-filtrations}

Suppose that $Y=Y_0\times\cnum_t$.
For simplicity, we assume that 
$Y_0$ is relatively compact and open 
in a larger complex manifold $Y_0'$.
Let $V_0\nbigr_Y$ be the sheaf of subalgebras
in $\nbigr_Y$ generated by 
$\lambda p^{\ast}\Theta_Y(\log t)$
over $\nbigo_{\nbigy}$.

Let $\nbign$ be an $\nbigr_Y$-module underlying 
a mixed twistor $\nbigd$-module on $Y$,
which can be extended to a mixed twistor $\nbigd$-module
on $Y_0'\times\cnum_t$.
Then, by the definition of mixed twistor $\nbigd$-module,
$\nbign(\ast t)$ is strictly specializable along $t$
as an $\nbigr_Y(\ast t)$-module.
Namely, 
for any $\lambda_0\in\cnum_{\lambda}$,
we have a neighbourhood 
$B(\lambda_0)\subset\cnum_{\lambda}$
of $\lambda_0$
and a unique filtration 
$V^{(\lambda_0)}\nbign(\ast t)_{|B(\lambda_0)\times Y}$
by coherent $V_0\nbigr_{Y}$-submodules
of $\nbign(\ast t)_{|B(\lambda_0)\times Y}$
satisfying the conditions
as in \cite[Definition 22.4.1]{Mochizuki-wild}.
Note that,
for each $a\in\real$,
we have the finite subset
$\nbigk(a,\lambda_0)\subset \real\times\cnum$
such that
$\prod_{u\in\nbigk(a,\lambda_0)}
 (-\deldel_{t}t+\eigenmap(\lambda,u))$
is nilpotent on 
$\Gr^{\Vzero}_{a}\bigl(\nbign(\ast t)_{|B(\lambda_0)\times Y}\bigr)$,
where 
we set
$\eigenmap(\lambda,(b,\beta))
=\beta-\lambda b-\betabar\lambda^2$
for $(b,\beta)\in\real\times\cnum$.

Suppose moreover that 
$\nbign$ is enhanced to an $\nbigrtilde_Y$-module.
As in \cite[Proposition 7.3.1]{sabbah2},
we have 
$\nbigk(a,\lambda_0)=\{(a,0)\}$,
and 
\[
 \lambda^2\del_{\lambda}
 \Vzero_a
 \bigl(\nbign(\ast t)_{|B(\lambda_0)\times Y}\bigr)
 \subset
 \Vzero_a
 \bigl(\nbign(\ast t)_{|B(\lambda_0)\times Y}\bigr)
\]
for any $a\in\real$.
We obtain a global filtration
$V_{a}\bigl(\nbign(\ast t)\bigr)$ $(a\in\real)$
of $V_0\nbigr_{Y}$-coherent
submodules of $\nbign(\ast t)$
by gluing $\Vzero$ $(\lambda_0\in\cnum)$,
which is uniquely characterized by the following conditions.
\begin{itemize}
\item
 We have 
 $\bigcup_{a\in\real}V_a\bigl(
 \nbign(\ast t)\bigr)
=\nbign(\ast t)$.
\item
For any $a\in\real$,
we have $\epsilon>0$
such that
$V_{a}\bigl(\nbign(\ast t)\bigr)
 =V_{a+\epsilon}\bigl(
 \nbign(\ast t)\bigr)$.
\item
We have
$t V_a\bigl(\nbign(\ast t)\bigr)
 =V_{a-1}\bigl(\nbign(\ast t)\bigr)$
and $\deldel_{t}V_a\bigl(\nbign(\ast t)\bigr)
\subset 
 V_{a+1}\bigl(\nbign(\ast t)\bigr)$
for any $a\in\real$.
\item
The induced endomorphisms
 $\deldel_{t}t+\lambda a$
 are nilpotent on 
\[
 \Gr^V_a(\nbign(\ast t))
 :=V_{a}\bigl(\nbign(\ast t)\bigr)
 \big/
 V_{<a}\bigl( \nbign(\ast t) \bigr)
\]
 for any $a\in\real$.
 Here, we set
 $V_{<a}\bigl(\nbign(\ast t)\bigr)
 =\bigcup_{b<a}
 V_b\bigl(\nbign(\ast t)\bigr)$.
\item
 $\Gr^V_a(\nbign(\ast t))$  are strict,
 i.e., flat over $\nbigo_{\cnum_{\lambda}}$.
\end{itemize}
We set
$\psitilde_a(\nbign):=
 \Gr^V_a(\nbign(\ast t))$.

\subsection{Mixed twistor $\nbigd$-modules}
\label{subsection;14.5.16.31}

\subsubsection{Wild harmonic bundles}

We use the notation in \S\ref{section;14.5.16.30}.
Let $E$ be a product line bundle on 
$X^{(1)}\setminus P_{\reduced}^{(1)}$
with a frame $e$.
The natural holomorphic structure is
denoted by $\delbar_E$.
The holomorphic line bundle 
$(E,\delbar_E)$ is equipped with
the Higgs field $\theta=d(\tau f)$,
and the metric $h$ given by $h(e,e)=1$.
Then, $(E,\delbar_E,\theta,h)$ is a wild harmonic bundle.
It is homogeneous 
with respect to the $S^1$-action on $X^{(1)}$
given by $t(\tau,Q)=(t\tau,Q)$
and $t^{\ast}e=e$,
in the sense of \cite{Mochizuki-Toda-lattice}.

We set $\nbigx^{(1)}:=\cnum_{\lambda}\times X^{(1)}$
and $\nbigp^{(1)}_{\reduced}:=
 \cnum_{\lambda}\times P^{(1)}_{\reduced}$.
Let $p:\nbigx^{(1)}\setminus\nbigp^{(1)}_{\reduced}
\lrarr
 X^{(1)}\setminus P^{(1)}_{\reduced}$
denote the projection.
Let us recall that a family of $\lambda$-flat bundles
is associated to a harmonic bundle
in this situation.
The $C^{\infty}$-bundle
$p^{-1}E$ is equipped with 
the holomorphic structure $\delbar_1$
given by
$\delbar_1p^{-1}(e)=p^{-1}(e)\cdot
\lambda d\bigl( \overline{\tau f}\bigr)$.
It is also equipped with the family of flat $\lambda$-connections
$\DD$ given by
$\DD p^{-1}(e)
=p^{-1}(e)\cdot
\bigl(
 d(\tau f) 
+\lambda d\bigl(\overline{\tau f}\bigr)
\bigr)$.
We set
$\upsilon:=p^{-1}(e)\,\exp(-\lambda \overline{\tau f})$.
It is a holomorphic frame of
$(p^{-1}E,\delbar_1)$.
We have
$\DD\upsilon=\upsilon\,d(\tau f)$.

Let $\nbige$ denote the sheaf of holomorphic sections of
$(p^{-1}E,\delbar_1)$.
Multiplying $\lambda^{-1}$
to the $(1,0)$-part of $\DD$,
we obtain the family of flat connections
\[
 \DD^f:
 \nbige
\lrarr
 \nbige\otimes
 \Bigl(
 \lambda^{-1}
 \Omega^1_{(\nbigx^{(1)}\setminus\nbigp^{(1)}_{\reduced})
 /\cnum_{\lambda}}
 \Bigr).
\]
In terms of the frame $\upsilon$,
it is given by
$\DD^f\upsilon=
 \upsilon
 \cdot \lambda^{-1}d(\tau f)$.

The bundle $p^{-1}E$ is naturally $S^1$-equivariant
with respect to the action given by
$t(\lambda,\tau,Q)=(t\lambda,t\tau,Q)$,
for which we have $t^{\ast}\upsilon=\upsilon$.
We have
$t^{\ast}\DD^f=\DD^f$.
Hence, $\DD^f$ is extended 
to a meromorphic flat connection
\[
 \nabla:
 \nbige\lrarr
 \nbige\otimes
\Bigl(
 \lambda^{-1}
 \Omega^1_{\nbigx^{(1)}\setminus\nbigp^{(1)}_{\reduced}}
 \bigl(
 \log(\{0\}\times(X^{(1)}\setminus P^{(1)}_{\reduced}))
 \bigr)
\Bigr).
\]
Namely, the derivative in the $\lambda$-direction
is induced.
(See \cite{Mochizuki-Toda-lattice}, for example.)
In terms of the frame $\upsilon$,
it is given by
$\nabla \upsilon
=\upsilon d(\lambda^{-1}\tau f)$.

The bundle $\nbige$ with $\DD^f$ gives
an $\nbigr_{X^{(1)}\setminus P^{(1)}_{\reduced}}$-module,
which is also denoted by $\nbige$.
Because $\DD^f$ is extended to 
the meromorphic connection $\nabla$ as mentioned above,
$\nbige$ is naturally enhanced to
an $\nbigrtilde_{X^{(1)}\setminus P^{(1)}_{\reduced}}$-module.

We set
$\vecS:=\{|\lambda|=1\}$.
Let $\sigma:\vecS\times (X^{(1)}\setminus P^{(1)}_{\reduced})
 \lrarr\vecS\times(X^{(1)}\setminus P^{(1)}_{\reduced})$
be given by 
$\sigma(\lambda,\tau,Q)=(-\lambda,\tau,Q)$.
We have the sesqui-linear pairing
induced by the metric $h$:
\[
 C_h:
 \nbige_{|\vecS\times (X^{(1)}\setminus D^{(1)})}
\times
 \sigma^{\ast}
\nbige_{|\vecS\times (X^{(1)}\setminus D^{(1)})}
\lrarr
\distribution_{\vecS\times (X^{(1)}\setminus D^{(1)})/\vecS}.
\]
(See \cite{sabbah2} for the concept of sesqui-linear pairings.)
In this case,
it is given by
$C_h(\upsilon,\sigma^{\ast}\upsilon)
=\exp\bigl(
 -\lambda (\overline{\tau f})
+\lambdabar(\tau f)
 \bigr)
=\exp\bigl(
 2\sqrt{-1}
 \Image(\lambdabar \tau f)
 \bigr)$.
Thus, we obtain 
an $\nbigr$-triple $(\nbige,\nbige,C_h)$
on $X^{(1)}\setminus P^{(1)}_{\reduced}$.
It is naturally enhanced to
an $\nbigrtilde$-triple.

\subsubsection{The associated pure twistor $\nbigd$-modules}

According to \cite{Mochizuki-wild},
$(\nbige,\nbige,C_h)$ is uniquely extended to 
a pure twistor $\nbigd$-module
$\gbigt=(\gbigm,\gbigm,\gbigc)$ on $X^{(1)}$
of weight $0$ with the polarization $(\id,\id)$
such that 
(i) the strict support of $\gbigt$ is $X^{(1)}$,
(ii) $\gbigt_{|X^{(1)}\setminus P^{(1)}_{\reduced}}$
is identified with $(\nbige,\nbige,C_h)$.
Let us describe the $\nbigr_{X^{(1)}}$-module $\gbigm$.

The holomorphic bundle $\nbige$ is naturally extended
to an $\nbigo_{\nbigx^{(1)}}(\ast \nbigp^{(1)}_{\reduced})$-module
$\nbigq\nbige:=
 \nbigo_{\nbigx^{(1)}}(\ast \nbigp^{(1)}_{\reduced})\upsilon$
with a meromorphic connection
\[
 \nabla:
 \nbigq\nbige
\lrarr
 \nbigq\nbige\otimes
 \Bigl(
 \lambda^{-1}\Omega_{\nbigx^{(1)}}
 \bigl(\log(\{0\}\times X^{(1)})\bigr)
 \Bigr)
\]
given by 
$\nabla\upsilon=\upsilon\cdot d(\lambda^{-1}\tau f)$.
It naturally induces
an $\nbigrtilde_{X^{(1)}}$-module
denoted by $\nbigq\nbige$.

\begin{lemma}
\label{lem;14.5.2.1}
We have a natural isomorphism
of $\nbigr_{X^{(1)}}$-modules
$\gbigm\simeq
 \nbigq\nbige$.
\end{lemma}
\pf
We take a projective birational morphism
$F:X^{(2)}\lrarr X^{(1)}$
such that 
(i) $P^{(2)}_{\reduced}:=
 F^{-1}(P^{(1)}_{\reduced})$ is simply normal crossing,
(ii) $X^{(2)}\setminus P^{(2)}\simeq
 X^{(1)}\setminus P^{(1)}$,
(iii) the zeroes and the poles of $F^{\ast}(\tau f)$
 are separated.
Let $(E',\delbar_{E'},\theta',h')$
be the pull back of the harmonic bundle
$(E,\delbar_E,\theta,h)$ by $F$.
Then, it is a good wild harmonic bundle
on $(X^{(2)},D^{(2)})$.
We have the associated pure twistor $\nbigd$-module
$\gbigt'=(\gbigm',\gbigm',\gbigc')$ on $X^{(2)}$.
By the construction of $\gbigm'$ in \cite{Mochizuki-wild},
we have a natural isomorphism
$\gbigm'(\ast \nbigp^{(2)})
\simeq
 F^{\ast}\nbigq\nbige$.
Then, we have a natural isomorphism
$F_{\dagger}\bigl(
 \gbigm'\bigr)(\ast \nbigp_{\reduced}^{(1)})
\simeq
 \nbigq\nbige$.
Note that
$\gbigm$ is a direct summand of
$F_{\dagger}(\gbigm')$
such that
$\gbigm_{|\nbigx^{(1)}\setminus \nbigp_{\reduced}^{(1)}}
=F_{\dagger}(\gbigm')_{|
 \nbigx^{(1)}\setminus\nbigp_{\reduced}^{(1)}}$,
and that $\gbigm$ is strictly $S$-decomposable.
Because the strict support of $\gbigm$ is $X^{(1)}$,
a local section of $\gbigm$ is $0$
if its support is contained in $\nbigp^{(1)}$.
Hence, we obtain 
$\gbigm\subset
 \gbigm(\ast \nbigp^{(1)}_{\reduced})
\simeq
F_{\dagger}(\gbigm')(\ast\nbigp^{(1)}_{\reduced})
\simeq
 \nbigq\nbige$.
Moreover, we have
$\nbigo_{\nbigx^{(1)}}(-N\nbigp_{\reduced}^{(1)})
 \upsilon
\subset
 \gbigm$
for some positive integer $N$.
By using $\deldel_{\tau}\upsilon=f\upsilon$,
we obtain that
$\upsilon$ is a section of $\gbigm$,
and hence
$\gbigm\supset
 \nbigq\nbige$.
\qed

\vspace{.1in}
In particular, $\gbigm$ is naturally
an $\nbigrtilde_{X^{(1)}}$-module.

\subsubsection{The associated mixed twistor $\nbigd$-modules}
\label{subsection;14.5.2.21}

Recall that $H$ is a hypersurface of $X$
such that $D=H\cup P_{\reduced}$
and $\codim (H\cap P_{\reduced})\geq 2$.
We set $H^{(1)}:=\cnum_{\tau}\times H$.
We set
$\nbigd^{(1)}:=\cnum_{\lambda}\times D^{(1)}$
and
$\nbigh^{(1)}:=\cnum_{\lambda}\times H^{(1)}$.

We have the localization 
of $\gbigt$ along $H^{(1)}$
in the category of mixed twistor $\nbigd$-modules on $X^{(1)}$.
(See \cite{Mochizuki-MTM}.)
It is denoted by $\gbigt[\ast H^{(1)}]$.
It consists of an $\nbigr_{X^{(1)}}$-triple
$\bigl(\gbigm[!H^{(1)}],\gbigm[\ast H^{(1)}],
 \gbigc[\ast H^{(1)}]\bigr)$
with a weight filtration $\nbigw$.
We have a natural morphism of $\nbigr_{X^{(1)}}$-triples
$\gbigt\lrarr\gbigt[\ast H^{(1)}]$.
Let us describe the $\nbigr_X(\ast \tau)$-module 
$\gbigm[\ast H^{(1)}](\ast \tau)$.

We have the $\nbigr_{X^{(1)}}(\ast\tau)$-module
$\nbigq\nbige(\ast\nbigh^{(1)})(\ast \tau)$
and the  $\nbigo_{\nbigx^{(1)}}(\ast\tau)$-submodule
$\nbigq\nbige(\nbigh^{(1)})(\ast\tau)$.
Let $\nbigmtilde$ denote
the $\nbigr_{X^{(1)}}(\ast \tau)$-submodule
of 
$\nbigq\nbige(\ast\nbigh^{(1)})(\ast\tau)$
generated by 
$\nbigq\nbige(\nbigh^{(1)})(\ast\tau)$
over $\nbigr_{X^{(1)}}$:
\[
\nbigmtilde:=
 \nbigr_{X^{(1)}}\cdot
 \Bigl(
 \nbigq\nbige(\nbigh^{(1)})(\ast\tau)
 \Bigr)
\subset
 \nbigq\nbige(\ast\nbigh^{(1)})(\ast\tau).
\]
It is naturally an $\nbigrtilde_{X^{(1)}}(\ast\tau)$-submodule of
$\nbigq\nbige(\ast\nbigh^{(1)})(\ast\tau)$.

\begin{lemma}
\label{lem;17.1.20.100}
We have a natural isomorphism 
of $\nbigr_{X^{(1)}}(\ast\tau)$-modules:
\begin{equation}
\label{eq;17.1.20.101}
 \gbigm[\ast H^{(1)}](\ast\tau)
 \simeq
 \nbigmtilde.
\end{equation}
\end{lemma}
\pf
We have a natural injection
$\gbigm\lrarr\gbigm[\ast H^{(1)}]$
whose restriction to
$\nbigx^{(1)}\setminus\nbigh^{(1)}$
is an isomorphism.
It induces 
$\gbigm(\ast\nbigh^{(1)})(\ast \tau)
\simeq 
 \bigl(
 \gbigm[\ast H^{(1)}]
 \bigr)(\ast\nbigh^{(1)})(\ast\tau)$.
We obtain a morphism
$\gbigm[\ast H^{(1)}](\ast\tau)\lrarr
 \gbigm(\ast\nbigh^{(1)})(\ast\tau)
=\nbigq\nbige(\ast\nbigh^{(1)})(\ast\tau)$.

Take a locally defined holomorphic function $g$
on $U\subset X^{(1)}$ such that $g^{-1}(0)=H^{(1)}$.
Let $\iota_g:U\lrarr U\times\cnum_t$
denote the graph of $g$.
We have
$\iota_{g\dagger}(\gbigm[\ast H^{(1)}])
\simeq
\bigl(
 \iota_{g\dagger}\gbigm
\bigr)[\ast t]
\subset
 \iota_{g\dagger}(\gbigm(\ast\nbigh^{(1)}))$.
(See \S3.1.2, \S3.3 and \S5.4.1 of \cite{Mochizuki-MTM}.)
It implies that
a local section of 
$\gbigm[\ast H^{(1)}]$ is $0$
if its support is contained in $\nbigh^{(1)}$.
Hence, we may regard 
$\gbigm[\ast H^{(1)}](\ast\tau)$
as an $\nbigr_{X^{(1)}}$-submodule of
$\nbigq\nbige(\ast\nbigh^{(1)})(\ast\tau)$.
By the construction of
$\nbigr$-module 
$\gbigm[\ast H^{(1)}]_{|
 \nbigx^{(1)}\setminus\{\tau=0\}}$
in \S5.3.3 of \cite{Mochizuki-MTM},
outside $\{\tau=0\}$,
$\gbigm[\ast H^{(1)}]$ is generated by
$\nbigq\nbige(\nbigh^{(1)})$
over $\nbigr_{X^{(1)}}$.
Hence,
we have
$\gbigm[\ast H^{(1)}](\ast\tau)_{|
 \nbigx^{(1)}\setminus\{\tau=0\}}
=
\nbigmtilde_{|
 \nbigx^{(1)}\setminus\{\tau=0\}}$.
Because the actions of $\tau$
on $\gbigm[\ast H^{(1)}](\ast\tau)$
and $\nbigmtilde$
are invertible,
we obtain 
$\gbigm[\ast H^{(1)}](\ast\tau)
=\nbigmtilde$
in $\nbigq\nbige(\ast\nbigh^{(1)})(\ast\tau)$.
\qed

\vspace{.1in}

Because the restriction of (\ref{eq;17.1.20.101})
to $\nbigx^{(1)}\setminus \nbigd^{(1)}$ is 
an isomorphism of 
$\nbigrtilde_{X^{(1)}\setminus D^{(1)}}$-modules,
we can easily deduce that
(\ref{eq;17.1.20.101}) is an isomorphism
of $\nbigrtilde_{X^{(1)}}$-modules.
We have the $V$-filtration of $\nbigmtilde$.
We put $U_{\alpha}\nbigmtilde:=V_{\alpha-1}\nbigmtilde$
for any $\alpha\in\real$,
which is a unique filtration satisfying the condition
in \S\ref{subsection;17.1.20.20}.
The filtration $U_{\alpha}$ $(\alpha\in\real)$
is also called the $V$-filtration in this paper.
We shall explicitly describe 
the filtration in \S\ref{subsection;14.5.14.1}.

\subsubsection{Push-forward}

Let $p_1:X^{(1)}\lrarr\cnum_{\tau}$ denote the projection.
As the cohomology of the push-forward,
we obtain $\nbigr_{\cnum_{\tau}}$-triples
for $j\in\seisuu$:
\[
 p_{1\dagger}^j\gbigt[\ast H^{(1)}]
=\Bigl(
 p_{1\dagger}^{-j}\gbigm[!H^{(1)}],
 p_{1\dagger}^j\gbigm[\ast H^{(1)}],
 p_{1\dagger}^j\gbigc[\ast H^{(1)}]
 \Bigr).
\]
They are mixed twistor $\nbigd$-modules
with the induced filtrations $\nbigw$.

\begin{lemma}
For any $i$,
$p_{1\dagger}^i\bigl(
\gbigt[\ast H^{(1)}]
 \bigr)(\ast \tau)$ 
are admissible variations of mixed twistor structure
on $(\cnum_{\tau},0)$ 
in the sense of {\rm\cite[\S9]{Mochizuki-MTM}}.
In particular,
$p^{i}_{1\dagger}\gbigm[\star H^{(1)}](\ast \tau)$
are locally free
$\nbigo_{\cnum^2_{\lambda,\tau}}(\ast\tau)$-modules.
\end{lemma}
\pf
We set $\gbigm^1[\star H^{(1)}]:=
 \gbigm[\star H^{(1)}]\big/
 (\lambda-1)
  \gbigm[\star H^{(1)}]$.
We have
\[
 p_{1\dagger}^i\bigl(
 \gbigm[\star H^{(1)}]
 \bigr)
\big/
 (\lambda-1)
 p_{1\dagger}^i\bigl(
 \gbigm[\star H^{(1)}]
 \bigr)
\simeq
 p_{1\dagger}^i\bigl(
 \gbigm^1[\star H^{(1)}]
 \bigr).
\]
Because 
$\gbigm^1[\ast H^{(1)}]
\simeq
\nbigm$,
$p_{1\dagger}^i\bigl(
 \gbigm^1[\ast H^{(1)}]
 \bigr)_{|\cnum_{\tau}^{\ast}}$
is a locally free
$\nbigo_{\cnum_{\tau}^{\ast}}$-module
according to Corollary \ref{cor;14.5.2.2}.
(We can check the claim directly
more easily than Corollary \ref{cor;14.5.2.2}.)
Hence,
$p_{1\dagger}^i\bigl(
\gbigt[\ast H^{(1)}]
 \bigr)_{|\cnum_{\tau}^{\ast}}$
is a smooth $\nbigr_{\cnum^{\ast}_{\tau}}$-triple.
Because
$\Gr^{\nbigw}_wp_{1\dagger}^i(\gbigt[\ast H^{(1)}])$
are polarizable pure twistor $\nbigd$-modules,
$\Gr^{\nbigw}_wp_{1\dagger}^i\bigl(\gbigt[\ast H^{(1)}]\bigr)(\ast\tau)$
are obtained as the canonical prolongation
of good wild polarizable variation of pure twistor structure
of weight $w$.
Hence,
the tuple of
the identity
$\id:\cnum_{\tau}\lrarr\cnum_{\tau}$,
the open subset
$\cnum_{\tau}^{\ast}\subset\cnum_{\tau}$,
and 
$p_{1\dagger}^i\bigl(\gbigt[\ast H^{(1)}]\bigr)(\ast\tau)$
with $\nbigw$ 
gives a cell of the mixed twistor $\nbigd$-module
$p_{1\dagger}^i\bigl(\gbigt[\ast H^{(1)}]\bigr)$
in the sense of \S11.1 of \cite{Mochizuki-wild}.
Then, according to Proposition 11.1 of \cite{Mochizuki-wild},
$p_{1\dagger}^i\bigl(\gbigt[\ast H^{(1)}]\bigr)(\ast \tau)$
is an admissible variation of mixed twistor structure.
\qed

\vspace{.1in}

The $\nbigr_{\cnum_{\tau}}(\ast \tau)$-modules
$p^i_{1\dagger}\gbigm[\ast H^{(1)}](\ast \tau)$
are strictly specializable along $\tau$.
As in the case of $\nbigmtilde$,
we have a unique filtration
$V_{a}p^i_{1\dagger}\bigl(
 \gbigm[\ast H^{(1)}](\ast\tau)
\bigr)$ $(a\in\real)$
by $\lefttop{\tau}V_0\nbigr_{\cnum_{\tau}}$-coherent 
submodules of 
$p^i_{1\dagger}\bigl(\gbigm[\ast H^{(1)}](\ast\tau)\bigr)$
satisfying the conditions in 
\S\ref{subsection;17.1.20.201}.
We set 
$U_{\alpha}p^i_{1\dagger}\bigl(
 \gbigm[\ast H^{(1)}](\ast\tau)
\bigr):=
 V_{\alpha-1}p^i_{1\dagger}\bigl(
 \gbigm[\ast H^{(1)}](\ast\tau)
\bigr)$.

\begin{corollary}
$U_{\alpha}\bigl(
 p^i_{1\dagger}\gbigm[\ast H^{(1)}](\ast \tau)
 \bigr)$
are locally free $\nbigo_{\cnum^2_{\lambda,\tau}}$-modules
for any $\alpha\in\real$.
\end{corollary}
\pf
Because $p_{1\dagger}^i\bigl(\gbigt[\ast H^{(1)}]\bigr)(\ast\tau)$
are admissible variations of mixed twistor structure,
the $\nbigr_{\cnum_{\tau}}(\ast \tau)$-modules
$p^i_{1\dagger}\gbigm[\ast H^{(1)}](\ast \tau)$
are good-KMS in the sense of
\cite[\S5.11]{Mochizuki-MTM}.
By using the $\cnum^{\ast}$-equivariance,
we can easily observe that
$p^i_{1\dagger}\gbigm[\ast H^{(1)}](\ast \tau)$
are regular-KMS in the sense of
\cite[\S5.11]{Mochizuki-MTM}.
Hence, 
$U_{\alpha}\bigl(
 p^i_{1\dagger}\gbigm[\ast H^{(1)}](\ast \tau)
\bigr)$
are locally free $\nbigo_{\cnum^2_{\lambda,\tau}}$-modules.
\qed

\subsubsection{Hodge structure on the nearby cycle sheaf}
\label{subsection;17.7.30.110}

We consider the $\cnum^{\ast}$-action on 
$\cnum_{\lambda}\times\cnum_{\tau}$
given by $t(\lambda,\tau)=(t\lambda,t\tau)$.
It induces a $\cnum^{\ast}$-action on $\nbigx^{(1)}$.
The $\nbigr_{X^{(1)}}$-modules
$\gbigm[\star H^{(1)}]$ $(\star=\ast,!)$
are $\cnum^{\ast}$-equivariant,
and the sesqui-linear pairing of
$\gbigm[!H^{(1)}]$
and 
$\gbigm[\ast H^{(1)}]$
is equivariant with respect to
the $S^1$-action.

We have the induced natural $\cnum^{\ast}$-actions
on $p^i_{1\dagger}\gbigm[\star H^{(1)}]$.
The induced $S^1$-action is compatible
with the sesqui-linear pairing
of 
$p^{-i}_{1\dagger}\gbigm[!H^{(1)}]$
and 
$p^{i}_{1\dagger}\gbigm[\ast H^{(1)}]$.

For $i\in\seisuu$,
we have the nearby $\nbigr_{\{0\}}$-triple
\[
 \psitilde_ap_{1\dagger}^i\gbigt[\ast H^{(1)}]
=\bigl(
 \psitilde_ap_{1\dagger}^{-i}\gbigm[!H^{(1)}],
 \psitilde_ap_{1\dagger}^i\gbigm[\ast H^{(1)}],
 \psitilde_ap_{1\dagger}^i\gbigc[\ast H^{(1)}]
 \bigr).
\]
Here,
$\psitilde_ap_{1\dagger}^i\gbigc[\ast H^{(1)}]$
denotes the induced sesqui-linear pairing.
With the relative monodromy filtration $W$,
it is a mixed twistor structure.
We have the induced $\cnum^{\ast}$-actions on
$\psitilde_ap_{1\dagger}^{-i}\gbigm[!H^{(1)}]$
and 
$\psitilde_ap_{1\dagger}^i\gbigm[\ast H^{(1)}]$,
which are compatible with the filtration $W$.
The sesqui-linear pairing
$\psitilde_ap_{1\dagger}^i\gbigc[\ast H^{(1)}]$
is compatible with the $S^1$-action.
Hence, the mixed twistor structure
comes from a complex mixed Hodge structure.

Namely, 
let $\gbigh^i_a$ denote the vector space
obtained as the fiber of
$\psitilde_ap_{1\dagger}^i\gbigm[\ast H^{(1)}]$
over $1\in\cnum_{\lambda}$.
Then, we have the decreasing filtrations $F$ and $G$
on $\gbigh^i_a$
such that 
the $\nbigr$-modules
$\psitilde_ap_{1\dagger}^i\gbigm[\ast H^{(1)}]$
and 
$\psitilde_ap_{1\dagger}^{-i}\gbigm[!H^{(1)}]$
with the $\cnum^{\ast}$-action
are isomorphic to
the analytification of the Rees modules
of the filtered vector spaces $(\gbigh^i_a,F)$
and $(\gbigh^i_a,F)^{\lor}$.
Moreover, we have the increasing filtration $W$
on $\gbigh^i_a$ induced by 
the filtration $W$ of
$\psitilde_ap_{1\dagger}^i\gbigm[\ast H^{(1)}]$
so that 
$(\gbigh^i_a,F,G,W)$ is a complex mixed Hodge structure.

For any $k\in\seisuu$,
let $N^k_a$ denote the morphisms
$\psitilde_ap_{1\dagger}^k\gbigm[\star H]
\lrarr
 \lambda^{-1}
\psitilde_ap_{1\dagger}^k\gbigm[\star H]$
$(\star=\ast,!)$
induced as the nilpotent part of $\tau\del_{\tau}$.
Let $\gbign^i_a$ denote the endomorphism 
of $\gbigh^i_a$ induced by $N^i_a$.
Recall that $(-N^{-i}_a,-N^i_a)$ gives a morphism
of mixed twistor structures
$\psitilde_ap_{1\dagger}^i(\gbigt[\ast H])
\lrarr
\psitilde_ap_{1\dagger}^i(\gbigt[\ast H])
\otimes
 \newTate(-1)$.
Hence,
$\gbign^i_a$ induces 
a morphism of mixed complex Hodge structure
$(\gbigh^i_a,F,G,W)
\longrightarrow
 (\gbigh^i_a,F,G,W)(-1)$,
where $(-1)$ indicates the Tate twist.
In particular, we obtain the following.
\begin{lemma}
\label{lem;17.7.30.100}
$\gbign^i_a(F^j)
=\Image(\gbign^i_a)\cap F^{j-1}$
for any $i\in\seisuu$, $j\in\seisuu$
and $a\in\real$.
\hfill\qed
\end{lemma}

\subsubsection{Functoriality of $V$-filtrations}

Let $\ptilde_1:\nbigx^{(1)}\lrarr \cnum_{\lambda}\times\cnum_{\tau}$
denote the projection.
Let $\Omega^1_{\nbigx^{(1)}/\cnum^2_{\lambda,\tau}}$
denote the kernel of
$\Omega^1_{\nbigx^{(1)}}
\lrarr
 \ptilde^{\ast}_{1}\Omega^1_{\cnum^2_{\lambda,\tau}}$.
We set 
$\Omegatilde^1_{\nbigx^{(1)}/\cnum^2_{\lambda,\tau}}
 :=
 \lambda^{-1}
 \Omega^1_{\nbigx^{(1)}/\cnum^2_{\lambda,\tau}}$
and 
$\Omegatilde^k_{\nbigx^{(1)}/\cnum^2_{\lambda,\tau}}
 :=
\bigwedge^k
\Omegatilde^1_{\nbigx^{(1)}/\cnum^2_{\lambda,\tau}}$.
Let $\pi:X^{(1)}\lrarr X$ denote the projection.
Put $\pitilde:=\id_{\cnum_{\lambda}}\times\pi$.
We set
$\pi^{\ast}\nbigr_X:=
 \nbigo_{\nbigx^{(1)}}\otimes_{\pitilde^{-1}\nbigo_{\nbigx}}
 \pitilde^{-1}\nbigr_X$.
For any $\pi^{\ast}\nbigr_{X}$-module $\nbigb$,
we have a naturally defined complex
$\nbigb\otimes
 \Omegatilde^{\bullet}_{\nbigx^{(1)}/\cnum^2_{\lambda,\tau}}$.

By the functoriality of mixed twistor $\nbigd$-modules
with respect to the push-forward via projective morphisms,
$p^j_{1\dagger}
 \Gr^U_{\alpha}\gbigm[\ast H^{(1)}]$
are strict for any $\alpha\in\real$.
Hence, we obtain the following isomorphism
(see the proof of Theorem 3.18 of \cite{sabbah2}):
\begin{equation}
\label{eq;17.7.30.30}
 U_{\alpha}\bigl(
 p_{1\dagger}^j\gbigm[\ast H^{(1)}](\ast \tau)
 \bigr)
\simeq
 \hyperr^j\ptilde_{1\ast}
 \Bigl(
 U_{\alpha}\bigl(
 \gbigm[\ast H^{(1)}](\ast \tau)
 \bigr)
 \otimes
 \Omegatilde^{\bullet}_{\nbigx^{(1)}/\cnum^2_{\lambda,\tau}}[\dim X]
 \Bigr).
\end{equation}
It is isomorphic to 
$\hyperr^j\ptilde_{1\ast}
 \bigl(
 U_{\alpha}\nbigmtilde
 \otimes
 \Omegatilde^{\bullet}_{\nbigx^{(1)}/\cnum^2_{\lambda,\tau}}[\dim X]
 \bigr)$.

\subsection{Main theorem}
\label{subsection;17.1.20.200}

\subsubsection{Quasi-isomorphism}

Let $q_X:\nbigx^{(1)}\lrarr X$ denote the projection.
We set
$\Omegatilde^{k}_{f,\lambda,\tau}(\alpha):=
 \lambda^{-k}q_X^{\ast}\Omega_f^{k}(\alpha)$.
We have the following morphism of sheaves:
\[
  d+\lambda^{-1}\tau df:
 \Omegatilde^k_{f,\lambda,\tau}(\alpha)
\lrarr
  \Omegatilde^{k+1}_{f,\lambda,\tau}(\alpha).
\]
We have 
$\bigl(d+\lambda^{-1}\tau df\bigr)^2=0$.
We shall prove the following theorem 
later in \S\ref{subsection;14.5.14.20}.
\begin{theorem}
\label{thm;14.5.2.20}
We have a natural quasi-isomorphism
\[
 \bigl(
 \Omegatilde_{f,\lambda,\tau}^{\bullet}(\alpha),
 d+\lambda^{-1}\tau df
\bigr)
\simeq
U_{\alpha}\bigl(
 \gbigm[\ast H^{(1)}](\ast \tau)
 \bigr)
 \otimes
 \Omegatilde^{\bullet}_{\nbigx^{(1)}/\cnum^2_{\lambda,\tau}}.
\]
It is equivariant with respect to
the $\cnum^{\ast}$-action.
\end{theorem}

\subsubsection{Consequences}
\label{subsection;14.12.28.10}

We shall explain some consequences
of Theorem \ref{thm;14.5.2.20}.

\paragraph{Proof of Theorem \ref{thm;14.4.24.20}}

We set 
$\nbigktilde^j_f(\alpha):=
 \hyperr^j\ptilde_{1\ast}\Omegatilde^{\bullet}_{f,\lambda,\tau}(\alpha)$.
By Theorem \ref{thm;14.5.2.20}
and the isomorphism (\ref{eq;17.7.30.30}),
we have
\[
\nbigktilde^j_f(\alpha)
\simeq
 U_{\alpha}\bigl(
 p_{1\dagger}^{j+\dim X}\gbigm[\ast H^{(1)}](\ast \tau)
 \bigr)
\]
for any $j$.
Because 
$U_{\alpha}\bigl(
 p_{1\dagger}^j\gbigm[\ast H^{(1)}](\ast \tau)
 \bigr)$ are locally free
$\nbigo_{\cnum^2_{\lambda,\tau}}$-modules
for any $j$,
we obtain that 
the fiber of 
$\nbigktilde^j_f(\alpha)$
at $(\lambda,\tau)$
is quasi-isomorphic to
the $j$-th hypercohomology group of
$(\Omega^{\bullet}_f(\alpha),\lambda d+\tau df)$
by using a general result on any perfect complexes
(Lemma \ref{lem;14.5.14.21} below).
In particular,
the dimension of the hypercohomology groups
are independent of $(\lambda,\tau)$,
i.e.,
Theorem \ref{thm;14.4.24.20} holds.

\begin{lemma}
\label{lem;14.5.14.21}
Let $Y$ be any complex manifold.
Let $C^{\bullet}$ be a bounded complex of
locally free $\nbigo_Y$-modules.
Let $\nbigf$ be any $\nbigo_Y$-module.
If the cohomology sheaves
$\nbigh^j(C^{\bullet})$ $(j\in\seisuu)$ are locally free
$\nbigo_Y$-modules,
we naturally have
$\nbigh^j(C^{\bullet})
 \otimes\nbigf
\simeq
 \nbigh^j(C^{\bullet}\otimes\nbigf)$
for any $j$.
\qed
\end{lemma}

\paragraph{Proof of Theorem \ref{thm;14.4.24.21}}

The $\nbigr_{\cnum_{\tau}}$-module
$p^j_{1\dagger}\gbigm[\ast H^{(1)}]$
is naturally an $\nbigrtilde_{\cnum_{\tau}}$-module,
i.e., it is equipped with the action of
$\lambda^2\del_{\lambda}$.
The action of $\lambda^2\del_{\lambda}$
preserves
$\nbigktilde^{j+\dim X}_f(\alpha)=
 U_{\alpha}\bigl(
p^j_{1\dagger}\gbigm[\ast H^{(1)}]
(\ast\tau)
 \bigr)$.
It means that $\nbigktilde^i_f(\alpha)$
is equipped with a meromorphic connection
$\nabla$ such that
$\lambda\nabla_{\tau\del_{\tau}}
 \nbigktilde^i_f(\alpha)
\subset
\nbigktilde^i_f(\alpha)$
and
$\nabla_{\lambda^2\del_{\lambda}}
 \nbigktilde^i_f(\alpha)
 \subset
 \nbigktilde^i_f(\alpha)$.
The holomorphic bundle $\nbigktilde^i_f(\alpha)$
is $\cnum^{\ast}$-equivariant
with respect to the action
$t(\lambda,\tau)=(t\lambda,t\tau)$.
The meromorphic connection $\nabla$ is
equivariant by the construction,
and the $\cnum^{\ast}$-action on the bundle
is equal to the parallel transport by the connection.
The induced vector bundle
on $\proj^1$
is isomorphic to
$\nbigk^i_f(\alpha)$.
It is equipped with the induced 
meromorphic connection $\nabla$
such that
$\nabla\nbigk^i_f(\alpha)
\subset
 \nbigk^i_f(\alpha)
 \otimes
 \Omega^1_{\proj^1}(\{0\}+2\{\infty\})$.
Here,
$0$ corresponds to
$[\lambda:\tau]=[1:0]$
and 
$\infty$ corresponds to
$[\lambda:\tau]=[0:1]$.
By the construction,
it is equal to the connection in
Corollary \ref{cor;14.5.2.2}
around $0$.
Thus,
we obtain Theorem \ref{thm;14.4.24.21}.

\paragraph{Proof of Theorem \ref{thm;17.7.30.20}}

For $\alpha-1<b\leq\alpha$,
let $U_{b}\bigl(
 \nbigktilde^i_f(\alpha)_{|\tau=0}
\bigr)$
denote the image of 
\[
U_{b}
 p^{i-\dim X}_{1\dagger}\bigl(
 \gbigm[\ast H^{(1)}]
 \bigr)_{|\tau=0}
\lrarr 
U_{\alpha}
 p^{i-\dim X}_{1\dagger}\bigl(
 \gbigm[\ast H^{(1)}]
 \bigr)_{|\tau=0}
=\nbigktilde^i_f(\alpha)_{|\tau=0}.
\]
We have
$\psitilde_{b-1}
 p^{i-\dim X}_{1\dagger}\bigl(
 \gbigm[\ast H^{(1)}]
\bigr)
 \simeq
 \Gr^{U}_{b}\bigl(
 \nbigktilde^i_f(\alpha)_{|\tau=0}
 \bigr)$.
We also have the filtration $U$
on $\nbigk^i_f(\alpha)_{|0}$
as in \S\ref{subsection;17.7.30.101}.
By the constructions,
we have the following natural identifications:
\[
 \gbigh^{i-\dim X}_{b-1}\simeq
 \Gr^U_{b}\bigl(
 \nbigktilde^i_f(\alpha)_{|\tau=0}
 \bigr)_{|\lambda=1}
\simeq
 \Gr^U_b
 \nbigk^i_f(\alpha)_{|0}.
\]
(See \S\ref{subsection;17.7.30.110} for $\gbigh^k_a$.)

By Theorem \ref{thm;14.4.24.20},
we have the filtration $F^{(1)}$
on $\nbigktilde^i_f(\alpha)$ 
induced by the truncations
as in \S\ref{subsection;17.7.30.101}
(see (\ref{eq;17.7.30.102})).
In particular,
it induces a filtration 
of $\nbigktilde^i_f(\alpha)_{|\cnum_{\lambda}\times\{0\}}$.
We also have a filtration $F^{(2)}$ 
on $\nbigktilde^i_f(\alpha)_{|\cnum_{\lambda}\times\{0\}}$
corresponding to the $\cnum^{\ast}$-action,
according to the Rees construction.
It is easy to observe that $F^{(1)}=F^{(2)}$.
Hence, 
the Hodge filtration $F$ on $\gbigh^{i-\dim X}_{b-1}$
is equal to the filtration induced by
$F^{(1)}$.
Then, the claim of Theorem \ref{thm;17.7.30.20}
follows from Lemma \ref{lem;17.7.30.100}.

\subsection{A description of the $V$-filtration}
\label{subsection;14.5.14.1}

\subsubsection{The $\nbigr_{X^{(1)}}(\ast \tau)$-module $\nbigmtilde$}

We consider
the $\nbigr_{X^{(1)}}(\ast \tau)$-module
$\nbigmtilde$
given in \S\ref{subsection;14.5.2.21}.
It has the global section $\upsilon$.
For a local coordinate system of $X$ around a point of $D$
as in \S\ref{subsection;14.4.23.1},
we have
$\deldel_i \upsilon=-k_i(\tau fx_i^{-1})\upsilon$
and
$\deldel_{\tau}\upsilon=f\upsilon$.
We have
\[
 \tau\deldel_{\tau}(x^{-\vecdelta+\vecm}\upsilon)
=x^{-\vecdelta+\vecm}\tau f\upsilon,
\quad
 \deldel_{i}x_ix^{-\vecdelta+\vecm}\upsilon
=x^{-\vecdelta+\vecm}
(\lambda m_i-k_i\tau f)\upsilon.
\]
Hence, for any $i$ with $k_i\neq 0$,
we have
\begin{equation}
\label{eq;14.5.2.22}
 (\tau\deldel_{\tau}+k_i^{-1}\deldel_ix_i)
 (x^{-\vecdelta+\vecm}\upsilon)
=x^{-\vecdelta+\vecm}
 (m_i/k_i)\lambda \upsilon.
\end{equation}

\subsubsection{The $V$-filtration of $\nbigmtilde$ along $\tau=0$}

Let $\pi:X^{(1)}\lrarr X$ denote the projection.
Let $\pitilde:\nbigx^{(1)}\lrarr\nbigx$
denote the induced morphism.
Set 
$\pi^{\ast}\nbigr_X:=
 \nbigo_{\nbigx^{(1)}}
 \otimes_{\pitilde^{-1}\nbigo_{\nbigx}}
 \pitilde^{-1}\nbigr_X $.
Let $\lefttop{\tau}V_0\nbigr_{X^{(1)}}$
denote the sheaf of subalgebras in 
$\nbigr_{X^{(1)}}$
generated by
$\pi^{\ast}\nbigr_{X}$
and
$\tau\deldel_{\tau}$.
We shall define
$\lefttop{\tau}V_0\nbigr_{X^{(1)}}$-coherent 
submodules
$U_{\alpha}\nbigmtilde$
of $\nbigmtilde$
for any $\alpha\in\real$.
For $0\leq\alpha\leq 1$,
we set 
\[
 U'_{\alpha}\nbigmtilde:=
 \pi^{\ast}\nbigr_X
 \Bigl(
 \nbigo_{\nbigx^{(1)}}(\nbigd^{(1)}+[\alpha \nbigp^{(1)}])\upsilon
 \Bigr)
\subset
 \nbigmtilde.
\]
For any real number $\alpha$,
we take the integer $n$
such that 
$0\leq \alpha+n<1$,
and we set
$U_{\alpha}\nbigmtilde:=
 \tau^{n}U'_{\alpha+n}\nbigmtilde$.

\begin{remark}
We do not need $U'_1(\nbigmtilde)$
for the construction of
$U_{\alpha}\nbigmtilde$.
We use it in the proof of Theorem
{\rm\ref{thm;14.5.2.40}}.
\end{remark}

We shall prove the following theorem 
in \S\ref{subsection;14.5.2.201}--\ref{subsection;14.5.2.202}.
\begin{theorem}
\label{thm;14.5.2.40}
$U_{\bullet}\nbigmtilde$
is the $V$-filtration of $\nbigmtilde$.
\end{theorem}

\begin{remark}
Theorem {\rm\ref{thm;14.5.2.40}}
is not precisely the twistor version of 
Proposition {\rm\ref{prop;14.5.1.10}}
because $\nbigmtilde_{|\{\lambda=1\}}=\nbigm(\ast \tau)$.
\end{remark}

\subsubsection{Easy property of $U_{\bullet}\nbigmtilde$}
\label{subsection;14.5.2.201}

We may prove the following lemma
by the argument in the proof of Lemma \ref{lem;14.5.20.10}.

\begin{lemma}
Outside $\{\tau=0\}$,
we have
$U_{\alpha}\nbigmtilde=\nbigmtilde$.
\qed
\end{lemma}

\begin{lemma}
We have a natural action of
$\lefttop{\tau}V_0\nbigr_{X^{(1)}}$
on $U_{\alpha}\nbigmtilde$.
\end{lemma}
\pf
We have only to check it around any point of 
$\nbigp_{\reduced}^{(1)}$
by using a coordinate system as in
\S\ref{subsection;14.4.23.1},
and the relation (\ref{eq;14.5.2.22}).
The argument is the same as that in
Lemma \ref{lem;14.5.1.3}.
\qed

\begin{lemma}
For $\alpha-1<\beta\leq \alpha$,
we have
$\tau U_{\alpha}\nbigmtilde
\subset
 U_{\beta}\nbigmtilde$.
As a result,
we have 
$U_{\alpha}\nbigmtilde
\subset
 U_{\alpha'}\nbigmtilde$
for any $\alpha\leq\alpha'$.
\end{lemma}
\pf
We have only to prove the claim
around any point of $\nbigp_{\reduced}^{(1)}$
by using a coordinate system in \S\ref{subsection;14.4.23.1}.
Set $\vecp:=[\alpha\veck]$.
It is enough to consider the case $0\leq \alpha<1$.
If $\alpha-1<\beta<0$,
we have 
$\tau U_{\alpha}\nbigmtilde
\subset
 U_{\beta}\nbigmtilde$
by the construction.
By using the relation
$\tau x^{-\vecdelta-\vecp}\upsilon
=\tau\deldel_{\tau}\bigl(
 x^{\veck-\vecdelta-\vecp}\upsilon
 \bigr)$,
we can prove the claim 
in the case $0\leq \beta<\alpha$.
\qed

\vspace{.1in}

The following lemma 
is proved as in the case of the $\nbigd$-modules
(Lemma \ref{lem;14.5.2.30}).

\begin{lemma}
The $\lefttop{\tau}V_0\nbigr_{X^{(1)}}$-module
$U_{\alpha}\nbigmtilde$
is coherent.
\qed
\end{lemma}

We set $X_0:=\{0\}\times X\subset X^{(1)}$.
We set
$U_{<\alpha}\nbigmtilde:=
 \bigcup_{\beta<\alpha}
 U_{\beta}\nbigmtilde$.
We obtain an $\nbigr_{X_0}$-module
$U_{\alpha}\nbigmtilde\big/
 U_{<\alpha}\nbigmtilde$.
It is equipped with an endomorphism
$[\tau\deldel_{\tau}]$
induced by $\tau\deldel_{\tau}$.
The following lemma is 
a generalization of Lemma \ref{lem;14.5.2.31}.

\begin{lemma}
\label{lem;14.5.2.35}
For $0<\alpha<1$,
the action of
$[\tau\deldel_{\tau}]+\alpha\lambda$
on $U_{\alpha}\nbigmtilde\big/U_{<\alpha}\nbigmtilde$
is nilpotent.
Indeed,
the action of
$\bigl([\tau\deldel_{\tau}]
+\alpha\lambda\bigr)^{\dim X}$
is $0$.
\end{lemma}
\pf
We have only to check the claim 
around any point of $\nbigp^{(1)}_{\reduced}$
by using a coordinate system as in \S\ref{subsection;14.4.23.1}.
Set $\vecp:=[\alpha\veck]$.
By using 
$\bigl(
 \tau\deldel_{\tau}+\lambda\alpha
\bigr)x^{-\vecdelta-\vecp}\upsilon
=-k_i^{-1}\deldel_{i}(x_ix^{-\vecdelta-\vecp}\upsilon)$,
we obtain
\[
 (\tau\deldel_{\tau}+\lambda\alpha)^N
 (x^{-\vecdelta-\vecp}\upsilon)
 \in U_{<\alpha}\nbigmtilde
\]
for any $N\geq \ell_1$
as in the proof of Lemma \ref{lem;14.5.2.31}.

We have
$\bigl(
 \tau\deldel_{\tau}+\lambda\alpha
 \bigr)(gx^{-\vecdelta-\vecp}\upsilon)
=g\cdot(\tau\deldel_{\tau}+\lambda\alpha)
(x^{-\vecdelta-\vecp}\upsilon)
+(\tau\deldel_{\tau}g)x^{-\vecdelta-\vecp}\upsilon$.
We have
$(\tau\deldel_{\tau}g)x^{-\vecdelta-\vecp}\upsilon
\in\tau U_{\alpha}\nbigmtilde
\subset U_{<\alpha}\nbigmtilde$.
Then, the claim of the lemma follows.
\qed

\vspace{.1in}
The following lemma can be also checked
as in the case of the $\nbigd$-modules
(Lemma \ref{lem;14.5.2.33}).

\begin{lemma}
\label{lem;14.5.2.34}
If $N\geq \dim X+1$,
we have 
$(\tau\deldel_{\tau})^NU_0\nbigmtilde
\subset
 \tau U_{<1}\nbigmtilde$.
\end{lemma}
\pf
We have only to check the claim 
around any point of $\nbigp^{(1)}_{\reduced}$
by using a coordinate system as in \S\ref{subsection;14.4.23.1}.
Let $\vecdelta_1$ be as in the proof of
Lemma \ref{lem;14.5.2.33}.
We have 
$(\tau\deldel_{\tau})^{\ell_1}(x^{-\vecdelta}\upsilon)
=\prod_{i=1}^{\ell_1}
 \bigl(-k_i^{-1}\deldel_i\bigr)\cdot x^{-(\vecdelta-\vecdelta_1)}\upsilon$.
Hence, we have
\begin{multline}
 (\tau\deldel_{\tau})^{\ell_1+1}
 (x^{-\vecdelta}\upsilon)
=\prod_{i=1}^{\ell_1}(-k_i^{-1}\deldel_i)
 (\tau f x^{-(\vecdelta-\vecdelta_1)}\upsilon)
 \\
=\tau
 \prod_{i=1}^{\ell_1}(-k_i^{-1}\deldel_i)
 (x^{-\vecdelta-(\veck-\vecdelta_1)}\upsilon)
\in \tau U_{<1}\nbigmtilde.
\end{multline}

For any section $s$ of $U_0\nbigmtilde$
and any holomorphic function $g$,
we have
$\tau\deldel_{\tau}(gs)
=\tau(\deldel_{\tau}g)\,s
+g\,\tau\deldel_{\tau}s$,
and $\tau(\deldel_{\tau}g)s\in \tau U_{<1}(\nbigmtilde)$.
Then, we obtain the claim of the lemma.
\qed

\vspace{.1in}

We obtain the following
from Lemma \ref{lem;14.5.2.35}
and Lemma \ref{lem;14.5.2.34}.
\begin{lemma}
$[\tau\deldel_{\tau}]+\alpha\lambda$
is nilpotent 
on $U_{\alpha}\nbigmtilde/U_{<\alpha}\nbigmtilde$
for any $\alpha\in\real$.
\qed
\end{lemma}

\begin{lemma}
We have
$\bigcup_{\alpha\in\real}
 U_{\alpha}\nbigmtilde
=\nbigmtilde$.
\end{lemma}
\pf
We have only to prove the claim
locally around any point of $\nbigd^{(1)}$.
By the construction,
$\nbigmtilde':=\bigcup_{\alpha\in\real}
 U_{\alpha}\nbigmtilde$
is an $\nbigr_{X^{(1)}}(\ast \tau)$-module.
By the construction,
we have
$U_{\alpha}\nbigmtilde
\subset
 \nbigmtilde$
for any $\alpha$,
and hence
$\nbigmtilde'\subset\nbigmtilde$.
We have
$\deldel_{\tau}^N
 x^{-\vecdelta-\vecp}
 \upsilon
=x^{-\vecdelta-\vecp-N\veck}\upsilon$.
Hence, we have
$\nbigq\nbige(\nbigh^{(1)})(\ast\tau)
\subset \nbigmtilde'$.
We obtain 
$\nbigmtilde\subset\nbigmtilde'$.
\qed

\vspace{.1in}

To prove Theorem \ref{thm;14.5.2.40},
it remains to prove that 
$U_{\alpha}\nbigmtilde\big/U_{<\alpha}\nbigmtilde$
are strict for any $\alpha$,
i.e.,
we need to establish that
the multiplication of $\lambda-\lambda_1$
on $U_{\alpha}\nbigmtilde\big/U_{<\alpha}\nbigmtilde$
is injective for any $\lambda_1\in\cnum$.
We shall prove it in \S\ref{subsection;14.5.2.202}
after some preparations.

\subsubsection{Preliminary}
\label{subsection;14.5.15.2}

Let $(U,x_1,\ldots,x_n)$ and $Y$
be as in \S\ref{subsection;14.5.1.20}.
We set
$\nbigy:=\cnum_{\lambda}\times Y$
and $\nbigy^{\lambda_0}:=\{\lambda_0\}\times Y$.
In the following,
$\nbigy^{(\lambda_0)}$
denotes a neighbourhood of
$\nbigy^{\lambda_0}$
in $\nbigy$,
and $\nbign^{(\lambda_0)}$
denote the product of
$\nbigy^{(\lambda_0)}$
and 
$\bigl\{(x_1,\ldots,x_{\ell})\,\big|\,|x_i|<\epsilon\bigr\}$
for some $\epsilon>0$.

For any section $s$ of
$\nbigo_{\nbigx^{(1)}}(\ast \nbigd^{(1)})$
on $\nbign^{(\lambda_0)}$,
we have the Laurent expansion
\[
 s=\sum_{\vecm\in\seisuu^{\ell}}
 \sum_{j\in\seisuu_{\geq 0}}
 h_{\vecm,j}x^{\vecm}\tau^j,
\]
where $h_{\vecm,j}$ are holomorphic functions
on $\nbigy^{(\lambda_0)}$.
As before,
we have the unique expansion
\[
 s=\sum_{\vecm\in\seisuu^{\ell}}
 \sum_{j\in\seisuu_{\geq 0}}
 h^{(1)}_{\vecm,j}x^{\vecm}(\tau f)^j,
\]
where $h^{(1)}_{\vecm,j}$ are holomorphic functions
on $\nbigy^{(\lambda_0)}$.
Indeed, we have
$h^{(1)}_{\vecm,j}
=h_{\vecm-j\veck,j}$.

\begin{definition}
A section $s$ of $\nbigo_{\nbigx^{(1)}}(\ast \nbigd^{(1)})$
on $\nbign^{(\lambda_0)}$
is called $(\vecm,j)$-primitive
if it is expressed as
$s=gx^{\vecm}(\tau f)^j$
for a holomorphic function $g$ on $\nbign^{(\lambda_0)}$
with $g_{|\nbigy^{(\lambda_0)}}\neq 0$.
\qed
\end{definition}

\begin{definition}
A primitive expression of 
a section $s$ of $\nbigo_{\nbigx^{(1)}}(\ast\nbigd^{(1)})$
on $\nbign^{(\lambda_0)}$
is an expression 
\[
 s=\sum_{(\vecm,j)\in\nbigs}
 g_{\vecm,j}x^{\vecm}(\tau f)^{j},
\]
where $\nbigs$ is a finite subset
in $\seisuu^{\ell}\times\seisuu_{\geq 0}$
and $g_{\vecm,j}$ are holomorphic functions
on $\nbign^{(\lambda_0)}$
with $g_{\vecm,j|\nbigy^{(\lambda_0)}}\neq 0$.
\qed
\end{definition}

We use the partial orders on $\seisuu^{\ell}$
and $\seisuu^{\ell}\times\seisuu_{\geq 0}$
as in \S\ref{subsection;14.5.1.20}.
We reword Lemma \ref{lem;14.5.1.21}
and Corollary \ref{cor;14.5.2.50}
in this context.

\begin{lemma}
\label{lem;14.5.2.131}
Any section $s$ of
$\nbigo_{\nbigx^{(1)}}(\ast\nbigd^{(1)})$
on $\nbign^{(\lambda_0)}$
has a primitive expression
$s=\sum_{(\vecm,j)\in\nbigs}
 g_{\vecm,j}x^{\vecm}(\tau f)^{j}$.
Moreover the following claims hold.
\begin{itemize}
\item
The set $\min\pi(\nbigs)$
is well defined for $s$.
For any $\vecm\in\min\pi(\nbigs)$,
$g_{\vecm,j|\nbigy^{(\lambda_0)}}$
is well defined for $s$.
\item
The set $\min(\nbigs)$ is well defined for $s$.
For any $(\vecm,j)\in\min(\nbigs)$,
$g_{\vecm,j|\nbigy^{(\lambda_0)}}$
is well defined for $s$.
\qed
\end{itemize}
\end{lemma}

We reword Lemma \ref{lem;14.5.15.1}
in this context.
\begin{lemma}
Let $\nbigt\subset\seisuu^{\ell}\times\seisuu_{\geq 0}$
be a finite subset.
Suppose 
we are given holomorphic functions
$g_{\vecm,j}$ $((\vecm,j)\in\nbigt)$
on $\nbign^{(\lambda_0)}$
such that 
$\sum_{(\vecm,j)\in\nbigt}
 g_{\vecm,j}x^{\vecm}(\tau f)^j=0$.
\begin{itemize}
\item
For any $\vecm\in\min\pi(\nbigt)$
and any $j\in\seisuu_{\geq 0}$,
we have $g_{\vecm,j|\nbigy^{(\lambda_0)}}=0$.
\item
For any $(\vecm,j)\in\min(\nbigt)$,
we have $g_{\vecm,j|\nbigy^{(\lambda_0)}}=0$.
\qed
\end{itemize}
\end{lemma}

\subsubsection{Primitive expressions}

Let $\Ptilde$ be any effective divisor of $X$
such that $\Ptilde_{\reduced}\subset P_{\reduced}$.
We set $\Ptilde^{(1)}:=\cnum_{\tau}\times \Ptilde$
and $\nbigptilde^{(1)}:=\cnum_{\lambda}\times \Ptilde^{(1)}$.
For any non-negative integer $N$,
we set
\[
 G_N\nbigmtilde_{\Ptilde}:=
 \sum_{j=0}^N
 \nbigo_{\nbigx^{(1)}}
 \bigl(\nbigd^{(1)}+\nbigptilde^{(1)}\bigr)
 (\tau f)^j\upsilon
\subset
\nbigmtilde.
\]
Let 
$\nbigmtilde_{\Ptilde}$ be the sheafification of
$\bigcup_{N}G_N\nbigmtilde_{\Ptilde}$.
We set $U_{\Ptilde}\nbigmtilde:=
 \pi^{\ast}\nbigr_X\cdot
 \nbigmtilde_{\Ptilde}$
in $\nbigmtilde$.

\vspace{.1in}

Let $V_0\nbigr_X\subset\nbigr_X$
be generated by
$\lambda\Theta_{\nbigx/\cnum}(\log \nbigd)$
over $\nbigo_{\nbigx}$.
Set
$\pi^{\ast}V_0\nbigr_X:=
 \nbigo_{\nbigx^{(1)}}
 \otimes_{\pitilde^{-1}\nbigo_{\nbigx}}
 \pitilde^{-1}V_0\nbigr_X$.
We naturally regard
$\pi^{\ast}V_0\nbigr_X\subset
 \nbigr_{X^{(1)}}$.
We can check the following lemma
as in the case of the $\nbigd$-modules
(Lemma \ref{lem;14.5.2.100}).
\begin{lemma}
We have $\nbigmtilde_{\Ptilde}
=\pi^{\ast}V_0\nbigr_X
 \cdot\Bigl(
 \nbigo_{\nbigx^{(1)}}
 \bigl(\nbigd^{(1)}+\nbigptilde^{(1)}\bigr)
 \upsilon
 \Bigr)$,
and hence
\[
 U_{\Ptilde}\nbigmtilde=
 \pi^{\ast}\nbigr_X\cdot
 \Bigl(
 \nbigo_{\nbigx^{(1)}}
 \bigl(\nbigd^{(1)}+\nbigptilde^{(1)}\bigr)
 \upsilon
 \Bigr). 
\]
\qed
\end{lemma}
By the lemma, we have
$U_{[\alpha P]}\nbigmtilde
=U'_{\alpha}\nbigmtilde$
for any $0\leq\alpha\leq 1$.

\vspace{.1in}

We use the notation in \S\ref{subsection;14.5.15.2}.
Take any $\lambda_0\in\cnum$.
We set $D_i=\{z_i=0\}$ on the neighbourhood $U$.
Let $\vecp\in\seisuu_{\geq 0}^{\ell}$ be determined by
$\sum_{i=1}^{\ell} p_iD_i=\Ptilde$ on $U$.
Note we have $p_i=0$ for $i>\ell_1$.
Let $s$ be a non-zero section of
$U_{\Ptilde}\nbigmtilde$
on $\nbign^{(\lambda_0)}$.
We have an expression
$s=\sum_{\vecn\in\seisuu_{\geq 0}^{\ell}}
 \deldel^{\vecn}s_{\vecn}$
as an essentially finite sum,
where $s_{\vecn}$ are
sections of $\nbigmtilde_{\Ptilde}$
on $\nbign^{(\lambda_0)}$.

\begin{definition}
A section $s$ of $U_{\Ptilde}\nbigmtilde$
on $\nbign^{(\lambda_0)}$ 
is called 
$(\vecm,j)$-primitive
if it is expressed as 
\[
 s=\deldel^{\vecm_-}\bigl(
 g x^{-\vecdelta-\vecp+\vecm_+}(\tau f)^j\upsilon\bigr)
\]
for a holomorphic function $g$ on $\nbign^{(\lambda_0)}$
with $g_{|\nbigy^{(\lambda_0)}}\neq 0$.
\qed
\end{definition}

\begin{definition}
A primitive expression of a section $s$ of 
$U_{\Ptilde}\nbigmtilde$
on $\nbign^{(\lambda_0)}$ is 
a decomposition
\[
 s=\sum_{(\vecm,j)\in\nbigs} s_{(\vecm,j)},
\]
where $\nbigs$ is a finite subset of
$\seisuu^{\ell}\times\seisuu_{\geq 0}$,
and $s_{(\vecm,j)}$ are $(\vecm,j)$-primitive 
sections of $U_{\Ptilde}\nbigmtilde$.
\qed
\end{definition}

\begin{lemma}
\label{lem;14.5.2.120}
Any section $s$ of
$U_{\Ptilde}\nbigmtilde$
on $\nbign^{(\lambda_0)}$
has a primitive expression
\[
 s=\sum_{(\vecm,j)\in\nbigs}
 \deldel^{\vecm_-}\bigl(
 g_{\vecm,j}x^{-\vecdelta-\vecp+\vecm_+}
 (\tau f)^j\upsilon
\bigr).
\]
\end{lemma}
\pf
We give an algorithm to obtain a primitive expression.
Let $s$ be any section of
$U_{\Ptilde}\nbigmtilde$
on $\nbign^{(\lambda_0)}$.
We have the following expression
as an essentially finite sum:
\begin{equation}
\label{eq;14.5.14.31}
 s=\sum_{\vecn\in\seisuu^{\ell}_{\geq 0}}
 \sum_{\vecq\in\seisuu_{\geq 0}^{\ell}}
 \sum_{j\in\seisuu_{\geq 0}}
 \deldel^{\vecn}
 \Bigl(
 g_{\vecn,\vecq,j}
 x^{-\vecdelta-\vecp+\vecq}
 (\tau f)^j\upsilon
 \Bigr).
\end{equation}
For $a\in\seisuu_{\geq 0}$,
we say that 
an expression (\ref{eq;14.5.14.31})
has the property $(R_a)$
if the following holds
for any $(\vecn,\vecq,j)$
with $|\vecn|\geq a$:
\begin{itemize}
\item
We have
$g_{\vecn,\vecq,j}=0$
unless
$\{i\,|\,n_i\neq 0,q_i\neq 0\}=\emptyset$.
\item
If $g_{\vecn,\vecq,j}\neq 0$,
then
$\deldel^{\vecn}\bigl(
 g_{\vecn,\vecq,j}
 x^{-\vecdelta-\vecp+\vecq}
 (\tau f)^j\upsilon
 \bigr)$
are $(\vecq-\vecn,j)$-primitive.
\end{itemize}
Take any expression
$s=
 \sum_{\vecn,\vecq,j}
 \deldel^{\vecn}
 \bigl(
 g_{\vecn,\vecq,j}
 x^{-\vecdelta-\vecp+\vecq}
 (\tau f)^j\upsilon
 \bigr)$.
If $a$ is sufficiently large,
the expression has the property $(R_a)$.
In general, if $q_i>0$,
we have 
\begin{multline}
\label{eq;14.5.2.130}
 \deldel_ix_i\bigl(
 g
 x^{-\vecdelta-\vecp+\vecq}x_i^{-1}
 (\tau f)^j\upsilon
 \bigr)
=
\Bigl(
x_i\deldel_ig
-\lambda(p_i-q_i+1+jk_i) g
 \Bigr)
 x^{-\vecdelta-\vecp+\vecq}x_i^{-1}
 (\tau f)^j\upsilon
 \\
-k_ig
 x^{-\vecdelta-\vecp+\vecq}x_i^{-1}
 (\tau f)^{j+1}\upsilon.
\end{multline}
Let $s=
 \sum_{\vecn,\vecq,j}
 \deldel^{\vecn}
 \bigl(
 g_{\vecn,\vecq,j}
 x^{-\vecdelta-\vecp+\vecq}
 (\tau f)^j\upsilon
 \bigr)$ be an expression 
with the property $(R_a)$
such that $a\geq 1$.
Applying (\ref{eq;14.5.2.130})
and Lemma \ref{lem;14.5.2.131}
to each
$\deldel^{\vecn}
 \bigl(
 g_{\vecn,\vecq,j}
 x^{-\vecdelta-\vecp+\vecq}
 (\tau f)^j\upsilon\bigr)$
with $|\vecn|=a-1$,
we can obtain an expression 
with the property $(R_{a-1})$.
We can arrive at an expression
with the property $(R_0)$,
which is a primitive expression of $s$.
\qed

\vspace{.1in}

Suppose that we are given
a finite set $\nbigs\subset\seisuu^{\ell}\times\seisuu_{\geq 0}$
and sections of 
$g_{\vecm,j}$ of
$\nbigo_{\nbigx^{(1)}}$
on $\nbign^{(\lambda_0)}$
for $(\vecm,j)\in\nbigs$,
such that
$0=\sum_{(\vecm,j)\in\nbigs}
 \deldel^{\vecm_-}
 \Bigl(
 g_{\vecm,j}x^{-\vecdelta-\vecp+\vecm_+}
 (\tau f)^j\upsilon
 \Bigr)$,
where $\vecm=\vecm_+-\vecm_-$
is the decomposition as in \S\ref{subsection;17.1.19.1}.

\begin{lemma}
\label{lem;14.5.14.100}
For any $\vecm\in\min\pi(\nbigs)$
and $j\in\seisuu_{\geq 0}$,
we have
$g_{\vecm,j|\nbigy^{(\lambda_0)}}=0$.
For any $(\vecm,j)\in\min\nbigs$,
we have
$g_{\vecm,j|\nbigy^{(\lambda_0)}}=0$.
\end{lemma}
\pf
In general, we have the following
for any section $g$ of $\nbigo_{\nbigx^{(1)}}$
on $\nbign^{(\lambda_0)}$:
\begin{multline}
\label{eq;14.5.15.40}
 \deldel_i\bigl(
 gx^{-\vecdelta-\vecp+\vecn}(\tau f)^j\upsilon\bigr)
=
\Bigl(
x_i\deldel_ig
-(1+p_i-n_i+jk_i)\lambda g
\Bigr)
x^{-\vecdelta-\vecp+\vecn}x_i^{-1}(\tau f)^j\upsilon
 \\
-k_igx^{-\vecdelta-\vecp+\vecn}x_i^{-1}
 (\tau f)^{j+1}\upsilon.
\end{multline}
Hence, we have the following expression:
\[
 \deldel^{\vecm_-}
 \bigl(g_{\vecm,j}x^{-\vecdelta-\vecp}x^{\vecm_+}
 (\tau f)^j\upsilon
 \bigr)
=\sum_{0\leq k\leq |\vecm_-|}
 h_{\vecm,j,k}x^{-\vecdelta-\vecp+\vecm}
 (\tau f)^{j+k}\upsilon,
\]
where
$h_{\vecm,j,k}$ are sections of
$\nbigo_{\nbigx^{(1)}}$ on $\nbign^{(\lambda_0)}$
such that 
\[
 h_{\vecm,j,k|\nbigy^{(\lambda_0)}}=
 C_{\vecm,j,k}\cdot g_{\vecm,j|\nbigy^{(\lambda_0)}}
\]
for some $C_{\vecm,j,k}\in\rnum$.
By 
$\{i\,|\,m_{+,i}\neq 0,m_{-,i}\neq 0\}=\emptyset$,
we have $C_{\vecm,j,0}\neq 0$.
We have the following 
in $\nbigo_{\nbigx^{(1)}}(\ast \nbigd^{(1)})$:
\[
0=\sum_{(\vecm,j)\in\nbigs}
 \sum_{0\leq k\leq |\vecm_-|}
 h_{\vecm,j,k}x^{-\vecdelta-\vecp+\vecm}
 (\tau f)^{j+k}.
\]
For $\vecm\in\min\pi(\nbigs)$
and $p\geq 0$,
we have
$\sum_{j+k=p}
 C_{\vecm,j,k}g_{\vecm,j|\nbigy^{(\lambda_0)}}=0$.
We obtain
$g_{\vecm,j|\nbigy^{(\lambda_0)}}=0$ by an ascending induction on $j$.
For $(\vecm,j)\in\min\nbigs$,
we have
$C_{\vecm,j,0}g_{\vecm,j|\nbigy^{(\lambda_0)}}=0$.
Thus, the proof of Lemma \ref{lem;14.5.14.100}
is finished.
\qed

\begin{corollary}
\label{cor;14.5.14.30}
Let $s$ be a section of
$U_{\Ptilde}\nbigmtilde$
on $\nbign^{(\lambda_0)}$
with a primitive expression
\[
 s=\sum_{(\vecm,j)\in\nbigs}
 \deldel^{\vecm_-}\bigl(
 g_{\vecn,j}x^{-\vecdelta-\vecp+\vecm_+}(\tau f)^j\upsilon
 \bigr).
\]
\begin{itemize}
\item
The set
$\min\pi(\nbigs)$ is well defined for $s$.
For any $\vecm\in\min\pi(\nbigs)$
and any $j\in\seisuu_{\geq 0}$,
$g_{\vecm,j|\nbigy^{(\lambda_0)}}$
is well defined for $s$.
\item
 The set
$\min\nbigs$ is well defined for $s$.
For any $(\vecm,j)\in\min(\nbigs)$,
$g_{\vecm,j|\nbigy^{(\lambda_0)}}$
is well defined for $s$.
\qed
\end{itemize}
\end{corollary}

\subsubsection{Variant of primitive expression}

Let $P_1$ be an irreducible component of $P_{\reduced}$.
We have the inclusion
$U_{\Ptilde}\nbigmtilde
\subset
 U_{\Ptilde+P_1}\nbigmtilde$.
We assume that $P_1=\{x_1=0\}$
on $(U,x_1,\ldots,x_n)$.
\begin{lemma}
\label{lem;14.5.15.20}
Any section $s$ of $U_{\Ptilde+P_1}\nbigmtilde$
on $\nbign^{(\lambda_0)}$
has an expression as an essentially finite sum
\[
 s=s'+\sum_{\vecn\in\seisuu^{\ell}}
 \deldel^{\vecn_-}\bigl(
 g_{\vecn}(\lambda,\tau,x)
 x^{-\vecdelta-\vecp+\vecn_+}x_1^{-1}
 \upsilon
 \bigr),
\]
such that
(i) $s'$ is a section of $U_{\Ptilde}\nbigmtilde$,
(ii) each 
$\deldel^{\vecn_-}\bigl(
 g_{\vecn}(\lambda,\tau,x)
 x^{-\vecdelta-\vecp+\vecn_+}x_1^{-1}
 \upsilon
 \bigr)$
is $(\vecn,0)$-primitive
as a section of $U_{\Ptilde+P_1}\nbigmtilde$.
\end{lemma}
\pf
We prepare two procedures
which are used to obtain an expression with the desired property.

\paragraph{(A)}
Suppose that 
a section $s$ of $U_{\Ptilde+P_1}\nbigm$
on $\nbign^{(\lambda_0)}$
has an expression as an essentially finite sum
\begin{equation}
\label{eq;14.5.15.41}
 s=\sum_{\vecn\in\seisuu^{\ell}_{\geq 0}}
 \sum_{j\geq 0}
 \deldel^{\vecn}\bigl(
 g_{\vecn,j}(\lambda,\tau,x)
 x^{-\vecdelta-\vecp}x_1^{-1}
 (\tau f)^j\upsilon
 \bigr).
\end{equation}
Then, 
by using the relations (\ref{eq;14.5.15.40}),
we obtain an expression of $s$
as an essentially finite sum:
\begin{equation}
\label{eq;14.5.15.42}
 s=s'+\sum_{\vecn\in\seisuu^{\ell}_{\geq 0}}
 \deldel^{\vecn}\bigl(
 g^{(1)}_{\vecn}(\lambda,\tau,x)
 x^{-\vecdelta-\vecp}x_1^{-1}\upsilon
 \bigr).
\end{equation}
Here, 
(i) $s'$ is a section of $U_{\Ptilde}\nbigm$,
(ii) $\max\bigl\{
 |\vecn|\,\big|\,
 g^{(1)}_{\vecn}\neq 0
 \bigr\}\leq 
\max\bigl\{
 |\vecn|\,\big|\,
 \exists j\,\,
 g_{\vecn,j}\neq 0
 \bigr\}$.

\paragraph{(B)}
Suppose that we are given a section
$\deldel^{\vecn}\bigl(
 g(\lambda,\tau,x)x^{-\vecdelta-\vecp}x_1^{-1}\upsilon
 \bigr)$.
By applying Lemma \ref{lem;14.5.2.131}
and (\ref{eq;14.5.2.130}),
we obtain an expression of 
$\deldel^{\vecn}\bigl(
 g(\lambda,\tau,x)x^{-\vecdelta-\vecp}x_1^{-1}\upsilon
 \bigr)$
as an essentially finite sum:
\begin{multline}
 \deldel^{\vecn}\bigl(
 g(\lambda,\tau,x)x^{-\vecdelta-\vecp}x_1^{-1}
 \upsilon
 \bigr)
=s'
+\sum_{\vecq}
 \deldel^{\vecn}(g^{(2)}_{\vecq}(\lambda,\tau,x)
 x^{-\vecdelta-\vecp+\vecq}x_1^{-1}\upsilon)
 \\
+\sum_{\vecm\in\seisuu^{\ell}_{\geq 0}}
 \sum_{j\geq 0}
 \deldel^{\vecm}\bigl(
 g^{(3)}_{\vecm,j}x^{-\vecdelta-\vecp}x_1^{-1}
 (\tau f)^j\upsilon
 \bigr).
\end{multline}
Here,
(i) $s'$ is a section of $U_{\Ptilde}\nbigm$,
(ii) 
if $g^{(2)}_{\vecq}\neq 0$,
$\{i\,|\,n_i\neq 0,q_i\neq 0\}$ is empty
and
$\deldel^{\vecn}\bigl(g^{(2)}_{\vecq}(\lambda,\tau,x')
 x^{-\vecdelta-\vecp+\vecq}x_1^{-1}\upsilon\bigr)$
is $(\vecq-\vecn,0)$-primitive,
(iii) if $g^{(3)}_{\vecm,j}\neq 0$,
we have $|\vecm|<|\vecn|$.

\vspace{.1in}
Suppose that we have an expression 
(\ref{eq;14.5.15.41}) of $s$
such that $|\vecn|\leq a$
for any $(\vecn,j)$
such that $g_{\vecn,j}\neq 0$.
By applying {\bf(A)},
we obtain an expression (\ref{eq;14.5.15.42}).
We have
$|\vecn|\leq a$
if $g^{(1)}_{\vecn}\neq 0$.
By applying {\bf(B)}
to $\deldel^{\vecn}\bigl(
 g^{(1)}_{\vecn}(\lambda,\tau,x) 
 x^{-\vecdelta-\vecp}x_1^{-1}\upsilon
 \bigr)$,
we obtain an expression of $s$
as follows:
\begin{multline}
 s=
s''+
 \sum_{|\vecn|=a}\sum_{\vecq}
 \deldel^{\vecn}\bigl(
 g^{(4)}_{\vecn}(\lambda,\tau,x)
 x^{-\vecdelta-\vecp+\vecq}x_1^{-1}\upsilon
 \bigr)
 \\
+\sum_{|\vecn|<a}\sum_j
 \deldel^{\vecn}\bigl(
 g^{(5)}_{\vecn,j}
 x^{-\vecdelta-\vecp}x_1^{-1}
 (\tau f)^j\upsilon
 \bigr).
\end{multline}
Here, 
(i) $s''$ is a section of $U_{\Ptilde}\nbigm$,
(ii) if $g^{(4)}_{\vecn}\neq 0$,
$\{i\,|\,n_i\neq 0,q_i\neq 0\}$ is empty
and 
$\deldel^{\vecn}\bigl(
 g^{(4)}_{\vecn}(\lambda,\tau,x)
 x^{-\vecdelta-\vecp+\vecq}x_1^{-1}\upsilon
 \bigr)$
is $(\vecq-\vecn,0)$-primitive.
Then, we can obtain an expression of 
$s$ with the desired property
by an inductive argument.
\qed

\begin{lemma}
Let $\nbigt_1$ be a finite subset in $\seisuu^{\ell}$
such that $m_1\leq 0$ for any $\vecm\in\nbigt_1$.
Suppose that we are given 
a section $s'$ of $U_{\Ptilde}\nbigm$
and holomorphic functions
$g_{\vecm}(\lambda,\tau,x)$ 
$(\vecm\in\nbigt_1)$
on $\nbign^{(\lambda_0)}$
such that the following holds:
\begin{equation}
\label{eq;14.5.15.50}
 s'=
 \sum_{\vecm\in\nbigt_1}
 \deldel^{\vecm_-}\bigl(
 g_{\vecm}(\lambda,\tau,x)
 x^{-\vecdelta-\vecp+\vecm_+}
 x_1^{-1}\upsilon
 \bigr).
\end{equation}
Then, 
$g_{\vecm|\nbigy^{(\lambda_0)}}=0$
for any $\vecm\in\min\nbigt_1$.
\end{lemma}
\pf
Take a primitive expression of $s'$
in $U_{\Ptilde}\nbigmtilde$:
\[
 s'=\sum_{\vecn\in\nbigt_2}\sum_{j=0}^{j_0(\vecn)}
 \deldel^{\vecn_-}
 \bigl(h_{\vecn,j}x^{-\vecdelta-\vecp+\vecn_+}
 (\tau f)^j\upsilon\bigr).
\]
Here, $\nbigt_2$ be a finite subset of $\seisuu^{\ell}$.
We set
$\nbigt_{2+}:=\bigl\{
 \vecn\in\nbigt_2\,\big|\,
 n_1\geq 0
 \bigr\}$
and 
$\nbigt_{2-}:=\bigl\{
 \vecn\in\nbigt_2\,\big|\,
 n_1< 0
 \bigr\}$.
For $\vecn$,
we set $\vecn_-':=\vecn_--(1,0,\ldots,0)$
and $\vecn_+':=\vecn_++(1,0,\ldots,0)$.
Then, 
by using (\ref{eq;14.5.2.130}),
we obtain an expression as follows:
\begin{multline}
 \label{eq;14.5.15.10}
 s'=
\sum_{\vecn\in\nbigt_{2+}}
\sum_{j=0}^{j_0(\vecn)}
 \deldel^{\vecn_-}\bigl(
 h_{\vecn,j}x^{-\vecdelta-\vecp+\vecn'_+}x_i^{-1}
 (\tau f)^j
 \upsilon
 \bigr)
 \\
+\sum_{\vecn\in\nbigt_{2-}}
 \sum_{j=0}^{j_0(\vecn)+1}
 \deldel^{\vecn'_-}
\bigl(
 h^{(1)}_{\vecn,j}x^{-\vecdelta-\vecp+\vecn_+}x_1^{-1}
 (\tau f)^j\upsilon
 \bigr).
\end{multline}
Here, we can observe that
$\deldel^{\vecn_-'}
\bigl(
 h^{(1)}_{\vecn,j_0(\vecn)+1}
 x^{-\vecdelta-\vecp+\vecn_+}x_1^{-1}
 (\tau f)^{j_0(\vecn)+1}\upsilon
 \bigr)$
is primitive.

We set 
$\nbigt'_2:=
 \{\vecn+(1,0,\ldots,0)\,|\,\vecn\in\nbigt_2\}$
and 
$\nbigt'_{2\pm}:=
 \{\vecn+(1,0,\ldots,0)\,|\,\vecn\in\nbigt_{2\pm}\}$.
Take any
$\vecm\in
 \min\bigl(\nbigt_1\cup\nbigt_2'\bigr)$.
If $\vecm\in\nbigt'_{2-}$,
i.e., $\vecm=\vecn+(1,0,\ldots,0)$
for $\vecn\in\nbigt_{2-}$,
we obtain
$h^{(1)}_{\vecn,j|\nbigy^{(\lambda_0)}}=0$
for any $j\in\seisuu_{> 0}$
by using (\ref{eq;14.5.15.50}),  (\ref{eq;14.5.15.10})
and Lemma \ref{lem;14.5.14.100}.
But, it contradicts with
$h^{(1)}_{\vecn,j_0(\vecn)+1|\nbigy^{(\lambda_0)}}\neq 0$.
Hence, we obtain 
$\vecm\not\in\nbigt_{2-}'$.
Then, we obtain 
$\min(\nbigt_1)
\subset
 \min\bigl(\nbigt_1\cup\nbigt_2'\bigr)$.
Then, we obtain the claim of the lemma
by using 
(\ref{eq;14.5.15.50}),  (\ref{eq;14.5.15.10})
and Lemma \ref{lem;14.5.14.100}
again.
\qed

\begin{corollary}
\label{cor;14.5.15.100}
Let $s$ be a section of 
$U_{\Ptilde+P_1}\nbigmtilde$
with an expression 
as in Lemma {\rm \ref{lem;14.5.15.20}}:
\begin{equation}
\label{eq;14.5.15.101}
 s=s'+\sum_{\vecn\in\nbigt}
 \deldel^{\vecn_-}\bigl(
 g_{\vecn}(\lambda,\tau,x)
 x^{-\vecdelta-\vecp+\vecn_+}x_1^{-1}
 \upsilon
 \bigr).
\end{equation}
Then, the set $\min(\nbigt)$
and the functions $g_{\vecn|\nbigy^{(\lambda_0)}}$
$(\vecn\in\min(\nbigt))$
are well defined for $s$.
\qed
\end{corollary}

We obtain the following corollary.
\begin{corollary}
\label{cor;14.5.15.30}
$U_{\Ptilde+P_1}\nbigmtilde\big/
 U_{\Ptilde}\nbigmtilde$
is strict.
\end{corollary}
\pf
We have only to consider the issue locally around
any point of $\nbigp^{(1)}_{\reduced}$
by using the notation above.
Let $s$ be a section of
$U_{\Ptilde+P_1}\nbigmtilde$
on $\nbign^{(\lambda_0)}$.
Take an expression (\ref{eq;14.5.15.101}) of $s$
as in Lemma {\rm \ref{lem;14.5.15.20}}.
By Corollary \ref{cor;14.5.15.100},
$s$ is a section of 
$U_{\Ptilde}\nbigmtilde$
if and only if $\nbigt=\emptyset$.
Then, for any $\lambda_1\in\cnum$,
the conditions for $s$
and $(\lambda-\lambda_1)s$
are equivalent.
\qed

\vspace{.1in}

Note that we have
$\tau U_{\Ptilde+P}
\subset
 U_{\Ptilde}$.

\begin{lemma}
Let $s$ be a section of
$U_{\Ptilde}\nbigmtilde$
on $\nbign^{(\lambda_0)}$
with a primitive expression
$s=\sum_{(\vecm,j)\in\nbigs}s_{\vecm,j}$.
It is a section of
$\tau U_{\Ptilde+P}\nbigmtilde$
if and only if
$j\neq 0$ for any 
$(\vecm,j)\in\nbigs$.
\end{lemma}
\pf
Note that the condition is independent 
of the choice of a primitive expression,
by Corollary \ref{cor;14.5.14.30}.
Then, the claim is clear.
\qed

\begin{corollary}
\label{cor;14.5.15.31}
$U_{\Ptilde}\nbigmtilde\big/
\tau U_{\Ptilde+P}\nbigmtilde$
is strict.
\qed
\end{corollary}

\subsubsection{The end of the proof of 
Theorem \ref{thm;14.5.2.40}}
\label{subsection;14.5.2.202}

Let us finish the proof of 
Theorem \ref{thm;14.5.2.40}.
It remains to prove that
$U_{\alpha}\nbigmtilde/U_{<\alpha}\nbigmtilde$
is strict
for any $\alpha\in\real$.
We have only to prove the claim
locally around any point of $\nbigd^{(1)}$
by using a coordinate system as in 
\S\ref{subsection;14.4.23.1}.
It is enough to consider the case $0\leq\alpha<1$.
Take any $\lambda_0\in\cnum_{\lambda}$.
Note that $U_{[\alpha P]}\nbigmtilde=U'_{\alpha}\nbigmtilde$
for $0\leq \alpha\leq 1$.
We have an effective divisor $P(<\alpha)$  for $0<\alpha\leq 1$
such that
(i) $U_{P(<\alpha)}\nbigmtilde=U'_{<\alpha}\nbigmtilde$,
(ii) $[\alpha P]-P(<\alpha)$ is effective.
Hence, by using Corollary \ref{cor;14.5.15.30} successively,
we obtain that
$U'_{\alpha}\nbigmtilde\big/U'_{<\alpha}\nbigmtilde$
is strict for $0<\alpha<1$.

We have
$U_{<0}\nbigmtilde=\tau U_{<1}\nbigmtilde
\subset
 \tau U'_1\nbigmtilde
\subset U_0\nbigmtilde$.
We obtain that
$(\tau U'_1\nbigmtilde)/U_{<0}\nbigmtilde$ is strict
as above.
By using Corollary \ref{cor;14.5.15.31},
$U_0\nbigmtilde\big/(\tau U'_1\nbigmtilde)$ is strict.
Thus, we are done.
\qed

\subsection{Proof of Theorem \ref{thm;14.5.2.20}}
\label{subsection;14.5.14.20}

\subsubsection{Quasi-isomorphisms}

We generalize the results in 
\S\ref{subsection;14.5.17.3}
in the context of mixed twistor $\nbigd$-modules.
Let 
$\Omega^1_{\nbigx^{(1)}/\cnum^{2}_{\lambda,\tau}}
(\log \nbigd^{(1)})$
denote the kernel of
$\Omega^1_{\nbigx^{(1)}}(\log \nbigd^{(1)})
\lrarr
 \ptilde_1^{\ast}
 \Omega^1_{\cnum^2_{\lambda,\tau}}$.
We set 
$\Omegatilde^1_{\nbigx^{(1)}/\cnum^{2}_{\lambda,\tau}}
(\log \nbigd^{(1)}):=
\lambda^{-1}
 \Omega^1_{\nbigx^{(1)}/\cnum^{2}_{\lambda,\tau}}
(\log \nbigd^{(1)})$
and 
$\Omegatilde^k_{\nbigx^{(1)}/\cnum^{2}_{\lambda,\tau}}
(\log \nbigd^{(1)}):=
\bigwedge^k\Omegatilde^1_{\nbigx^{(1)}/\cnum^{2}_{\lambda,\tau}}
(\log \nbigd^{(1)})$.

Let $0\leq\alpha<1$.
We set
$\nbigmtilde_{[\alpha P]}(-\nbigd^{(1)}):=
 \nbigmtilde_{[\alpha P]}
 \otimes\nbigo(-\nbigd^{(1)})$.
We have an inclusion of complexes:
\begin{equation}
\label{eq;14.4.23.10}
 \nbigmtilde_{[\alpha P]}(-\nbigd^{(1)})
\otimes
 \Omegatilde^{\bullet}_{\nbigx^{(1)}/\cnum^2_{\lambda,\tau}}
 (\log \nbigd^{(1)})
\lrarr
 U_{\alpha}\nbigmtilde
 \otimes
 \Omegatilde^{\bullet}_{\nbigx^{(1)}/\cnum^2_{\lambda,\tau}}
\end{equation}
The following proposition
is an analogue of Proposition \ref{prop;14.5.1.30},
which we shall prove in \S\ref{subsection;14.12.28.1}.

\begin{proposition}
\label{prop;14.5.2.210}
The morphism {\rm(\ref{eq;14.4.23.10})}
is a quasi-isomorphism.
\end{proposition}

The correspondence $1\longmapsto \upsilon$ induces 
a natural morphism of complexes:
\begin{equation}
\label{eq;14.4.23.11}
 \Omegatilde^{\bullet}_{f,\lambda,\tau}(\alpha)
\lrarr
 \nbigmtilde_{[\alpha P]}(-\nbigd^{(1)})
 \otimes
 \Omegatilde^{\bullet}_{\nbigx^{(1)}/\cnum^2_{\lambda,\tau}}
 (\log \nbigd^{(1)}).
\end{equation}

We obtain the following proposition as 
in the case of Proposition \ref{prop;14.4.24.10},
which we shall prove in \S\ref{subsection;14.12.28.2}.
\begin{proposition}
\label{prop;14.5.2.211}
The morphism {\rm(\ref{eq;14.4.23.11})}
is a quasi-isomorphism.
\end{proposition}

We immediately obtain Theorem \ref{thm;14.5.2.20}
from Proposition \ref{prop;14.5.2.210}
and Proposition \ref{prop;14.5.2.211}.

\subsubsection{Proof of Proposition \ref{prop;14.5.2.210}}
\label{subsection;14.12.28.1}

We essentially repeat the argument in 
\S\ref{subsection;14.5.17.2}.
We have only to check the claim
around any point 
$(\lambda_0,\tau_0,Q)$ of $\nbigd^{(1)}$.
We use a coordinate system
as in \S\ref{subsection;14.4.23.1}.
We set $\nbigd_i^{(1)}:=\{x_i=0\}$
on the coordinate neighbourhood.
Set $\vecp:=[\alpha\veck]$.

Let us consider the case $\tau_0\neq 0$.
For $\ell_1\leq p\leq \ell$,
we set
$S(p):=\{\overbrace{(0,\ldots,0)}^{\ell_1}\}
 \times\seisuu^{p-\ell_1}_{\geq 0}
 \times\{\overbrace{(0,\ldots,0)}^{\ell-p}\}
\subset\seisuu^{\ell_1}\times\seisuu^{p-\ell_1}
 \times\seisuu^{\ell-p}
=\seisuu^{\ell}$.
We consider
\[
 \nbigmtilde^{\leq p}:=
 \sum_{\vecn\in S(p)}
 \deldel^{\vecn}
 \nbigmtilde_{[\alpha P]}.
\]
For any $0\leq \alpha<1$,
it is equal to the following
in $\nbigo_{\nbigx^{(1)}}(\ast\nbigd^{(1)})$
around $(\lambda_0,\tau_0,Q)$:
\[
 \sum_{\vecn\in S(p)}
 \deldel^{\vecn}
 \bigl(
 \nbigo_{\nbigx^{(1)}}(\nbigd^{(1)})
 (\ast\nbigp_{\reduced}^{(1)})
 \upsilon
 \bigr)
=\sum_{\vecn\in S(p)}
 \lambda^{|\vecn|}
 \nbigo_{\nbigx^{(1)}}
 \Bigl(\nbigd^{(1)}
 +\sum n_i\nbigd^{(1)}_i\Bigr)
 (\ast\nbigp_{\reduced}^{(1)})
 \upsilon.
\]
We have 
$\nbigmtilde^{\leq \ell_1}=\nbigmtilde_{[\alpha P]}$.

\begin{lemma}
\label{lem;14.5.15.110}
For $\ell_1+1\leq p\leq \ell$,
the complexes 
$x_p\nbigmtilde^{\leq p-1}
\stackrel{\deldel_p}{\lrarr}
 \nbigmtilde^{\leq p-1}$
and
$\nbigmtilde^{\leq p}
\stackrel{\deldel_p}{\lrarr}
 \nbigmtilde^{\leq p}$
are quasi-isomorphic
with respect to the inclusion.
\end{lemma}
\pf
We set $\nbigd^{(1)}_{<p}:=
 \bigcup_{i=1}^{p-1}
 \{x_i=0\}$
and 
$\nbigd^{(1)}_{\leq p}:=
 \bigcup_{i=1}^{p}
 \{x_i=0\}$
on the neighbourhood.
For $j\in\seisuu$,
we consider the following:
\[
 \lefttop{p}F_j\nbigl^{\leq p}:=
 \nbigo_{\nbigx^{(1)}}\bigl(
 \nbigd^{(1)}+j\nbigd_p^{(1)}
 \bigr)(\ast\nbigd_{< p}^{(1)}).
\]
We set
$\nbigl^{\leq p}:=
\bigcup_{j}\lefttop{p}F_j\nbigl^{\leq p}
=\nbigo_{\nbigx^{(1)}}\bigl(
 \nbigd^{(1)}
 \bigr)(\ast\nbigd_{\leq p}^{(1)})$.
We have 
$\lefttop{p}F_j\nbigl^{\leq p}
\lrarr
 \lefttop{p}F_{j+1}\nbigl^{\leq p}$
induced by $\deldel_p$.
If $j\geq 0$,
the induced morphisms
\[
 \lefttop{p}F_j\nbigl^{\leq p}
 \big/\lefttop{p}F_{j-1}\nbigl^{\leq p}
\lrarr 
\lefttop{p}F_{j+1}\nbigl^{\leq p}
 \big/\lefttop{p}F_{j}\nbigl^{\leq p}
\]
are injective.
Indeed, the induced morphism
$x_p\deldel_{p}$ on
$\lefttop{p}F_j\nbigl^{\leq p}
 \big/\lefttop{p}F_{j-1}\nbigl^{\leq p}$
is the multiplication of
$(j+1)\lambda$.
By using the case $j=0$,
we obtain that the kernel of
$\lefttop{p}F_0
\lrarr
 \lefttop{p}F_1\big/\lefttop{p}F_0$
is $x_p\lefttop{p}F_0$.

For $j\in\seisuu_{\geq 0}$,
we consider the following:
\[
 \lefttop{p}F_j\nbigmtilde^{\leq p}
:=\sum_{\substack{\vecn\in S(p)\\ n_p\leq j}}
 \deldel^{\vecn}
 \nbigmtilde_{[\alpha P]}
=\sum_{\substack{\vecn\in S(p)\\ n_p\leq j}}
 \lambda^{|\vecn|}
 \nbigo_{\nbigx^{(1)}}
 \Bigl(\nbigd^{(1)}
 +\sum n_i\nbigd^{(1)}_i\Bigr)(\ast \nbigp_{\reduced}^{(1)})
 \upsilon
\]
We have
$\lefttop{p}F_0\nbigm^{\leq p}
=\nbigm^{\leq p-1}$.
We have
$\lefttop{p}F_j\nbigmtilde^{\leq p}
\subset 
 \lefttop{p}F_j\nbigl^{\leq p}$
and 
$x_p\lefttop{p}F_0\nbigmtilde^{\leq p}
\subset
 \lefttop{p}F_{-1}\nbigl^{\leq p}$.

Let us consider the following commutative diagram:
\[
 \begin{CD}
 \lefttop{p}F_0\nbigmtilde^{\leq p}\big/
 x_p\lefttop{p}F_0\nbigmtilde^{\leq p}
 @>{a_1}>>
 \lefttop{p}F_{0}\nbigltilde^{\leq p}\big/
 \lefttop{p}F_{-1}\nbigltilde^{\leq p}
 \\
 @V{b_1}VV @V{b_2}VV \\
 \lefttop{p}F_j\nbigmtilde^{\leq p}\big/
 \lefttop{p}F_{j-1}\nbigmtilde^{\leq p}
 @>{a_2}>>
 \lefttop{p}F_{j}\nbigltilde^{\leq p}\big/
 \lefttop{p}F_{j-1}\nbigltilde^{\leq p}.
 \end{CD}
\]
Here, $a_i$ are induced by the inclusions,
and $b_i$ are induced by 
$\deldel_p^j$.
By using an explicit description,
we can check that $a_1$ is a monomorphism.
As mentioned, $b_2$ is a monomorphism.
By the construction, $b_1$ is epimorphism.
Hence, we obtain that 
$b_1$ is an isomorphism,
and that $a_2$ is a monomorphism.

We can deduce the following claims
on the morphisms induced by $\deldel_p$:
\begin{itemize}
\item
If $j\geq 1$,
the induced morphisms
$\lefttop{p}F_j\big/\lefttop{p}F_{j-1}
\stackrel{\deldel_p}{\lrarr}
\lefttop{p}F_{j+1}\big/\lefttop{p}F_{j}$
are isomorphisms.
\item
The kernel of
$\lefttop{p}F_0\nbigmtilde^{\leq p}
=\nbigmtilde^{\leq p-1}
\lrarr
 \lefttop{p}F_1\big/\lefttop{p}F_0$
is $x_p\nbigmtilde^{\leq p-1}$.
\end{itemize}
Then, the claim of Lemma \ref{lem;14.5.15.110}
follows.
\qed

\vspace{.1in}

For $\ell_1\leq p\leq \ell$,
we set
$\nbigh^{(1)}_{>p}:=
 \bigcup_{i=p+1}^{\ell}\{x_i=0\}$.
By using Lemma \ref{lem;14.5.15.110},
the following inclusions of the complexes of sheaves
are quasi-isomorphisms
(see the argument after Lemma \ref{lem;14.5.14.10}):
\[
 \nbigmtilde^{\leq p-1}(-\nbigh^{(1)}_{>p-1})
 \otimes
 \Omegatilde^{\bullet}_{\nbigx^{(1)}/\cnum^2_{\lambda,\tau}}
 (\log \nbigh^{(1)}_{>p-1})
\lrarr
 \nbigmtilde^{\leq p}(-\nbigh^{(1)}_{>p})
 \otimes
 \Omegatilde^{\bullet}_{\nbigx^{(1)}/\cnum^2_{\lambda,\tau}}
 (\log \nbigh^{(1)}_{>p}).
\]
We have 
$\nbigmtilde^{\leq \ell}(-\nbigh^{(1)}_{>\ell})
 \otimes
 \Omegatilde^{\bullet}_{\nbigx^{(1)}/\cnum^2_{\lambda,\tau}}
 (\log \nbigh^{(1)}_{>\ell})
=\nbigmtilde\otimes
 \Omegatilde^{\bullet}_{\nbigx^{(1)}/\cnum^2_{\lambda,\tau}}$,
and 
\begin{multline}
\nbigm^{\leq \ell_1}(-\nbigh^{(1)}_{>\ell_1})
 \otimes
 \Omegatilde^{\bullet}_{\nbigx^{(1)}/\cnum^2_{\lambda,\tau}}
 (\log \nbigh^{(1)}_{>\ell_1})
=\nbigmtilde_{[\alpha P]}(-\nbigh^{(1)})
 \otimes\Omegatilde^{\bullet}
 _{\nbigx^{(1)}/\cnum^2_{\lambda,\tau}}
 (\log \nbigh^{(1)})
 \\
=\nbigmtilde_{[\alpha P]}(-\nbigd^{(1)})
 \otimes\Omegatilde^{\bullet}_{\nbigx^{(1)}/\cnum^2_{\lambda,\tau}}
 (\log \nbigd^{(1)}).
\end{multline}
Hence, we are done in the case $\tau_0\neq 0$.

\vspace{.1in}

Let us consider the case $\tau_0=0$.
For $1\leq p\leq \ell$,
we regard 
 $\seisuu^{p}=\seisuu^{p}\times
 \{\overbrace{(0,\ldots,0)}^{\ell-p}\}
\subset
 \seisuu^{\ell}$.
We set
\[
 U_{\alpha}\nbigmtilde^{\leq p}:=
 \sum_{\vecn\in\seisuu^p_{\geq 0}}
 \deldel^{\vecn}\nbigmtilde_{[\alpha P]},
\quad\quad
 \lefttop{p}F_jU_{\alpha}\nbigmtilde^{\leq p}
 =\sum_{\substack{
 \vecn\in\seisuu^p_{\geq 0}\\
 n_p\leq j }}
 \deldel^{\vecn}\nbigmtilde_{[\alpha P]}.
\]
We have
$U_{\alpha}\nbigmtilde^{\leq \ell}
=U_{\alpha}\nbigmtilde$
and 
$\lefttop{p}F_0U_{\alpha}\nbigmtilde^{\leq p}
=U_{\alpha}\nbigmtilde^{\leq p-1}$.
We consider the following maps:
\[
 \deldel_p:
 \lefttop{p}F_jU_{\alpha}\nbigmtilde^{\leq p}
\lrarr
 \lefttop{p}F_{j+1}U_{\alpha}\nbigmtilde^{\leq p}.
\]

The following lemma is easy to see.
\begin{lemma}
Let $s$ be a section of 
$U_{\alpha}\nbigmtilde$
on $\nbign^{(\lambda_0)}$
with a primitive expression
$s=\sum_{(\vecm,j)\in\nbigs}s_{\vecm,j}$.
Then,
$s$ is a section of
$\lefttop{p}F_jU_{\alpha}\nbigmtilde^{\leq p}$
if and only if
we have 
$m_i\geq 0$ $(i>p)$
and $m_p\geq -j$
for any $\vecm\in\min\pi(\nbigs)$.
\qed
\end{lemma}

\begin{lemma}
\label{lem;14.5.15.120}
If $j\geq 1$,
the following induced morphism 
of $\nbigo_{\nbigx^{(1)}}/x_p\nbigo_{\nbigx^{(1)}}$-modules
is an isomorphism:
\[
 \frac{\lefttop{p}F_jU_{\alpha}\nbigmtilde^{\leq p}}
 {\lefttop{p}F_{j-1}U_{\alpha}\nbigmtilde^{\leq p}}
\stackrel{\deldel_p}{\lrarr}
 \frac{\lefttop{p}F_{j+1}U_{\alpha}\nbigmtilde^{\leq p}}
 {\lefttop{p}F_{j}U_{\alpha}\nbigmtilde^{\leq p}}
\]
\end{lemma}
\pf
It is surjective by construction.
Let $s$ be a non-zero section of
$\lefttop{p}F_jU_{\alpha}\nbigmtilde^{\leq p}$
on $\nbign^{(\lambda_0)}$
with a primitive decomposition 
$s=\sum_{(\vecm,j)\in\nbigs}s_{\vecm,j}$
such that
$\deldel_ps$ is also a section of 
$\lefttop{p}F_jU_{\alpha}\nbigmtilde^{\leq p}$.
We set
$s':=\sum_{m_p=-j}s_{\vecm,j}$
and 
$s'':=\sum_{m_p>-j}s_{\vecm,j}$.
Because
$\deldel_ps''\in \lefttop{p}F_{j}U_{\alpha}\nbigmtilde^{\leq p}$,
we obtain 
$\deldel_ps'\in\lefttop{p}F_jU_{\alpha}\nbigmtilde^{\leq p}$.
If $s'$ is non-zero,
$\deldel_ps'=\sum_{m_p=-j} \deldel_ps_{\vecm,j}$
is a primitive expression of $\deldel_ps'$.
We obtain
$\deldel_ps'\not\in \lefttop{p}F_{j}U_{\alpha}\nbigmtilde^{\leq p}$,
and thus we arrive at a contradiction.
Hence, $s'=0$,
i.e.,
$s\in \lefttop{p}F_{j-1}U_{\alpha}\nbigmtilde^{\leq p}$.
\qed

\begin{lemma}
\label{lem;14.5.15.121}
The kernel of the following induced morphism
is $x_pU_{\alpha}\nbigmtilde^{\leq p-1}$:
\[
 U_{\alpha}\nbigmtilde^{\leq p-1}
\stackrel{\deldel_p}{\lrarr}
\frac{\lefttop{p}F_1U_{\alpha}\nbigmtilde^{\leq p}}
 {U_{\alpha}\nbigmtilde^{\leq p-1}}
\]
\end{lemma}
\pf
Let $s$ be a section of 
$U_{\alpha}\nbigmtilde^{\leq p-1}$,
and we take a primitive expression
$s=\sum_{(\vecm,j)\in\nbigs}s_{\vecm,j}$
such that 
$\delbar_ps\in U_{\alpha}\nbigmtilde^{p-1}$.
We set
$s':=\sum_{m_p=0}s_{\vecm,j}$
and 
$s'':=\sum_{m_p>0}s_{\vecm,j}$.
Because
$s''\in x_pU_{\alpha}\nbigmtilde^{\leq p-1}$,
we have
$\deldel_ps''\in U_{\alpha}\nbigmtilde^{\leq p-1}$
as in the proof of Lemma \ref{lem;14.5.1.101}.
Hence, we have
$\deldel_ps'\in U_{\alpha}\nbigmtilde^{\leq p-1}$.
If $s'\neq 0$,
$\deldel_ps'=\sum_{m_p=0}\deldel_ps_{\vecm,j}$
is a primitive expression of $s'$.
We obtain
$\deldel_ps'\not\in U_{\alpha}\nbigmtilde^{\leq p-1}$,
and we arrive at a contradiction.
Hence, we have $s'=0$,
i.e.,
$s\in x_p U_{\alpha}\nbigm^{\leq p-1}$.
\qed

\vspace{.1in}

By Lemma \ref{lem;14.5.15.120}
and  Lemma \ref{lem;14.5.15.121},
the inclusions of the complexes
\[
 \bigl(
 x_pU_{\alpha}\nbigmtilde^{\leq p-1}
\stackrel{\deldel_p}{\lrarr}
 U_{\alpha}\nbigmtilde^{\leq p-1}
 \bigr)
\lrarr
 \bigl(
 U_{\alpha}\nbigmtilde^{\leq p}
\stackrel{\deldel_p}{\lrarr}
 U_{\alpha}\nbigmtilde^{\leq p}
 \bigr)
\]
are quasi-isomorphisms.
For $0\leq p \leq \ell$,
we set $\nbigd^{(1)}_{>p}:=\bigcup_{i=p+1}^{\ell}\{x_i=0\}$
on the neighbourhood.
We obtain that the following inclusions of the complexes of sheaves
are quasi-isomorphisms:
\begin{equation}
 U_{\alpha}\nbigmtilde^{\leq p-1}(-\nbigd^{(1)}_{>p-1})
 \otimes
 \Omegatilde_{\nbigx^{(1)}/\cnum^2_{\lambda,\tau}}^{\bullet}
 (\log \nbigd^{(1)}_{>p-1})
\lrarr 
  U_{\alpha}\nbigmtilde^{\leq p}(-\nbigd^{(1)}_{>p})
 \otimes
 \Omegatilde_{\nbigx^{(1)}/\cnum^2_{\lambda,\tau}}^{\bullet}
 (\log \nbigd^{(1)}_{>p})
\end{equation}
We have 
$U_{\alpha}\nbigmtilde^{\leq \ell}(-\nbigd^{(1)}_{>\ell})
 \otimes
 \Omegatilde_{\nbigx^{(1)}/\cnum^2_{\lambda,\tau}}^{\bullet}
 (\log \nbigd^{(1)}_{>\ell})
=U_{\alpha}\nbigmtilde\otimes
 \Omegatilde_{\nbigx^{(1)}/\cnum^2_{\lambda,\tau}}^{\bullet}$.
We also have 
\[
 U_{\alpha}\nbigmtilde^{\leq 0}(-\nbigd^{(1)}_{>0})
 \otimes
 \Omegatilde_{\nbigx^{(1)}/\cnum^2_{\lambda,\tau}}^{\bullet}
 (\log \nbigd^{(1)}_{>0})
=\nbigmtilde_{[\alpha P]}(-\nbigd^{(1)})
 \otimes
 \Omegatilde_{\nbigx^{(1)}/\cnum^2_{\lambda,\tau}}^{\bullet}
 (\log \nbigd^{(1)}).
\]
Hence, Proposition \ref{prop;14.5.2.210} is proved.
\qed

\subsubsection{Proof of Proposition \ref{prop;14.5.2.211}}
\label{subsection;14.12.28.2}

We essentially repeat the argument in \S\ref{subsection;14.5.17.1}.
We have only to check the claim
around any point of $\nbigp_{\reduced}^{(1)}$.
We use the coordinate system as in \S\ref{subsection;14.4.23.1}.
We set $\vecp:=[\alpha\veck]$.
For any non-negative integer $N$,
we set
\[
 G_N\bigl(
 \nbigmtilde_{[\alpha P]}(-\nbigd^{(1)})
 \bigr)
:=\sum_{j=0}^N
 \nbigo_{\nbigx^{(1)}}x^{-\vecp}
 (\tau f)^j\upsilon.
\]
We define
$G_N\bigl(
 \nbigmtilde_{[\alpha P]}(-\nbigd^{(1)})
 \otimes\Omegatilde^k_{\nbigx^{(1)}/\cnum^2_{\lambda,\tau}}
 (\log \nbigd^{(1)})
 \bigr)$
as 
\[
 G_N\bigl(
 \nbigmtilde_{[\alpha P]}(-\nbigd^{(1)})
 \bigr)
 \otimes
 \Omegatilde^k_{\nbigx^{(1)}/
 \cnum^2_{\lambda,\tau}}
 (\log \nbigd^{(1)}).
\]
We set
$G_{-1}\bigl(
 \nbigmtilde_{[\alpha P]}(-\nbigd^{(1)})
\otimes
 \Omegatilde^k_{\nbigx^{(1)}/\cnum^2_{\lambda,\tau}}(\log \nbigd^{(1)})
 \bigr)
:=\Omegatilde_{f,\lambda,\tau}^{k}(\alpha)\upsilon
=\Omegatilde_{f,\lambda,\tau}^kx^{-\vecp}\upsilon$,
where 
$\Omegatilde_{f,\lambda,\tau}^k:=
 \Omegatilde^k_{f,\lambda,\tau}(0)$.
Let $N\geq 0$.
Take a section
\[
 \omega
 =\sum_{j=0}^N
 \omega_j
 x^{-\vecp}(\tau f)^j\cdot \upsilon
 \in 
 G_N\bigl(
 \nbigmtilde_{[\alpha P]}(-\nbigd^{(1)})\otimes
 \Omegatilde^k
 _{\nbigx^{(1)}/\cnum^2_{\lambda,\tau}}(\log\nbigd^{(1)})
 \bigr),
\]
where
$\omega_j\in
 \Omegatilde^k_{\nbigx^{(1)}/\cnum^2_{\lambda,\tau}}(\log \nbigd^{(1)})$.
If
$d\omega
 \in
 G_{N}\bigl(
  \nbigmtilde_{[\alpha P]}(-\nbigd^{(1)})
 \otimes\Omegatilde^{k+1}
 _{\nbigx^{(1)}/\cnum^2_{\lambda,\tau}}(\log\nbigd^{(1)})
 \bigr)$,
then we have
\[
 \lambda^{-1}\tau df\wedge \omega_N\,(\tau f)^Nx^{-\vecp}\upsilon
 \in
 G_N\bigl(
 \nbigmtilde_{[\alpha P]}(-\nbigd^{(1)})
\otimes
 \Omegatilde^{k+1}_{\nbigx^{(1)}/\cnum^2_{\lambda,\tau}}(\log
 \nbigd^{(1)})
 \bigr).
\]
As in Lemma \ref{lem;14.5.21.1},
we obtain
$df\wedge\omega_N
 \in
 \Omegatilde^{k+1}
 _{\nbigx^{(1)}/\cnum^2_{\lambda,\tau}}(\log \nbigd^{(1)})$,
i.e.,
$\omega_N$ is a section of
$\Omegatilde^k_{f,\lambda,\tau}$.
\begin{lemma}
We have an expression
\[
 \omega_N=
 (df/f)\wedge \kappa_1
+f^{-1}\kappa_2,
\]
where
$\kappa_1$ and $\kappa_2$
are sections of
$\Omegatilde^{k-1}_{\nbigx^{(1)}/\cnum^2_{\lambda,\tau}}
 (\log\nbigd^{(1)})$
and 
$\Omegatilde^{k}_{\nbigx^{(1)}/\cnum^2_{\lambda,\tau}}
 (\log \nbigd^{(1)})$,
respectively.
\qed
\end{lemma}

If $N\geq 1$,
we have 
$f^{-1}\kappa_2(\tau f)^N
=\tau \kappa_2(\tau f)^{N-1}$.
Hence, we have
\[
 \omega-d\bigl(\kappa_1(\tau f)^{N-1}x^{-\vecp}\upsilon\bigr)
\in
 G_{N-1}\bigl(
 \nbigmtilde_{[\alpha P]}(-\nbigd^{(1)})
 \otimes\Omegatilde^k_{\nbigx^{(1)}/\cnum^2_{\lambda,\tau}}
 (\log \nbigd^{(1)})
 \bigr).
\]
We also have
$f^{-1}\kappa_1(\tau f)^Nx^{-\vecp}\upsilon
 \in G_{N-1}$.

Let $\omega$ be a local section of
$G_N\bigl(
 \nbigmtilde_{[\alpha P]}(-\nbigd^{(1)})
 \otimes\Omegatilde^k_{\nbigx^{(1)}/\cnum^2_{\lambda,\tau}}
 (\log \nbigd^{(1)})
 \bigr)$
such that 
$d\omega$
is a local section of
$G_{-1}\bigl(
 \nbigmtilde_{[\alpha P]}(-\nbigd^{(1)})
 \otimes\Omegatilde^k_{\nbigx^{(1)}/\cnum^2_{\lambda,\tau}}
 (\log \nbigd^{(1)})
 \bigr)$.
Then, by applying the previous argument successively,
we can find 
a local section
$\tau$ of 
$G_{N-1}\bigl(
 \nbigmtilde_{[\alpha P]}(-\nbigd^{(1)})
 \otimes\Omegatilde^{k-1}_{\nbigx^{(1)}/\cnum^2_{\lambda,\tau}}
 (\log \nbigd^{(1)})
 \bigr)$
such that
$\omega-d\tau$
is a local section of 
$G_{-1}\bigl(
 \nbigmtilde_{[\alpha P]}(-\nbigd^{(1)})
 \otimes\Omegatilde^{k}_{\nbigx^{(1)}/\cnum^2_{\lambda,\tau}}(\log\nbigd^{(1)})
 \bigr)$.

\begin{itemize}
\item
If a local section $\omega$
of $G_N\bigl(
 \nbigmtilde_{[\alpha P]}(-\nbigd^{(1)})
 \otimes\Omegatilde^k_{\nbigx^{(1)}/\cnum^2_{\lambda,\tau}}
 (\log \nbigd^{(1)})
 \bigr)$
satisfies $d\omega=0$,
we can find 
a local section
$\tau$ of 
$G_{N-1}\bigl(
 \nbigmtilde_{[\alpha P]}(-\nbigd^{(1)})
 \otimes\Omegatilde^{k-1}_{\nbigx^{(1)}/\cnum^2_{\lambda,\tau}}
 (\log \nbigd^{(1)})
 \bigr)$
such that
$\omega-d\tau$
is a local section of 
$G_{-1}\bigl(
 \nbigmtilde_{[\alpha P]}(-\nbigd^{(1)})
 \otimes\Omegatilde^{k}_{\nbigx^{(1)}/\cnum^2_{\lambda,\tau}}(\log\nbigd^{(1)})
 \bigr)$.
\item
Let $\omega$ be a local section of
$G_{-1}\bigl(
 \nbigmtilde_{[\alpha P]}(-\nbigd^{(1)})
 \otimes\Omegatilde^{k}_{\nbigx^{(1)}/\cnum^2_{\lambda,\tau}}
 (\log \nbigd^{(1)})
 \bigr)$
such that $d\omega=0$.
Suppose that we have
a local section $\tau$ of
$G_{N}\bigl(
 \nbigmtilde_{[\alpha P]}(-\nbigd^{(1)})
 \otimes
 \Omegatilde^{k-1}_{\nbigx^{(1)}/\cnum^2_{\lambda,\tau}}
 (\log\nbigd^{(1)})
 \bigr)$
such that
$\omega=d\tau$.
Then, we can find 
a local section $\sigma$
of 
$G_{N-1}\bigl(
 \nbigmtilde_{[\alpha P]}(-\nbigd^{(1)})
 \otimes\Omegatilde^{k-2}_{\nbigx^{(1)}/\cnum^2_{\lambda,\tau}}
 (\log\nbigd^{(1)})
 \bigr)$
such that 
$\tau-d\sigma$ is a local section of
$G_{-1}\bigl(
 \nbigmtilde_{[\alpha P]}(-\nbigd^{(1)})
 \otimes
 \Omegatilde^{k-1}
 _{\nbigx^{(1)}/\cnum^2_{\lambda,\tau}}(\log\nbigd^{(1)})
 \bigr)$.
We have
$\omega=d(\tau-d\sigma)$.
\end{itemize}
Then, we obtain the claim of Proposition \ref{prop;14.5.2.211}.
The proof of Theorem \ref{thm;14.5.2.20}
is also completed.
\qed

\vspace{.1in}
\noindent
Address: Research Institute for Mathematical Sciences,
Kyoto University, Kyoto 606-8502, Japan\\
Email: takuro@kurims.kyoto-u.ac.jp

\end{document}